\documentclass{article}
\usepackage{geometry}
\usepackage{amssymb}
\usepackage{amsmath}
\usepackage{amsthm}
\usepackage{amsrefs}
\usepackage{enumitem}
\usepackage{caption}
\usepackage{subcaption}
\usepackage{float}
\usepackage{tikz}
\usepackage{etoolbox, calc}
\usetikzlibrary{patterns}
\usetikzlibrary{arrows}
\usetikzlibrary{calc}

\newtheorem{theorem}{Theorem}[section]
\newtheorem{corollary}[theorem]{Corollary}
\newtheorem{lemma}[theorem]{Lemma}
\newtheorem{proposition}[theorem]{Proposition}
\newtheorem{conjecture}[theorem]{Conjecture}
\newtheorem{problem}[theorem]{Problem}

\theoremstyle{definition}
\newtheorem{definition}[theorem]{Definition}
\newenvironment{example1}[1][]{\refstepcounter{theorem}\par\medskip
   \noindent \textbf{Example~\thetheorem. #1} \rmfamily}{\leavevmode\unskip\penalty9999 \hbox{}\nobreak\hfill\quad\hbox{$\Diamond$}\medskip}

\theoremstyle{remark}

\newtheorem{remark}[theorem]{Remark}

\numberwithin{equation}{section}

\newcommand{\R}{\mathbf{R}}

\newcommand{\Z}{\mathbf{Z}}

\newcommand{\m}{m_G}
\newcommand{\cb}{\mathcal{C}_B}
\newcommand{\cw}{\mathcal{C}_W}
\newcommand{\mc}{\mathcal{C}}
\newcommand{\mat}{\mathcal{M}}
\newcommand{\bp}{\mathcal{B}}

\DeclareMathOperator{\nul}{null}

\DeclareMathOperator{\ch}{ch}
\DeclareMathOperator{\chdeg}{chdeg}
\DeclareMathOperator{\lcm}{lcm}

\newcounter{colsum}

\newcommand\vertices[4]{%
	\foreach \x in {#1,...,#3}
	\foreach \y in {#2,...,#4}
	{
		\setcounter{colsum}{\x+\y};
		
		\ifnumodd{\value{colsum}}{\def\altcol{white}}{\def\altcol{black}}
		\draw[fill=\altcol] (\x,\y) circle (.1);
	}
}

\newcommand{\edge}[2]{\overline{#1 #2}}

\title{Channels, Billiards, and Perfect Matching 2-Divisibility}
\author{Grant T. Barkley and Ricky Ini Liu}
\date{}
\begin{document}

\maketitle

\begin{abstract}
Let $m_G$ denote the number of perfect matchings of the graph $G$. 
We introduce a number of combinatorial tools for determining the parity of $m_G$ and giving a lower bound on the power of 2 dividing $m_G$. In particular, we introduce certain vertex sets called channels, which correspond to elements in the kernel of the adjacency matrix of $G$ modulo $2$. A result of Lov\'asz states that the existence of a nontrivial channel is equivalent to $m_G$ being even. We give a new combinatorial proof of this result and strengthen it by showing that the number of channels gives a lower bound on the power of $2$ dividing $m_G$ when $G$ is planar. 
We describe a number of local graph operations which preserve the number of channels. 
We also establish a surprising connection between 2-divisibility of $m_G$ and dynamical systems by showing an equivalency between channels and billiard paths. 
We exploit this relationship to show that $2^{\frac{\gcd(m+1,n+1)-1}{2}}$ divides the number of domino tilings of the $m\times n$ rectangle. We also use billiard paths to give a fast algorithm for counting channels (and hence determining the parity of the number of domino tilings) in simply connected regions of the square grid.
\end{abstract}

\section{Introduction}

Given a graph $G$, a \emph{perfect matching} of $G$ is a set of edges $\mu$ such that each vertex of $G$ is contained in a unique edge in $\mu$. We let $\m$ denote the number of distinct perfect matchings of $G$. The problem of determining $\m$ arises in various mathematical contexts, particularly in tiling problems, but also in statistical mechanics \cite{Kasteleyn1961}, spectral graph theory \cite{Kenyon1999}, network analysis \cite{Eppstein2018}, total positivity \cite{Lam2015}, and representation theory \cite{Fraser2017}.
Exact formulas for $\m$ over an infinite family of graphs are quite rare. One notable exact formula is for $G=\mathcal R_{m\times n}$, the rectangular subgraph of the square lattice with $m$ rows of $n$ vertices. In this case, a famous result of Kasteleyn \cite{Kasteleyn1961} gives 
\[\m^4 = \prod_{j=1}^{n}\prod_{k=1}^{m} \left(4\cos^2\frac{j\pi }{n+1} + 4\cos^2\frac{k\pi}{m+1}\right). \]

From this product we may extract certain number-theoretic information: for example, $\m$ is always divisible by $2^{\frac{\gcd(n+1,m+1)-1}{2}}$ \cite{Pachter1997}. Studying similar 2-divisibility patterns is a common theme in the literature on domino tilings, which are equivalent to perfect matchings of subgraphs of the square lattice (see, e.g., \cites{Ciucu1996, Cohn1999,
Pachter1997, Propp1999, Tenner2005, Tenner2009}). It is often the case that the 2-component of the prime factorization of $\m$ follows a predictable pattern, even when an exact formula for $\m$ is elusive or unwieldy. In Propp's perfect matching problem anthology \cite{Propp1999}, he gives a number of conjectured and known power of 2 patterns for various graphs. For example, the following is a refinement by Pachter \cite{Pachter1997} of one of these conjectures.

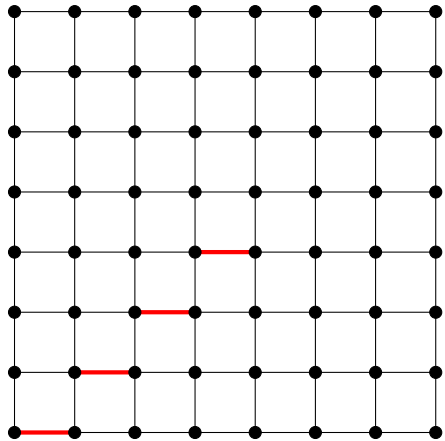
\begin{figure}[htb]                                       
	\centering
	\begin{tikzpicture}[scale=0.8]
	\draw (0,0) grid (7,7);
	\pgfmathsetmacro{\radius}{0.1}
	\foreach \x in {0, 1, 2, 3} {
		\draw[ultra thick, red] (\x,\x) -- ++(1,0);
	}
	
	\foreach \x in {0, 1,..., 7} {
		\foreach \y in {0, 1,..., 7} {
			\draw[fill=black] (\x,\y) circle (\radius);
		}
	}
	\end{tikzpicture}
	\caption{The step-diagonal for $r=4$.}
	\label{fig:first}
\end{figure}

\begin{conjecture}[Deleting from step-diagonals] \label{step-diag-conj}
	Let $\mathcal R_{2r\times 2r}$ be the $2r\times 2r$ grid graph shown in Figure \ref{fig:first}. If $G$ is a subgraph with the same vertex set, constructed by deleting any $k$ of the edges highlighted in the figure, then 
	$$\m = 2^{r-k}b$$
	for some odd integer $b$.
\end{conjecture}

The following special case of a theorem of Ciucu is perhaps the most widely used 2-divisibility result in the literature.
\begin{proposition} [Ciucu's Factorization Theorem \cite{Ciucu1996}]
	\label{ciucu}
	If a bipartite planar graph has a line of symmetry containing $2r$ vertices, and no edges lie on the line or connect two vertices on opposite sides of the line, then the number of perfect matchings of the graph is divisible by $2^r$.
\end{proposition}

However, the symmetry requirement in Ciucu's theorem means that it does not apply to graphs such as those described in Conjecture \ref{step-diag-conj} or to $\mathcal R_{m\times n}$ for $m\neq n$ (though it is quite important for studying the square $\mathcal R_{m\times m}$). The results we introduce here provide a uniform (partial) explanation of power of 2 patterns in terms of the geometry of the graph. Our foundational construction is based on the following result known to Lov\'asz. 

\begin{proposition}[\cite{Lovasz1993}, Problem 5.18]\label{even}
	Let $G$ be a graph. Then $\m$ is even if and only if there is a nonempty vertex set $C\subseteq V$ such that every vertex in $G$ is adjacent to an even number of vertices in $C$.
\end{proposition}

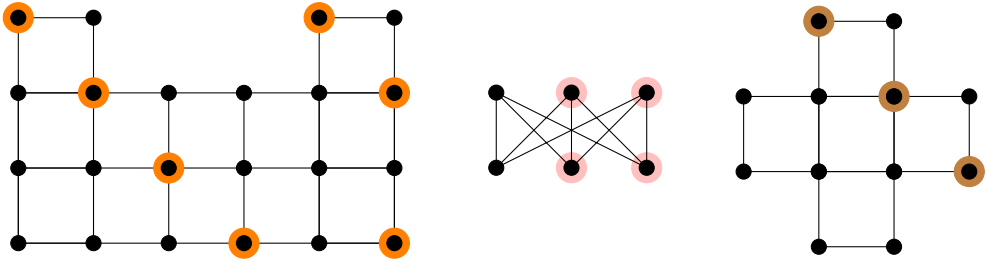
\begin{figure}
    \begin{subfigure}{.48\textwidth}
    \centering
    \begin{tikzpicture}
        \draw (0,0) grid (1, 3);
        \draw (0,0) grid (5, 2);
        \draw (4,0) grid (5,3);
        \foreach \x in {0, 1,..., 3} {
		    \draw[orange,fill] (\x,3-\x) circle(.2);
	    }
        \draw[orange,fill] (5,0) circle (.2);
        \draw[orange,fill] (5,2) circle (.2);
        \draw[orange,fill] (4,3) circle (.2);
        \foreach \x in {0, 1,..., 5} {
		\foreach \y in {0, 1, 2} {
			\draw[fill=black] (\x,\y) circle (.1);
		}}
        \draw[fill=black] (0,3) circle (.1) (1,3) circle (.1);
        \draw[fill=black] (4,3) circle (.1) (5,3) circle (.1); 
    \end{tikzpicture}
\end{subfigure}
\begin{subfigure}{.21\textwidth}
    \begin{tikzpicture}
        \draw[pink,fill] (1,0) circle (.2);
        \draw[pink,fill] (2,0) circle (.2);
        \draw[pink,fill] (1,1) circle (.2);
        \draw[pink,fill] (2,1) circle (.2);
        \draw[fill=black] (0,0) circle (.1) (1,0) circle (.1) (2,0) circle (.1);
        \draw[fill=black] (0,1) circle (.1) (1,1) circle (.1) (2,1) circle (.1);
        \draw (0,0) -- (0,1) -- (1,0) -- (1,1) -- (2,0) -- (2,1) -- (0,0);
        \draw (0,0) -- (1,1)    (0,1) -- (2,0)    (1,0) -- (2,1);
    \end{tikzpicture}
\end{subfigure}
\begin{subfigure}{.25\textwidth}
    \begin{tikzpicture}
        \draw (0,-1) grid (1,2);
		\draw (-1,0) grid (2,1);
        \draw[brown,fill] (0,2) circle (.2);
        \draw[brown,fill] (1,1) circle (.2);
        \draw[brown,fill] (2,0) circle (.2);
        \foreach \x in {-1, 0,..., 2} {
		\foreach \y in {0, 1} {
			\draw[fill=black] (\x,\y) circle (.1);
			\draw[fill=black] (\y,\x) circle (.1);
		}}
    \end{tikzpicture}
\end{subfigure}
    \caption{Three graphs, each with a channel $C$ depicted by colored rings around the vertices contained in $C$.}\label{fig:exchannel}
\end{figure}

We will call a vertex set $C$ satisfying the hypothesis of Proposition~\ref{even} a \emph{channel}. We also count the empty set as a trivial channel. Figure \ref{fig:exchannel} shows examples of channels in different graphs. (The name ``channel'' is intended to evoke the image of a river winding its way around $G$, which may be redirected by ``digging'' edges and vertices out of $G$.) The symmetric difference of two channels is again a channel, so they form a vector space over $\Z/2\Z$. This implies that the number of distinct channels in a graph is always a power of 2. The space of channels can be identified with the kernel of the adjacency matrix of $G$ modulo 2 (see Lemma \ref{ker2chan}). 

Lov\'asz's result already shows the importance of channels for determining the parity of $\m$. The main theorem of this paper shows that channels have even more to say for planar graphs.

\begin{theorem}[Channeling 2s]\label{2ch}
    Let $G$ be a planar graph. Then the number of distinct channels in $G$ divides $\m^2$.
\end{theorem}

Since the number of channels will always be a power of 2, Theorem \ref{2ch} gives a lower bound on the power of 2 dividing $\m$ for any planar graph. We will prove this theorem in greater generality (for any graph admitting a Kasteleyn signing) in Theorem \ref{twos}. We demonstrate the strength of this theorem in a number of examples throughout the article. In particular, we show 2-divisibility results related to graphs described above: that $2^{\frac{\gcd(n+1,m+1)-1}{2}}$ divides $m_{\mathcal R_{m\times n}}$ and that $2^{r-k}$ divides $\m$ for the graph described in Conjecture \ref{step-diag-conj}. These are Corollary \ref{corarithmetic} and Proposition \ref{step-diag}, respectively.

Because of their utility, the majority of this paper is dedicated to studying the structure of channels and methods for finding them. Many of our results are tailored for subgraphs of the square lattice, where perfect matchings are equivalent to domino tilings of a region. (When possible, however, we will state results in greater generality.) Our most fascinating result is a characterization of channels in terms of dynamical systems. We state the result here for subgraphs of the square lattice and show the general case in Section \ref{secbounce}.

\begin{figure}[H]
    \centering
    \begin{tikzpicture}
        \draw (0,0) grid (3,4);
        \draw (3,1) grid (4,4);
        \begin{scope}[shift={(-1,0)}]
        \vertices{1}{0}{4}{4};
        \vertices{5}{1}{5}{4};
        \end{scope}
    \end{tikzpicture}
    \caption{A graph $G$ which is an induced subgraph of the square lattice bounded by a simple cycle.}
    \label{unitsqgraph}
\end{figure}
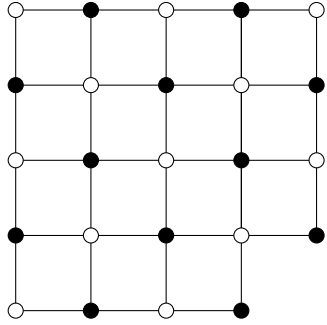

Consider any simple cycle in the square lattice (considered as a grid graph). This cycle divides the square lattice into an interior and an (unbounded) exterior. Let $G$ be the subgraph of the square lattice consisting of all vertices and edges in the cycle or its interior. 
In the dual language of domino tilings, graphs constructed in this way correspond roughly to simply connected subsets of a square grid. Since $G$ is bipartite, we 2-color the vertices of $G$ black and white. An example of such a graph is shown in Figure \ref{unitsqgraph}.

Now we define a \emph{billiard nest} on $G$ to be any collection of paths traced out by billiard balls placed on black vertices of $G$ and launched at 45 degree angles. When a billiard ball reaches a wall, it reflects at a 90 degree angle and proceeds in its new direction, continuing until it is caught by a corner or returns to its start position. (If the billiard brushes past a corner or hits one head-on, the situation is more complicated---the path splits into two paths continuing in different directions. See Section \ref{secbounce} for more examples and a precise definition.) 

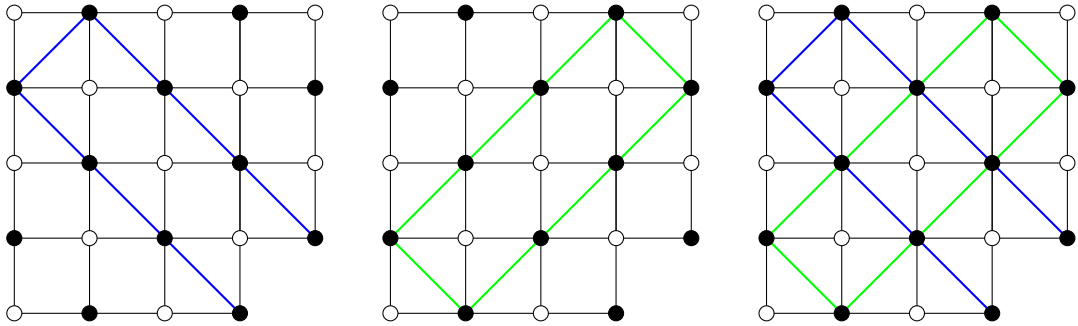
\begin{figure}[H]
    \centering
    \begin{tikzpicture}
        \draw (0,0) grid (3,4);
        \draw (3,1) grid (4,4);
        \draw[blue, thick] (3,0) -- (0,3) -- (1,4) -- (4,1);
        \begin{scope}[shift={(-1,0)}]
        \vertices{1}{0}{4}{4};
        \vertices{5}{1}{5}{4};
        \end{scope}
        
        \begin{scope}[shift={(5,0)}]
        \draw (0,0) grid (3,4);
        \draw (3,1) grid (4,4);
        \draw[green, thick] (1,0) -- (0,1) -- (3,4) -- (4,3) -- (1,0);
        \begin{scope}[shift={(-1,0)}]
        \vertices{1}{0}{4}{4};
        \vertices{5}{1}{5}{4};
        \end{scope}
        \end{scope}
        
        \begin{scope}[shift={(10,0)}]
        \draw (0,0) grid (3,4);
        \draw (3,1) grid (4,4);
        \draw[blue, thick] (3,0) -- (0,3) -- (1,4) -- (4,1);
        \draw[green, thick] (1,0) -- (0,1) -- (3,4) -- (4,3) -- (1,0);
        \begin{scope}[shift={(-1,0)}]
        \vertices{1}{0}{4}{4};
        \vertices{5}{1}{5}{4};
        \end{scope}
        \end{scope}
    \end{tikzpicture}
    \caption{The three nonempty billiard nests in $G$.}
\end{figure}

Remarkably, channels and billiard nests are intrinsically connected. Let $G'$ be the \emph{inner subgraph} of $G$, the subgraph formed by removing all vertices of $G$ which are incident to the unbounded face and all edges incident to those vertices. 

\begin{theorem}\label{billintro}
Let $G$ satisfy the assumptions described above, and let $G'$ be the inner subgraph of $G$. Then the number of billiard nests in $G$ is twice the number of channels in $G'$ that contain only black vertices from $G'$.

As a result, the number of billiard nests in $G$ divides $2m_{G'}$. Moreover,
if $G$ has an equal number of black and white vertices, then
$G$ has exactly one nonempty billiard nest if and only if $m_{G'}$ is odd.
\end{theorem}

\begin{figure}[H]
    \centering
    \begin{tikzpicture}
        \begin{scope}[opacity=.1]
        \draw (0,0) grid (3,4);
        \draw (3,1) grid (4,4);
        \begin{scope}[shift={(-1,0)}]
        \vertices{1}{0}{4}{4};
        \vertices{5}{1}{5}{4};
        \end{scope}
        \end{scope}
        
        \draw (1,1) grid (2,3);
        \draw (3,2) grid (2,3);
        \begin{scope}[shift={(-1,0)}]
        \vertices{2}{1}{3}{3};
        \vertices{3}{2}{4}{3};
        \end{scope}
    \end{tikzpicture}
    \caption{The inner subgraph $G'$ of the graph $G$ defined in Figure \ref{unitsqgraph}.}
    \label{innerunsq}
\end{figure}
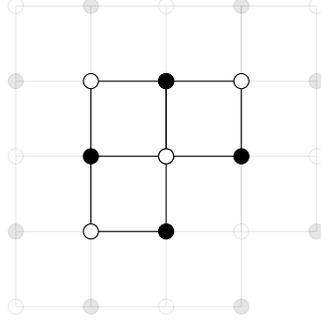

For example, for the graph $G$ in Figure \ref{unitsqgraph}, the inner subgraph $G'$ is shown in Figure \ref{innerunsq}. Since $G$ has 4 billiard nests and satisfies the hypothesis of Theorem \ref{billintro}, it follows that $m_{G'}$ is divisible by 2. In this case we may simply count the perfect matchings of $G'$; we find that $m_{G'}=4$, which is divisible by $2$ as anticipated by the theorem. 
In Theorem \ref{bounce} and its corollary, we prove this billiard--channel correspondence in much greater generality (for \emph{inner semi-Eulerian} graphs).

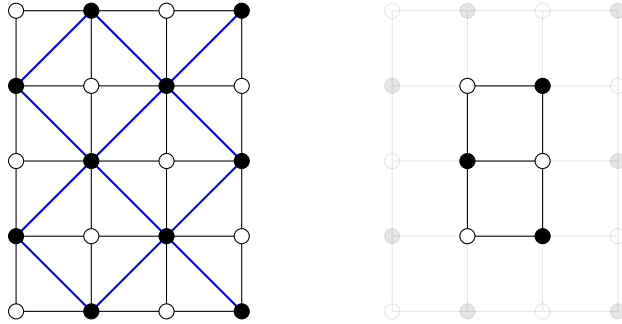
\begin{figure}[H]
    \centering
    \begin{tikzpicture}
        \draw (0,0) grid (3,4);
        \draw[blue, thick] (3,0) -- (0,3) -- (1,4) -- (3,2) -- (1,0) -- (0,1) -- (3,4);
        \begin{scope}[shift={(-1,0)}]
        \vertices{1}{0}{4}{4};
        \end{scope}
        
        \begin{scope}[shift={(5,0)}]
        \begin{scope}[opacity=.1]
        \draw (0,0) grid (3,4);
        \begin{scope}[shift={(-1,0)}]
        \vertices{1}{0}{4}{4};
        \end{scope}
        \end{scope}
        \draw (1,1) grid (2,3);
        \begin{scope}[shift={(-1,0)}]
        \vertices{2}{1}{3}{3};
        \end{scope}
        \end{scope}
    \end{tikzpicture}
    \caption{A graph $G$ which has only one nonempty billiard nest, along with its inner subgraph $G'$.}
    \label{graphno2}
\end{figure}

The connection between 2-divisibility, channels, and dynamical systems explains both the sensitivity and the regularity of perfect matching 2-divisibility. Small changes to $G$ can result in entirely different billiard dynamics, with the effects visible in $m_{G'}$. For instance, if we take $G$ to be the graph in Figure \ref{graphno2}, then there is only one nonempty billiard nest.

Since (counting the empty nest) $G$ has 2 billiard nests, by Theorem \ref{billintro} the inner subgraph $G'$ has no nonempty channels on black vertices. Since $G'$ has the same number of black and white vertices, this will imply (by Lemma \ref{sqbipart}) that $G'$ has no nonempty channels at all. Thus by Proposition \ref{even}, $G'$ has an odd number of perfect matchings---in this case $3$. 

The dynamics involved can also induce a regularity in the 2-divisibility of $\m$. The well-known theory of arithmetic billiards describes billiard nests for rectangles in terms of divisibility properties of the rectangle side lengths. In Section \ref{arithmetic}, we exploit these results to explain the factor of $2^{\frac{\gcd(m+1,n+1)-1}{2}}$ dividing $\m$ for the $m\times n$ grid graph.

Billiard nests give a global explanation of channel structure for many graphs. Sometimes we are instead interested in local behavior. For instance, we may have a family of graphs that are globally similar but that differ locally. To relate these graphs, we introduce a set of channel-preserving graph operations and show that they may be applied repeatedly to reduce many graphs to a graph with no edges, whose channels are apparent.

The paper is organized as follows. The proof of Theorem \ref{2ch} is algebraic and given in Sections 2 and 3. Sections 4, 5, and 6 are independent and may be read in any order. Section 4 describes billiard nests for a large class of graphs called \emph{inner semi-Eulerian graphs}. The results described in the introduction are applied to the rectangle grid graph, connecting its 2-divisibility to the theory of arithmetic billiards. Section 4 concludes with a fast algorithm that constructs the billiard nests and therefore the channels for certain graphs. In Section 5, we give a combinatorial proof of Proposition \ref{even} for bipartite graphs. In the course of this proof we introduce a graph move called \emph{channel routing}, which involves removing adjacent vertex pairs from a graph while tracking the effect on channels. Channel routing is not always well-behaved, but certain other graph moves always preserve the number of channels in a graph. These are called \emph{channel-preserving moves} and are the focus of Section 6. Section 6 also introduces a useful graph move called \emph{diagonal contraction}. We give results on diagonal contraction and channels which generalize a number of known domino tiling parity results. We wrap up in Section 7 with remarks and directions for future work.

\section{Preliminaries}\label{secprelim}

All graphs in this paper are undirected, finite, and (unless otherwise indicated) contain no self-loops. Until Section 6, we further assume that there is at most one edge between any pair of vertices. This convention is used only to simplify notation; using suitable definitions, all results hold without it. If a graph is bipartite, we will consider its vertices to be colored black and white. Additionally, all matchings discussed will be perfect matchings, and thus the word ``perfect'' will be omitted in the future for brevity. For a graph $G=(V,E)$, $V$ denotes the vertex set, $E$ denotes the edge set, and $A$ denotes the adjacency matrix. We write $\edge{v}{w}$ to denote an (undirected) edge between $v$ and $w$; this is to contrast with the notation $\{v,w\}$ indicating a set of two vertices. Given a vertex $v$, the neighborhood of $v$ is 
\[ N(v) = \{v' \mid \edge{v}{v'}\in E  \}. \] 
For an edge $e$, we use the notation $G-e$ to denote the subgraph $(V,E-e)$. For a subset $S \subseteq V$, we use the notation $G-S$ to denote the subgraph of $G$ induced on $V-S$.  Recall that $\m$ is the number of matchings of $G$. As an exercise in this notation, we have the following proposition.
\begin{proposition}\label{matchprop}
	Let $G$ be any graph, and let $e=\edge{v_1}{v_2}$. Then
	$$m_G = m_{G-e} + m_{G-\{v_1,v_2\}} .$$
	Also, fix any vertex $v$. Then 
	$$m_G = \sum_{v'\in N(v)} m_{G-\{v,v'\}}.$$
\end{proposition}
\begin{proof}
	For the first relation, notice that matchings of $G-e$ are just those matchings of $G$ that do not use the edge $e$. The other matchings of $G$ do use $e$, and therefore for these matchings the vertices $v_1$ and $v_2$ are never in an edge with any vertex other than each other. Thus such matchings are equivalent to matchings of $G-\{v_1,v_2\}$, plus the edge $e$, and the first equation is shown.
	
	For the second relation, partition the set of matchings of $G$ based on the vertex that pairs with $v$ in the matching. By the same reasoning as in the last paragraph, the number of matchings in which $v$ pairs with $v'$ is $m_{G-\{v,v'\}}$. Summing over the possible pairings shows the claim.
\end{proof}

Sometimes we will be interested in connected planar graphs $G$. Such graphs admit a dual graph, with vertices given by the faces of $G$ and edges between faces separated by an edge in $G$. If the same face is on both sides of an edge of $G$, then that edge corresponds to a self-loop in the dual graph. The \emph{external face} of $G$ is the face which is unbounded, and all other faces are \emph{internal faces} of $G$. The \emph{reduced dual graph} of $G$ is the dual graph of $G$ with the vertex corresponding to the external face of $G$ removed. We say a vertex of $G$ is \emph{external} if it is incident to the external face, and we say it is \emph{internal} otherwise.

We begin by recalling Lov\'asz's original proof of Proposition \ref{even}. Later, in Section \ref{seccomb}, we will reprove this combinatorially for bipartite graphs.
\begin{proposition}[\cite{Lovasz1993}, Problem 5.18]\label{algproof}
    Let $G$ be a finite graph with no self-loops. Then $\m$ is even if and only if there is a nonempty vertex set $C\subseteq V$ such that every vertex in $G$ is adjacent to an even number of vertices in $C$.
\end{proposition}
\begin{proof}
    Let $A_2=(a_{ij})_{i,j=1}^n$ be the adjacency matrix of $G$ modulo $2$. Note that vertex sets as described in the proposition exist if and only if the kernel of $A_2$ is nonzero (since an element of the kernel is a sequence of $1$s and $0$s indexed by vertices, which may be treated as an indicator function for a vertex set that will necessarily have the desired property). 
    
    The determinant of $A_2$ is
    \[\det A_2 = \sum_{\sigma \in S_n} \mathrm{sgn}(\sigma) a_{1\sigma(1)}\cdots a_{n\sigma(n)}
    =\sum_{\sigma\in S_n} a_{1\sigma(1)}\cdots a_{n\sigma(n)},\]
    where the sum is over all permutations of $\{1,...,n\}$. (We may ignore the sign since we work over $\Z/2\Z$.) Because the adjacency matrix is symmetric, 
    \[a_{1\sigma(1)}\cdots a_{n\sigma(n)} =  a_{\sigma^{-1}(1)1}\cdots a_{\sigma^{-1}(n)n}
    = a_{1\sigma^{-1}(1)}\cdots a_{n\sigma^{-1}(n)}.\]
    Thus we may remove from the sum all pairs $\sigma\neq \sigma^{-1}$ as such pairs yield terms that sum to 0. We are left with the sum over involutions
    \[\sum_{\sigma=\sigma^{-1}} a_{1\sigma(1)}\cdots a_{n\sigma(n)}.\]
    Now, because $G$ has no self loops, if $\sigma$ satisfies $\sigma(i)=i$ for some $i$, then the corresponding term in the sum contains $a_{ii}=0$. Thus we are left with a sum over involutions of the vertices with no fixed points. For such an involution $\sigma$,
    $a_{1\sigma(1)}\cdots a_{n\sigma(n)}$ is $1$ if there exists a matching of $G$ that pairs vertex $i$ with vertex $\sigma(i)$ for all $i$, and $0$ otherwise. Hence we get the equality (in $\Z/2\Z$)
    \[\m = \det A_2.\]
    Thus $\m= 0$ in $\Z/2\Z$ if and only if $\ker A_2 \neq 0$,
    and we are done.
\end{proof}

Now we discuss some important algebraic constructions. Our main tools for the remainder of the section are the Kasteleyn matrix of a graph and the Smith decomposition of a matrix. 

\begin{proposition}[Kasteleyn \cite{Kasteleyn1963}]\label{kast}
	Let $G$ be a planar graph with adjacency matrix $A=(a_{ij})$. Then there exists a matrix $K$ such that $K=(\pm a_{ij})$, and
	$$\det K = m_G^2.$$
\end{proposition}

The matrix $K$ is called a \textit{Kasteleyn matrix}, and if a (not necessarily planar) graph $G$ admits such a $K$, then $G$ is said to have a \textit{Kasteleyn signing}.

For bipartite graphs we can be more specific.

\begin{definition}
	Let $G$ be a bipartite graph with adjacency matrix $A$. The \textit{bipartite adjacency matrix} of $G$ is the submatrix of $A$ formed by selecting rows from $A$ associated to white vertices and columns from $A$ associated to black vertices. 
\end{definition}

When $G$ is bipartite, the adjacency matrix $A$ can be written in block matrix form as
\[ A = \begin{pmatrix}0 & B \\ B^T & 0  \end{pmatrix}, \]
where $B$ is the bipartite adjacency matrix.
If in addition $G$ is planar and has the same number of black and white vertices, then there is an analogue to the Kasteleyn matrix called the \textit{bipartite Kasteleyn matrix}.

\begin{proposition}[Percus \cite{Percus1969}]
	Let $G$ be a bipartite planar graph with bipartite adjacency matrix $B=(b_{ij})$. If $B$ is square, then there exists a matrix $H$, the \emph{bipartite Kasteleyn matrix}, such that $H=(\pm b_{ij})$, and
	$$\det H = m_G.$$
\end{proposition}

We will want to diagonalize these matrices over the integers. The canonical tool for doing so is called \textit{Smith normal form}.

\begin{proposition}\label{snf}
	Let $A$ be a matrix over a principal ideal domain (PID) $\mathcal R$. Then there exist matrices $S,D,T$ over $\mathcal R$ with the following properties:
	
	\begin{enumerate}[label=(\roman*)]
		\item $A = SDT$.
		
		\item $S$ and $T$ are invertible over $\mathcal R$. For $\mathcal R=\Z$, this means $\det S, \det T = \pm 1$.
		
		\item $D$ is diagonal, with diagonal entries $\alpha_1,...,\alpha_n$ satisfying 
		$$\alpha_1 \text{ divides } \alpha_2 \text{ divides } ... \text{ divides } \alpha_n .$$
	\end{enumerate}

	The matrices $S,D,T$ are called a \emph{Smith decomposition} for $A$.
	
\end{proposition}

For our purposes, $\mathcal R$ will be the integers or a finite field. Smith decompositions have many useful properties. See, e.g., \cite{Kuperberg2002} or \cite{Stanley2016} for more background and combinatorial applications of Smith decompositions.  The next result follows directly from Proposition \ref{snf}(i)--(ii).

\begin{proposition}\label{smithprop}
	Let $A = SDT$ be a Smith decomposition for a square matrix $A$ over a PID $\mathcal R$. Then
	$$ \ker A \cong \ker D \text{ as abelian groups and } \det A = u\det D \text{ for some unit $u\in \mathcal R$.} $$
\end{proposition}

We are now prepared to study 2-divisibility of $\m$ for graphs $G$ with a Kasteleyn matrix $K$. As described above, planarity is a sufficient condition for $G$ to have a Kasteleyn signing. Using the Smith normal form of $K$, we will find 2-divisibility results for such graphs that we develop further in the next section.

\begin{definition}
	For an integral matrix $A$, define the \textit{reduction of $A$ modulo $2$} to be the matrix $A_2$ over $\Z/2\Z$ given by reducing the entries of $A$ modulo $2$ and considering them as elements of $\Z/2\Z$. 
	
	Define $\ker_2A$, the \textit{$2$-kernel} of $A$, to be the kernel of $A_2$ as a vector space over $\Z/2\Z$. 
\end{definition}

Notice that for an adjacency matrix $A$ and corresponding Kasteleyn matrix $K$, we have $A_2=K_2$. This follows from the definition of $K$ as a signed version of $A$. This is a key observation that will allow us to translate our algebraic results in this section into geometric results in Section \ref{secchannels}. Before that, let us see what we can learn from reducing the Kasteleyn matrix mod 2. Let the \textit{$2$-nullity} of a matrix be the dimension of its $2$-kernel:
$$\nul_2A := \dim\ker_2 A.$$

\begin{lemma}\label{2nul}
	Let $A$ be a square matrix with integer entries. Then
	$$2^{\nul_2 A} \mid \det A.$$
\end{lemma}
\begin{proof}
	Let $A= S D T$ be a Smith decomposition of $A$. Reducing modulo 2 gives $A_2 = S_2 D_2 T_2$. One may check that this is a Smith normal form of $A_2$. Then by Proposition \ref{smithprop}, the kernel of $A_2$ is isomorphic to the kernel of $D_2$. Set $k:=\nul_2 A = \nul_2 D$. Since $D_2$ is diagonal, its kernel has a basis consisting of standard basis vectors that indicate columns where the diagonal entry is $0$. Therefore there are exactly $k$ such entries. Since $0$ entries in $D_2$ correspond to even integral entries in $D$, there are exactly $k$ even entries on the diagonal of $D$. Thus the determinant of $D$ contains at least $k$ factors of $2$, and the result follows by Proposition \ref{smithprop}.
\end{proof}

Applying this lemma to the Kasteleyn matrix or the bipartite Kasteleyn matrix will let us use the 2-kernel to find powers of two in the number of matchings of a graph.

\begin{theorem}\label{2div}
	Let $G$ be a graph with a Kasteleyn signing (e.g.\ a planar graph). If $A$ is the adjacency matrix of $G$, then 
	\begin{equation*}\label{2diveq}
	2^{\nul_2 A} \text{ divides } m_G^2.
	\end{equation*}
\end{theorem}
\begin{proof}
	As remarked above, $A_2=K_2$, so in particular $\nul_2 A=\nul_2K$. Additionally, by Proposition \ref{kast}, the determinant of $K$ is $m_G^2$. Thus by Lemma \ref{2nul}, 
	\[2^{\nul_2A} = 2^{\nul_2 K} \text{ divides } \det K = m_G^2 .\qedhere\]
\end{proof}

\begin{theorem}\label{2divbi}
	Let $G$ be a bipartite graph with a Kasteleyn signing. If $B$ is the bipartite adjacency matrix of $G$, then 
	\begin{equation*}
	2^{\nul_2 B} \text{ divides } m_G.
	\end{equation*}
\end{theorem}
\begin{proof}
	If $B$ is not square, then $\m=0$ and the claim holds trivially. Otherwise, the proof is the same as the previous theorem using the bipartite Kasteleyn matrix.
\end{proof}

We therefore may deduce powers of 2 dividing $m_G$ by finding elements of the 2-kernel of $A$. The remainder of the paper details how this can be done.

\section{Channels}\label{secchannels}
Let $G$ be a graph with adjacency matrix $A$. Then a vector $x$ in $\ker_2 A$ has entries in $\Z/2\Z$ and can be lifted to a vector $\tilde x$ with entries $0,1\in \Z$. The condition $A_2x = 0$ then becomes $A\tilde x = 2y$ for some integral vector $y$. Because each row of $x$ corresponds to a vertex in $G$, we may interpret $x$ as the indicator function for a vertex set $C$, where a row with a $1$ indicates the vertex is in $C$ and a row with a $0$ indicates the vertex is not in $C$. This leads to the following interpretation of 2-kernel elements.

\begin{definition}
	Let $G = (V, E)$ be any graph. A \textit{channel} is a set $C$ of vertices such that every vertex in $G$ is adjacent to an even number of vertices in $C$. In other words, letting $N(v)$ denote the neighborhood of $v$, a channel satisfies
	$$|N(v)\cap C| \text{ is even, for all $v\in V$}.$$
	Let the set of channels in $G$ be denoted $\mathcal{C}(G)$. If $G$ is bipartite, then let $\cb(G)$ (resp. $\cw(G)$) be the subset of $\mc(G)$ consisting of channels that use only black (resp. white) vertices from $G$.
\end{definition}

The 2-kernel also has an additive structure as a $\Z/2\Z$ vector space. This transfers to $\mc(G)$ by defining the sum of $C_1,C_2\in\mc(G)$ to be
$$C_1 \oplus C_2 := (C_1\cup C_2)-(C_1\cap C_2)\text{, the symmetric difference of $C_1$ and $C_2$} .$$
Note that in the bipartite case, both $\cb(G)$ and $\cw(G)$ are subspaces of $\mc(G)$.

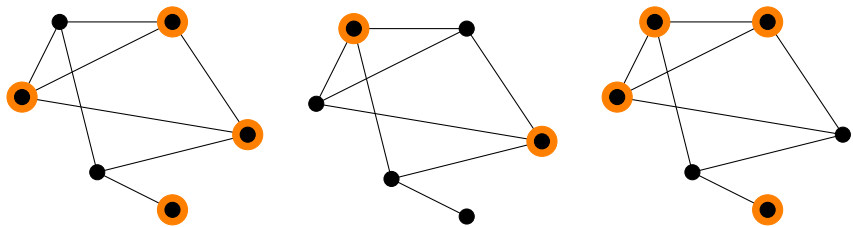
\begin{figure}[htb]
    \centering
    \begin{tikzpicture}
        \pgfmathsetmacro{\radius}{0.1}
        \node (a) at (0,0){};
        \node (b) at (0.5,1){};
        \node (c) at (2,1){};
        \node (d) at (3,-0.5){};
        \node (e) at (2,-1.5){};
        \node (f) at (1,-1){};
        \draw (a.center) -- (c.center) -- (d.center) -- (a.center);
        \draw (a.center) -- (b.center);
        \draw (b.center) -- (c.center);
        \draw (b.center) -- (f.center);
        \draw (e.center) -- (f.center);
        \draw (f.center) -- (d.center);
        \draw[fill,orange] (a) circle (2*\radius);
        \draw[fill=black] (a) circle (\radius);
        \draw[fill=black] (b) circle (\radius);
        \draw[fill,orange] (c) circle (2*\radius);
        \draw[fill,orange] (d) circle (2*\radius);
        \draw[fill,orange] (e) circle (2*\radius);
        \draw[fill=black] (c) circle (\radius);
        \draw[fill=black] (d) circle (\radius);
        \draw[fill=black] (e) circle (\radius);
        \draw[fill=black] (f) circle (\radius);
    \end{tikzpicture} \hspace{1em}
    \begin{tikzpicture}
        \pgfmathsetmacro{\radius}{0.1}
        \node (a) at (0,0){};
        \node (b) at (0.5,1){};
        \node (c) at (2,1){};
        \node (d) at (3,-0.5){};
        \node (e) at (2,-1.5){};
        \node (f) at (1,-1){};
        \draw (a.center) -- (c.center) -- (d.center) -- (a.center);
        \draw (a.center) -- (b.center);
        \draw (b.center) -- (c.center);
        \draw (b.center) -- (f.center);
        \draw (e.center) -- (f.center);
        \draw (f.center) -- (d.center);
        \draw[fill,orange] (b) circle (2*\radius);
        \draw[fill,orange] (d) circle (2*\radius);
        \draw[fill=black] (a) circle (\radius);
        \draw[fill=black] (b) circle (\radius);
        \draw[fill=black] (c) circle (\radius);
        \draw[fill=black] (d) circle (\radius);
        \draw[fill=black] (e) circle (\radius);
        \draw[fill=black] (f) circle (\radius);
    \end{tikzpicture} \hspace{1em}
    \begin{tikzpicture}
        \pgfmathsetmacro{\radius}{0.1}
        \node (a) at (0,0){};
        \node (b) at (0.5,1){};
        \node (c) at (2,1){};
        \node (d) at (3,-0.5){};
        \node (e) at (2,-1.5){};
        \node (f) at (1,-1){};
        \draw (a.center) -- (c.center) -- (d.center) -- (a.center);
        \draw (a.center) -- (b.center);
        \draw (b.center) -- (c.center);
        \draw (b.center) -- (f.center);
        \draw (e.center) -- (f.center);
        \draw (f.center) -- (d.center);
        \draw[fill,orange] (a) circle (2*\radius);
        \draw[fill,orange] (b) circle (2*\radius);
        \draw[fill,orange] (c) circle (2*\radius);
        \draw[fill,orange] (e) circle (2*\radius);
        \draw[fill=black] (a) circle (\radius);
        \draw[fill=black] (b) circle (\radius);
        \draw[fill=black] (c) circle (\radius);
        \draw[fill=black] (d) circle (\radius);
        \draw[fill=black] (e) circle (\radius);
        \draw[fill=black] (f) circle (\radius);
    \end{tikzpicture}
    \caption{A graph $G$ with its three nonempty channels indicated by shading. Any two of these form a basis for the space $\mc(G)$.}
    \label{fig:channelexample}
\end{figure}

The discussion above then implies the following, which was hinted at in the proof of Proposition \ref{algproof}.
\begin{lemma}\label{ker2chan}
	The spaces $\mathcal{C}(G)$ and $\ker_2 A$ are isomorphic vector spaces over $\Z/2\Z$.
\end{lemma}

If $G$ is additionally bipartite with bipartite adjacency matrix $B$, then we can refine Lemma \ref{ker2chan} by accounting for vertex colors.
\begin{lemma}\label{sqbipart}
	If $G$ is bipartite, then there are $\Z/2\Z$-vector space isomorphisms
	\[ \cb(G) \cong \ker_2 B \quad \text{ and } \quad \cw(G) \cong \ker_2 B^T  .\]
	Furthermore,
	\[  \mathcal{C}(G) = \cb(G) \oplus \cw(G) . \] 
\end{lemma}
\begin{proof}
The identification of channels using only black vertices with elements of $\ker_2 B$ follows the same lines as the discussion at the beginning of the section. The claim for white vertices is the same, noting that taking the transpose of $B$ is equivalent to swapping the two vertex colors. 


For the final claim, the set of black vertices in any channel $C$ is also a channel: in a bipartite graph, the set of black vertices adjacent to a given vertex is either the empty set (which has even cardinality) if the vertex is black, or its entire neighborhood if the vertex is white. Similarly the white vertices in $C$ also form a channel. As a result, any channel may be split into its black part and its white part. Since $\cb(G)$ and $\cw(G)$ have trivial intersection, this proves the direct sum decomposition.
\end{proof}

\begin{lemma}
Let $|V_B|$ and $|V_W|$ be the number of black vertices and white vertices in a bipartite graph $G$, respectively. Then 
\[  \dim \cb(G) - \dim\cw(G) = |V_B| - |V_W|.  \]
In particular, if $G$ has the same number of black and white vertices, then
\[ \dim\cb(G) = \dim\cw(G) =  \frac{1}{2}\dim\mc(G). \] 
\end{lemma}
\begin{proof}
Let $B$ be the bipartite adjacency matrix of $G$. We have 
\[\dim \cb(G) = \nul_2 B \tag{1}\label{eqn:kereq}\]
and
\[\dim \cw(G) = \nul_2 B^T = \dim \mathrm{coker}_2 B = |V_W| - \dim \mathrm{im}_2 B, \tag{2}\label{eqn:cokereq}\]
where $\mathrm{coker}_2$ and $\mathrm{im}_2$ refer to cokernel and image of $B$ as a matrix over $\Z/2\Z$. The rank-nullity theorem for $B$ may be stated as
\[ \nul_2 B + \dim \mathrm{im}_2 B = |V_B|. \]
The claim then follows by subtracting equation (\ref{eqn:kereq}) from equation (\ref{eqn:cokereq}) and applying rank-nullity.
\end{proof}

Combining the observations above with the 2-divisibility results from the last section, we can now prove our main theorem.

\begin{theorem}[Channeling 2s]\label{twos}
	If $C_1,...,C_n$ are linearly independent channels in a graph $G$ with a Kasteleyn signing, then 
	$$2^n \text{ divides } m_G^2.$$
	If additionally $G$ is bipartite, and $C_1,...,C_n\in \cb(G)$, then
	$$2^n \text{ divides } m_G.$$
\end{theorem}
\begin{proof}
Let $A$ be the adjacency matrix of $G$. By Lemma \ref{ker2chan},  $C_1,...,C_n$ may be viewed as $n$ linearly independent elements of $\ker_2 A$. It follows that $\nul_2 A \geq n$. Then the first claim results from applying Theorem \ref{2div}.

The second claim follows similarly, using Lemma \ref{sqbipart} and Theorem \ref{2divbi}.
\end{proof}

\begin{remark}\label{kastsignnec}
	Despite the fact that Proposition \ref{algproof} holds for an arbitrary graph, Theorem \ref{twos} does not. For example, the complete bipartite graph $K_{3,3}$ has $|\cb(K_{3,3})| = 2^2$, but $m_{K_{3,3}}=6$ is not divisible by $2^2$. Thus the assumption of a Kasteleyn signing for $G$ cannot be weakened much further.
\end{remark}

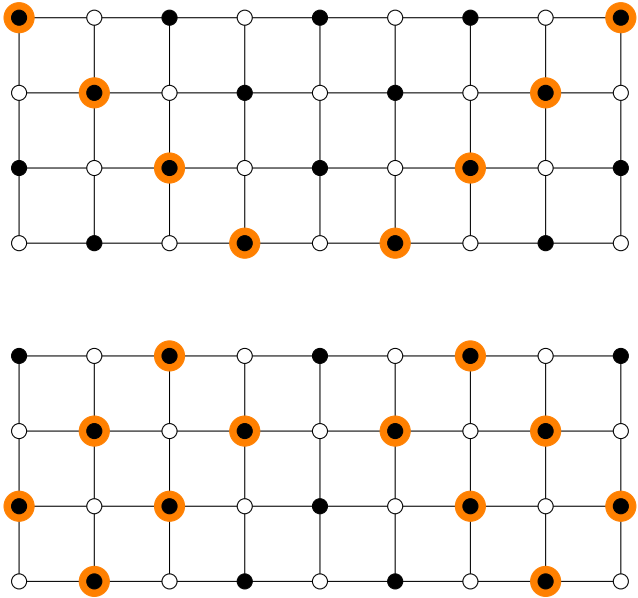
\begin{figure}[H]
    \centering
    \begin{tikzpicture}
        \draw[thin, black] (0,0) grid (8,3);
        \foreach \x in {0, 1, 2, 3} {
        	\draw[orange,fill] (\x,3-\x) circle (.2);
        	\draw[orange,fill] (8-\x,3-\x) circle (.2);
        }
    	
    	\begin{scope}[yscale=-1,shift={(0,-3)}]
    	\vertices{0}{0}{8}{3};
    	\end{scope}
    \end{tikzpicture} \\ \vspace{3em}
    \begin{tikzpicture}
        \draw[thin, black] (0,0) grid (8,3);
        \foreach \x in {0, 1, 2} {
        	\draw[orange,fill] (\x,1+\x) circle (.2);
        	\draw[orange,fill] (\x+1,\x) circle (.2);
        	\draw[orange,fill] (7-\x,\x) circle (.2);
        	\draw[orange,fill] (8-\x,1+\x) circle (.2);
        	
        }
    	\begin{scope}[yscale=-1,shift={(0,-3)}]
		\vertices{0}{0}{8}{3};
		\end{scope}
    \end{tikzpicture}
    \caption{A basis for $\cb(\mathcal R_{4\times 9})$.}
    \label{fig:application}
\end{figure}

\begin{example1} \label{rectapp}
	Let $\mathcal R_{m\times n}$ denote the $m\times n$ rectangular grid graph. The shading in Figure \ref{fig:application} shows a basis for $\cb(\mathcal R_{4\times 9})$. By channeling 2s, we have that $2^2$ divides $m_{\mathcal R_{4\times 9}}$. And indeed, $m_{\mathcal R_{4\times 9}}=6336$ is divisible by 4. Note, however, that 6336 is also divisible by $2^6$---Theorem \ref{twos} gives only a \textit{lower bound} on the power of 2 dividing $m_G$.
\end{example1}

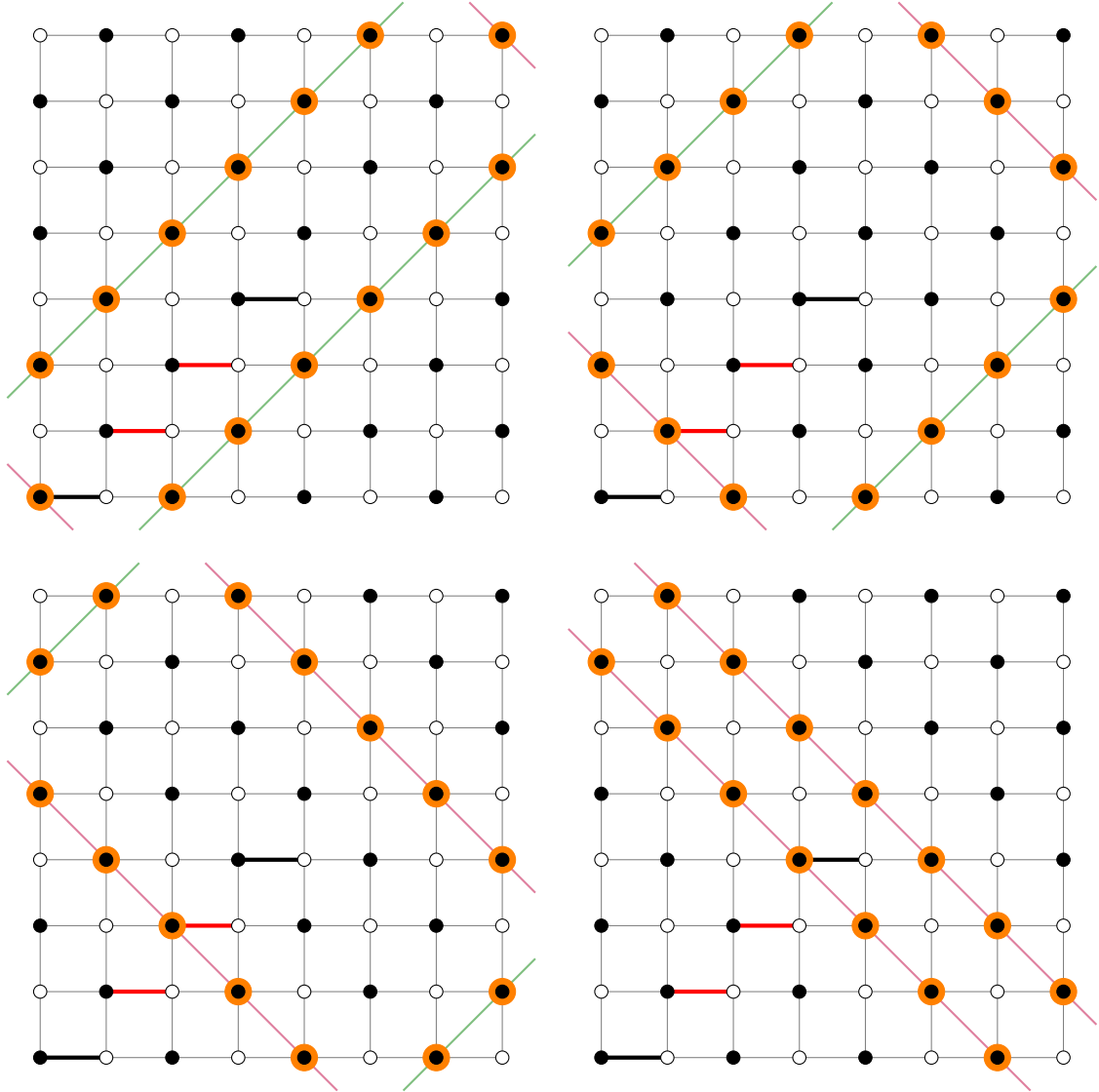
\begin{figure}
	\centering
	\begin{tikzpicture}[scale=.9]
	\draw[thin,opacity=.5] (0,0) grid (7,7);
	\draw[red, ultra thick] (1,1) -- (2,1);
	\draw[red, ultra thick] (2,2) -- (3,2);
	\draw[ultra thick] (0,0) -- (1,0);
	\draw[ultra thick] (3,3) -- (4,3);

	\draw[thick, purple, opacity=.5] (-.5,.5) -- (.5,-.5);
	\draw[thick, purple, opacity=.5] (6.5,7.5) -- (7.5,6.5);
	\draw[thick, black!50!green, opacity=.5] (-.5,1.5) -- (5.5,7.5);
	\draw[thick, black!50!green, opacity=.5] (1.5,-.5) -- (7.5,5.5);
	
	\foreach \x in {0, 1,..., 5} {
		\draw[orange,fill] (\x,\x+2) circle (.2);
	}
	\foreach \x in {2, 3,..., 7} {
		\draw[orange,fill] (\x,\x-2) circle (.2);
		.,}
	\draw[orange,fill] (0,0) circle (.2);
	\draw[orange,fill] (7,7) circle (.2);
	\vertices{0}{0}{7}{7};
	
	\begin{scope}[shift={(8.5,0)}]
		\draw[thin,opacity=.5] (0,0) grid (7,7);
		\draw[red, ultra thick] (1,1) -- (2,1);
		\draw[red, ultra thick] (2,2) -- (3,2);
		\draw[ultra thick] (0,0) -- (1,0);
		\draw[ultra thick] (3,3) -- (4,3);

		\draw[thick, purple, opacity=.5] (-.5,2.5) -- (2.5,-.5);
		\draw[thick, purple, opacity=.5] (4.5,7.5) -- (7.5,4.5);
		\draw[thick, black!50!green, opacity=.5] (-.5,3.5) -- (3.5,7.5);
		\draw[thick, black!50!green, opacity=.5] (3.5,-.5) -- (7.5,3.5);
		
		\foreach \x in {0, 1,..., 3} {
			\draw[orange,fill] (\x,\x+4) circle (.2);
		}
		\foreach \x in {4, 5,..., 7} {
			\draw[orange,fill] (\x,\x-4) circle (.2);
			.,}
		\foreach \x in {0, 1, 2} {
			\draw[orange,fill] (\x,2-\x) circle (.2);
			.,}
		\foreach \x in {0, 1, 2} {
			\draw[orange,fill] (\x+5,7-\x) circle (.2);
			.,}
		\vertices{0}{0}{7}{7};
	\end{scope}
	
	\begin{scope}[shift={(8.5,-8.5)}]
	\draw[thin,opacity=.5, black] (0,0) grid (7,7);
	\draw[red, ultra thick] (1,1) -- (2,1);
	\draw[red, ultra thick] (2,2) -- (3,2);
	\draw[ultra thick] (0,0) -- (1,0);
	\draw[ultra thick] (3,3) -- (4,3);

	\draw[thick, purple, opacity=.5] (-.5,6.5) -- (6.5,-.5);
	\draw[thick, purple, opacity=.5] (0.5,7.5) -- (7.5,0.5);
		
	\foreach \x in {0, 1,..., 6} {
		\draw[orange,fill] (\x,6-\x) circle (.2);
	}
	\foreach \x in {0,1,..., 6} {
		\draw[orange,fill] (\x+1,7-\x) circle (.2);
		.,}
	\vertices{0}{0}{7}{7};
	\end{scope}
	
	\begin{scope}[shift={(0,-8.5)}]
	\draw[thin, opacity=.5, black] (0,0) grid (7,7);
	\draw[red, ultra thick] (1,1) -- (2,1);
	\draw[red, ultra thick] (2,2) -- (3,2);
	\draw[ultra thick] (0,0) -- (1,0);
	\draw[ultra thick] (3,3) -- (4,3);

	\draw[thick, purple, opacity=.5] (-.5,4.5) -- (4.5,-.5);
	\draw[thick, purple, opacity=.5] (2.5,7.5) -- (7.5,2.5);
	\draw[thick, black!50!green, opacity=.5] (-.5,5.5) -- (1.5,7.5);
	\draw[thick, black!50!green, opacity=.5] (5.5,-.5) -- (7.5,1.5);

	\foreach \x in {0, 1} {
		\draw[orange,fill] (\x,6+\x) circle (.2);
	}
	\foreach \x in {6, 7} {
		\draw[orange,fill] (\x,-6+\x) circle (.2);
		.,}
	\foreach \x in {0, 1, ..., 4} {
		\draw[orange,fill] (\x,4-\x) circle (.2);
		.,}
	\foreach \x in {0, 1, ..., 4} {
		\draw[orange,fill] (\x+3,7-\x) circle (.2);
		.,}
	\vertices{0}{0}{7}{7};
	\end{scope}
	\end{tikzpicture}
	\caption{A basis for $\cb(\mathcal R_{8\times 8})$. The diagonal segments which intersect the step diagonal are shown in purple, while the diagonal segments which are parallel to the step diagonal are in green. Note that the channels in the top left and bottom right are also channels of the graph with the two red edges removed.}\label{stepchannel}
\end{figure}

 In Section \ref{secbounce}, we shall employ billiard paths to count channels in a rectangle grid graph of arbitrary size.	We already have the results we need, however, to give a lower bound supporting Pachter's conjecture.
	
\begin{proposition}\label{step-diag}
	Let $\mathcal{R}_{2r\times 2r}$ be the $2r\times 2r$ grid graph shown in Figure \ref{fig:first}. If $G$ is constructed from $\mathcal R_{2r\times 2r}$ by deleting any $k$ of the highlighted edges in the figure, then
    $$2^{r-k} \text{ divides } \m.$$
\end{proposition}
\begin{proof}
 For a $2r\times 2r$ square grid graph, we will find that the space $\cb(\mathcal R_{2r\times 2r})$ has $r$ independent channels which each intersect exactly one of vertex pairs removed from the step diagonal. Thus removing $k$ of the step diagonal vertex pairs interrupts just $k$ of the channels, so the other $r-k$ channels will still be present. The result then follows from channeling 2s once we construct these channels. 
 
 Consider Figure \ref{stepchannel}. To construct each channel, pick a black vertex $b$ on the step diagonal. In the figure, these are the vertices along the diagonal from the bottom left to the top right. Our channel will consist of four (possibly empty) diagonal segments of vertices. Two segments will intersect the step diagonal transversely, one at $b$ in the bottom left and one at $b$'s mirror image in the top right. The other two segments will be parallel to the step diagonal and placed so that the channel vertices along the sides of the grid graph each have one intervening vertex between them. Each such intervening vertex is adjacent to one endpoint from two diagonal segments. All other vertices are adjacent to either 0 or 2 vertices from each diagonal segment. Hence each vertex of the graph is adjacent to an even number of vertices from our diagonal segments, so these vertices form a channel.
 
 By following this construction for each $b$ lying along the lower $r$ black vertices of the step diagonal, we create $r$ channels that each contain a different vertex from the step diagonal, and therefore they must be independent. 
 \end{proof}
\section{Billiards and Channels}\label{secbounce}

For an arbitrary graph, it is not at all clear how to construct or count its channels without solving the requisite linear system. In this section we give a geometric approach to channel construction based on a phenomenon that can be observed in the channels of a rectangle grid graph. 

\subsection{Billiard nests}

In the rectangle, we note that channels tend to form along diagonal lines as in the following figure. This pattern was studied by Tomei and Vieira in \cite{Tomei2002}, where they described it in terms of polygonal tilings of the rectangle. We propose an alternative description. 

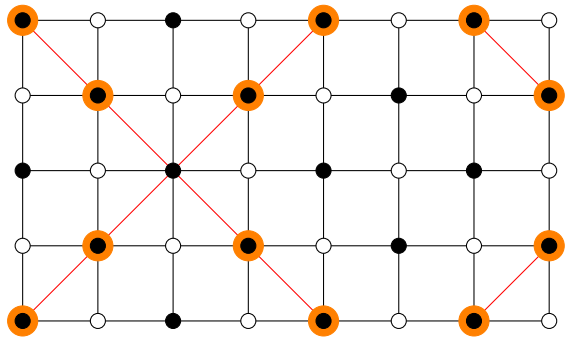
\begin{figure}[H]
	\centering
\begin{tikzpicture}
	\draw (0,0) grid (7,4);
	\pgfmathsetmacro{\chr}{0.2}
	
	\draw[red] (0,0) -- (4,4);
	\draw[red] (0,4) -- (4,0);
	\draw[red] (6,0) -- (7,1);
	\draw[red] (6,4) -- (7,3);
	
	\foreach \coord in %
	{(0,0),(1,1),(1,3),(0,4),(3,1),(3,3),(4,0),(4,4),(6,0),(6,4),(7,1),(7,3)}
	\draw[fill,orange] \coord circle (\chr);
	
	\vertices{0}{0}{7}{4}
	
\end{tikzpicture}
\end{figure}
By extending these diagonals, we find that these lines form a path which reflects off the edges of a larger rectangle, as shown below. The channel vertices are vertices in the interior of this larger rectangle which intersect exactly one line from this path.
\begin{figure}[H]
	\centering
	\begin{tikzpicture}
	\draw[dotted, opacity=.5] (-1,-1) grid (8,5);
	\draw (0,0) grid (7,4);
	\pgfmathsetmacro{\chr}{0.2}
	
	\draw[red] (-1,-1) -- (5,5)--(8,2)--(5,-1)--(-1,5);
	
	\foreach \coord in %
	{(0,0),(1,1),(1,3),(0,4),(3,1),(3,3),(4,0),(4,4),(6,0),(6,4),(7,1),(7,3)}
	\draw[fill,orange] \coord circle (\chr);
	
	\vertices{0}{0}{7}{4}
	
	\end{tikzpicture}
\end{figure}

Such paths either form loops or start and end on distinct corners. Notice that we may recover the path by remembering only the faces that it passes through and the fact that it passes through black vertices. The face information is shown in the next figure. 

\begin{figure}[H]
	\centering
	\begin{tikzpicture}
	\draw[dotted, opacity=.5] (-1,-1) grid (8,5);
	\draw (0,0) grid (7,4);
	\pgfmathsetmacro{\chr}{0.2}
	
	\foreach \x in {-1, 0, ..., 4}
	\fill[opacity = .1, red] (\x,\x) rectangle ++(1,1);
	\foreach \x in {-1, 0, ..., 4}
	\fill[opacity = .1, red] (\x,3-\x) rectangle ++(1,1);
	\foreach \x in {-1, 0, ..., 1}
	\fill[opacity = .1, red] (6+\x,\x) rectangle ++(1,1);
	\foreach \x in {-1, 0, ..., 1}
	\fill[opacity = .1, red] (6+\x,3-\x) rectangle ++(1,1);
	
	\foreach \coord in %
	{(0,0),(1,1),(1,3),(0,4),(3,1),(3,3),(4,0),(4,4),(6,0),(6,4),(7,1),(7,3)}
	\draw[fill,orange] \coord circle (\chr);
	
	\vertices{0}{0}{7}{4}
	
	\end{tikzpicture}
\end{figure}

It turns out to be natural to consider arbitrary collections of billiard paths---what we will call a \emph{billiard nest}. Viewing the nest as a set of faces will allow us to isolate the abstract properties of billiards which are of combinatorial importance. From this viewpoint, billiard nests can be defined on a large class of graphs, called \emph{inner semi-Eulerian graphs}. 

\begin{definition}\label{def:bp}
		We say that a bipartite planar graph $G$ is \textit{inner semi-Eulerian} if every internal black vertex of $G$ has even degree. Let $G$ be inner semi-Eulerian and let $F$ denote the set of internal faces of $G$. We say that a subset of faces $B \subseteq F$ of $G$ is a \textit{billiard nest} if the following hold:
		\begin{itemize}
		    \item if $b$ is an internal black vertex, then either all faces incident to $b$ are in $B$, no faces incident to $b$ are in $B$, or every second face incident to $b$ is in $B$.
			
			\item if $b$ is an external black vertex, then either all internal faces incident to $b$ are in $B$ or no internal faces incident to $b$ are in $B$. 
		\end{itemize}
		Denote the set of billiard nests in $G$ by $\bp(G)$.
\end{definition}
When $G$ is an induced subgraph of the square lattice with all internal faces being unit squares, this roughly agrees with the intuitive notion of billiard nests as (the faces containing) paths traced out by a collection of billiard balls since at an internal black vertex, the ball will pass from one face to the opposite one at that vertex. (At an external black vertex the billiard ball can be thought of as splitting into all possible directions.) Note that $\varnothing$ and $F$ are trivially billiard nests for every graph. As with channels, we may define the sum of two billiard nests to be their symmetric difference, making $\bp(G)$ a vector space over $\Z/2\Z$.

There is a canonical basis for $\bp(G)$ such that the basis billiard nests are mutually disjoint. Indeed, define a graph $G_\bp$ with vertex set $F$ and edges between $f$ and $f'$ if they satisfy one of the following:
\begin{itemize}
	\item $f$ and $f'$ are incident to the same internal black vertex $b$ and are separated by an even number of edges incident to $b$.
	
	\item $f$ and $f'$ are incident to the same external black vertex $b$.
\end{itemize}
Then the connected components of $G_\bp$ are independent billiard nests that span $\bp(G)$. This is called the \textit{path basis} for $\bp(G)$, and its elements are called \emph{billiard paths}. Later we will demonstrate an efficient algorithm to find the path basis for certain graphs. This is particularly useful since, as we shall soon see, billiard nests in $G$ are equivalent to channels in a certain subgraph of $G$.

\begin{definition}
	Let $G$ be inner semi-Eulerian. The \textit{inner subgraph} of $G$, denoted $G'$, is the induced subgraph on the internal vertices of $G$. Similarly, the \textit{outer subgraph} of $G$ is the induced subgraph on the external vertices of $G$. 
	
	Given any inner semi-Eulerian graph $H$, an \textit{outer completion} of $H$ is an inner semi-Eulerian graph $G$ such that $G' = H$ and such that the outer subgraph of $G$ is a simple cycle.
\end{definition}

We can always construct an outer completion for an inner semi-Eulerian graph $H$ by taking a copy of the boundary of $H$, expanding the copy so that $H$ lies within it, and adding edges between the two copies of the boundary as needed to make the graph inner semi-Eulerian. This is described in detail in the next proposition and the following example.

\begin{proposition}\label{outercomp}
	Any inner semi-Eulerian graph $H$ has an outer completion.
\end{proposition}
\begin{proof}
	Without loss of generality, let $H$ be connected. Set $n$ to be the degree of the external face $f$ of $H$ in the dual graph $H^\vee$ of $H$. The cyclic ordering on the faces adjacent to $f$ in $H^\vee$ induces a cyclic ordering on the external vertices of $H$ (possibly with repetition). Construct a simple cycle $Y$ with $n$ vertices and edges, and embed $Y$ in the plane such that $H$ is in the interior of $Y$. Then the cyclic order of the external vertices of $H$ surjectively assigns to each vertex of $Y$ an external vertex of $H$. For each external vertex $v$ of $H$, if $v$ has odd degree, then place an edge between $v$ and one of the vertices of $Y$ to which $v$ is assigned.
	
	Call the resulting graph $G$. Then $G$ is bipartite since we may color each vertex in $Y$ with the color opposite of the vertex assigned to it. Also $G$ is planar: the edges we introduced between $Y$ and $H$ may be embedded so that they do not cross since they connect to the boundary of $H$ in the same order as they do to $Y$. Finally, $G$ is inner semi-Eulerian since the construction of $G$ makes all external vertices of $H$ have even degree, and all other internal black vertices of $G$ are internal vertices of $H$.
\end{proof}

\begin{example1}
	Let $H$ be the following graph.
	
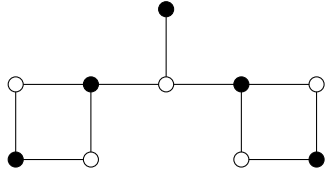
\begin{figure}[H]
		\centering
		\begin{tikzpicture}
		\draw (0,0) grid (1,1);
		\draw (3,0) grid (4,1);
		\draw (1,1) -- (3,1);
		\draw (2,1) -- (2,2);
		\vertices{0}{0}{1}{1};
		\vertices{2}{1}{2}{2};
		\vertices{3}{0}{4}{1};
		\end{tikzpicture}
	\end{figure}
The degree of the outer face $f$ is 14, and the sequence of external vertices visited as we traverse the boundary is indicated by the numbering below. Construct a simple cycle $Y$ of length 14 and embed it so $H$ lies in the interior face of $Y$.
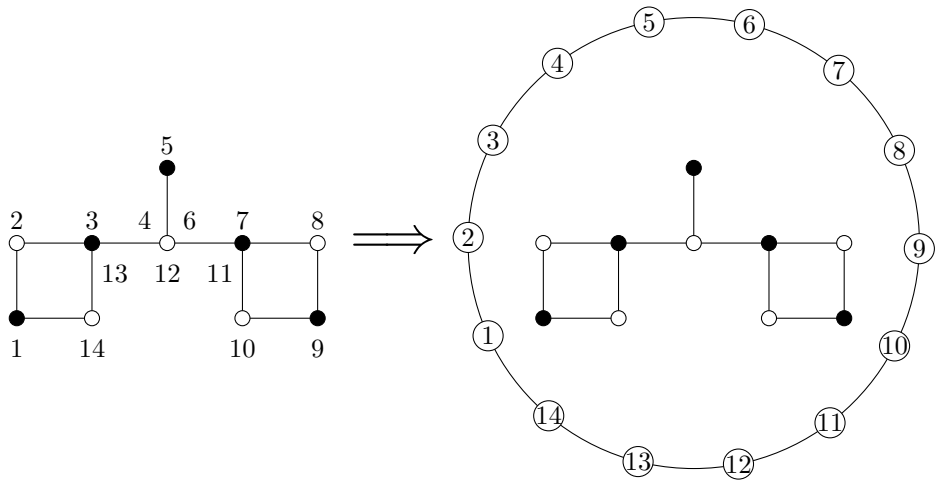
\begin{figure}[H]
	\centering
	\begin{tikzpicture}
	\draw (0,0) grid (1,1);
	\draw (3,0) grid (4,1);
	\draw (1,1) -- (3,1);
	\draw (2,1) -- (2,2);
	\vertices{0}{0}{1}{1};
	\vertices{2}{1}{2}{2};
	\vertices{3}{0}{4}{1};
	
	\draw (0,-.4) node{1};
	\draw (0,1.3) node{2};
	\draw (1,1.3) node{3} (1.7,1.3) node{4};
	\draw (2,2.3) node{5} (2.3,1.3) node{6};
	\draw (3,1.3) node{7} (4,1.3) node{8} (4,-.4) node{9} (3,-.4) node {10};
	\draw (2.7,.6) node{11} (2, .6) node {12} (1.3, .6) node {13};
	\draw (1,-.4) node {14};
	
	\begin{scope}[shift={(7,0)}]
	
	\draw (0,0) grid (1,1);
	\draw (3,0) grid (4,1);
	\draw (1,1) -- (3,1);
	\draw (2,1) -- (2,2);
	\vertices{0}{0}{1}{1};
	\vertices{2}{1}{2}{2};
	\vertices{3}{0}{4}{1};
	
	\draw (-2, 1) node[scale=2] {$\implies$};
	\draw (2,1) circle (3);
	
	\foreach \x in {1,2,...,14} {
		\node (\x) at ($(2,1) + (230-\x * 360/14 :3)$){};
		\draw[fill=white] (\x) circle (.2) node{\x};
	}
	\end{scope}
	
	\end{tikzpicture}
\end{figure}
For each external vertex of $H$ with odd degree, pick one of the vertices associated to it in $Y$ and add an edge between them. The resulting graph $G$ is inner semi-Eulerian.

\begin{figure}[H]
	\centering
	\begin{tikzpicture}
	
	\draw (0,0) grid (1,1);
	\draw (3,0) grid (4,1);
	\draw (1,1) -- (3,1);
	\draw (2,1) -- (2,2);
	
	\draw (2,1) circle (3);
	
	\foreach \x in {1,2,...,14} {
		\node (\x) at ($(2,1) + (230-\x * 360/14 :3)$){};
	}

	\draw (1,1) -- (3.center);
	\draw (2,1) -- (12.center);
	\draw (3,1) -- (7.center);
	\draw (2,2) -- (5.center);

	\vertices{0}{0}{1}{1};
	\vertices{2}{1}{2}{2};
	\vertices{3}{0}{4}{1};
	
	\foreach \x in {1,2,...,14} {
		\setcounter{colsum}{\x};
	
		\ifnumodd{\value{colsum}}{\def\altcol{white}}{\def\altcol{black}}
		\draw[fill=\altcol] (\x) circle (.1);
	}
	\end{tikzpicture}
	\caption{An inner semi-Eulerian graph $G$ with inner subgraph $H$.}
\end{figure}
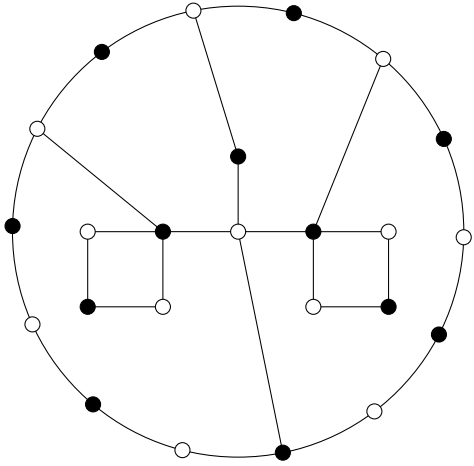

\end{example1}

\begin{remark}
    Often there is a more natural choice of outer completion for $H$ than the construction described above. In particular, many of our examples that are subgraphs of the square lattice use an outer completion that is also a subgraph of the square lattice.
\end{remark}

Given a billiard nest in $G$, we may construct an associated channel in $G'$ as follows. Let $B \in \bp(G)$ be a billiard nest. Define
$$\ch\colon \bp(G) \longrightarrow \cb(G') $$
by setting the vertices in $\ch(B)$ to be the internal black vertices $b$ of $G$ for which exactly half of the faces incident to $b$ are in $B$.

\begin{lemma}
	The map $\ch$ is a group homomorphism from $\bp(G)$ to $\cb(G')$.
\end{lemma}
\begin{proof}
	Let $B\in \bp(G)$. We first check that $\ch(B)$ is in fact a channel. Let $w$ be a white vertex of $G'$, i.e., an internal white vertex of $G$. Consider the edges incident to $w$. We wish to show that the number of these edges that contain a channel vertex is even. Let $f_1,...,f_n$ be the incident faces to $w$ in cyclic order. If $f_i$ and $f_{i+1}$ are both in $B$ or both not in $B$, then the edge between them does not contain a vertex in $\ch(B)$; otherwise, the edge contains a vertex in $\ch(B)$. Thus the number of incident edges to $w$ containing vertices from $\ch(B)$ is the number of times $f_i$ changes from being in $B$ to not being in $B$ or vice versa as we traverse the incident faces. Since after traversing all of the faces incident to $w$ we must arrive back at the starting face, we must change state an even number of times. Thus $\ch(B)$ is a channel.
	
	Now we check that $\ch$ preserves symmetric differences. Let $B_1, B_2 \in \bp(G)$. Consider the incident faces of $b$ for an arbitrary black vertex $b\in G'$. If $b\in \ch(B_1)$, then every other of these faces is in $B_1$. If $b\not\in \ch(B_1)$, then all or none of these faces are in $B_1$. The same holds for $B_2$. Then the claim follows since the symmetric difference of the incident faces to $b$ in $B_1$ with the incident faces to $b$ in $B_2$ is every second face if and only if exactly one of $B_1$ and $B_2$ contains every second face incident to $b$, which occurs if and only if $b$ is in exactly one of $\ch(B_1)$ and $\ch(B_2)$.
\end{proof}

We can now state our main result on billiard nests. 
Recall that the reduced dual graph of $G$ is the dual graph with the vertex corresponding to the external face of $G$ removed.

\begin{theorem}\label{bounce}
	Let $G$ be inner semi-Eulerian and let $G'$ be the inner subgraph of $G$. Assume that the outer subgraph of $G$ is a simple cycle.
	Then 
	$$|\bp(G)| = 2|\cb(G')|. $$
\end{theorem}
\begin{proof}
	 
	First we show that 
	$$|\ker \ch| = 2.$$
	Let $B\in \bp(G)$ be such that $\ch(B) = 0$. Then $B$ is locally constant, i.e., at each black vertex $b$, either all faces incident to $b$ are in $B$ or all faces incident to $b$ are not in $B$. In a connected component of the reduced dual graph of $G$, there is a path between any two faces. Any edge between two faces on that path contains a black vertex, which enforces the locally constant condition. Thus $B$ is constant along the entire path and therefore across any connected component of the reduced dual of $G$. So the claim follows if we can show that the reduced dual of $G$ is connected. Indeed, since $G$ is surrounded by a simple cycle, we may embed it such that it fills a convex region in the plane. Then a generic line connecting two internal faces will induce a path between them in the dual graph, which will avoid the external face by convexity. This proves that $B\in \bp(G)$ if and only if $B$ is constant, i.e., equal to either $F$ or $\varnothing$. Hence $|\ker\ch|=2$.
	
	
	To complete the proof, we show that $\ch$ is surjective.
	Let $C\in \cb(G')$ be a channel. We will construct a billiard nest $B$ such that $\ch(B)=C$. 
	Fix a spanning tree $T$ of the reduced dual of $G$ and some internal face $f_0$ of $G$. We decide whether an internal face $f$ is in $B$ as follows. If $f=f_0$, then $f\in B$. Otherwise, there is a unique path in $T$ from $f_0$ to $f$. If the path crosses an even number of edges in $G$ containing vertices in $C$, then $f\in B$. Otherwise, $f\not\in B$. 
	
	\textbf{Claim.} If $f_1,f_2$ are two adjacent internal faces of $G$, then $f_1$ and $f_2$ are both in $B$ or both not in $B$ if and only if the edge $e_1$ separating $f_1$ and $f_2$ does not contain a channel vertex. 
	
	If the dual edge $e_1^\vee$ is a part of $T$, then this follows by definition of $B$. Otherwise, adding $e_1^\vee$ to $T$ creates a simple planar cycle 
	$$Y=\{e_1^\vee,...,e_n^\vee\}\subseteq T\cup \{e_1^\vee\}.$$
	We wish to show that the number of channel vertices in $e_1$ has the same parity as the total number of channel vertices in $e_2,...,e_n$. We will be done if we can show that the number of edges in 
	$$Y^\vee = \{e_1,..., e_n\} $$
	that contain a channel vertex is even.
	
	Let $H$ be the induced subgraph of $G$ on the vertices in the interior of $Y$. Each $e_i$ connects a vertex $v_i \in H$ to a vertex $v_i'\in G\backslash H$. Set
	$$ \partial H = \{v_1, ..., v_n\} .$$
	
	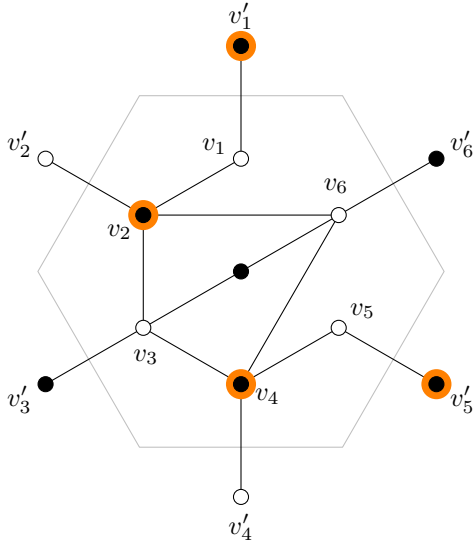
\begin{figure}[H]
	\centering
	\begin{tikzpicture}
	\newdimen\R
    \R=2.7cm
	\newdimen\Ro
    \Ro=3cm
    \newdimen\RR
    \RR=1.5cm
	\draw[thin,gray,opacity=.5] (0:\R) -- (60:\R) -- (120:\R) -- (180:\R) -- (240:\R)
	        -- (300:\R) -- cycle;
	\draw  (90:\RR) -- (150:\RR) -- (210:\RR) -- (270:\RR)
	        -- (330:\RR);
	        
	\foreach \x in {30, 90, ..., 330} {
	    \draw (\x:\RR) -- (\x:3);
	}
	\draw (30:\RR) -- (30:\RR) -- (150:\RR);
	\draw (30:\RR) -- (270:\RR);
	\draw (30:\RR) -- (210:\RR);
	
	\fill[orange] (90:\Ro) circle (.2) (150:\RR) circle (.2) (270:\RR) circle (.2) (330:\Ro) circle (.2);
	
	\draw[fill=black] (150:\RR) circle (0.1) (270:\RR) circle (.1) (30:\Ro) circle (.1)
	(90:\Ro) circle (.1) (210:\Ro) circle (.1) (330:\Ro) circle (.1) (0:0) circle (.1);
	
	\draw[fill=white] (150:\Ro) circle (0.1) (270:\Ro) circle (.1) (30:\RR) circle (.1)
	(90:\RR) circle (.1) (210:\RR) circle (.1) (330:\RR) circle (.1);
	
	\foreach \x in {1, 2, ..., 6} {
	    \node at (30+ \x * 60: 3.4) {$v'_{\x}$};
	}
	
	\foreach \x in {1, 2, ..., 6} {
	    \node at (42 + \x * 60: 1.7) {$v_{\x}$};
	}
	\end{tikzpicture}
	\caption{The graph $G$ near a dual graph cycle $Y$, shown as a gray hexagon. The vertices in the interior of $Y$ form the subgraph $H$. The channel $C$ is shown in orange.}
	\end{figure}
	
	Since $H$ is bipartite and $C$ uses only black vertices,
	\begin{equation}\label{bipeq}
	\sum_{b\in C\cap H} \deg_H b = \sum_{w \in H} \chdeg_H w ,
	\end{equation}
	where the sum on the left is over black vertices, the sum on the right is over white vertices, and $\chdeg_H w$ is the size of the neighborhood of $w$ in $C\cap H$. 
	
	Note that for $b\in H\backslash\partial H$ we have
	$$\deg_H b = \deg_G b \equiv_2 0 $$
	since $G$ is inner semi-Eulerian and all vertices of $H$ are internal to $G$. (Here we write $n\equiv_2 m$ to mean that $n-m$ is even.) Thus
	$$\sum_{b\in C\cap H} \deg_H b \equiv_2 \sum_{b\in C\cap \partial H} \deg_H b
	\equiv_2 \#\{e_i\mid v_i \text{ is in $C$} \}. $$
	The last equivalence follows since 
	$$\deg_H b = \deg_G b - \#\{e_i\mid b\in e_i \}
	\equiv_2 \#\{e_i\mid b \in e_i \}. $$
	Working now with the other side of (\ref{bipeq}), for $w\in H\backslash \partial H$ we have that
	$$\chdeg_H w = \chdeg_G w = \chdeg_{G'}\equiv_2 0 $$
	by definition of a channel. Thus
	$$\sum_{w\in H} \chdeg_H w \equiv_2 \sum_{w\in \partial H} \chdeg_H w
	\equiv_2 \#\{e_i\mid v_i' \text{ is in $C$}\}. $$
	The last equivalence here follows since 
	$$\chdeg_H w = \chdeg_G w - \#\{e_i\mid w\in e_i \text{ and } v_i'\in C \}
	\equiv_2 \#\{e_i\mid w\in e_i \text{ and } v_i'\in C \}. $$
	
	Substituting our results into (\ref{bipeq}), we find
	$$\#\{e_i\mid v_i \text{ is in $C$} \}\equiv_2\#\{e_i\mid v_i' \text{ is in $C$}\}.$$
	Thus
	$$ \#\{e_i \mid \text{$e_i$ contains a vertex from $C$}\}
	=\#\{e_i\mid v_i \text{ is in $C$} \}+\#\{e_i\mid v_i' \text{ is in $C$}\}
	\equiv_2 0$$
	as desired, proving the claim.
	
	We may now verify that $B$ is indeed a billiard nest such that $\ch(B)=C$. 
	Let $b\in G$ be a black vertex. If $b\in C$, then every edge incident to $b$ contains a channel vertex. Thus, by the preceding claim, the faces incident to $b$ must alternate between being in $B$ and not being in $B$. If $b\not\in C$, then every edge incident to $b$ does not contain a channel vertex. Thus, by the claim, the internal faces incident to $b$ are either all in $B$ or all not in $B$. In particular, external black vertices are not in $C$ (since they are not in $G'$), so $B$ meets the conditions for being a billiard nest. Furthermore, it is clear from this description that $\ch(B)= C$. Thus $\ch$ is a surjective homomorphism
	$$\bp(G) \longrightarrow \cb(G') .$$
	
	Combining this with the size of the kernel computed earlier, we find that
	\[|\bp(G)| = 2 |\cb(G')|. \qedhere\]
	
\end{proof}

\begin{corollary}
	Let $G$ be inner semi-Eulerian and let $G'$ be the inner subgraph of $G$. Assume that the outer subgraph of $G$ is a simple cycle and that $G$ has the same number of black and white vertices.
	Then 
	$$\frac{1}{2}|\bp(G)| \text{ divides the number of matchings of $G'$}. $$
	Moreover, $m_{G'}$ is odd if and only if $G$ has exactly one nonempty billiard nest.
\end{corollary}
\begin{proof}
    By Theorem~\ref{bounce}, we have
    \[ \frac{1}{2}|\bp(G)| = |\cb(G')|. \]
    By channeling 2s (Theorem \ref{twos}), we know that $|\cb(G')|$ divides the number of matchings of $G'$. By Proposition~\ref{algproof}, $m_{G'}$ is odd if and only if $\dim \mc(G') = 2 \dim \cb(G') = 0$, that is, if and only if $|\cb(G')|=1$. The claim follows.
\end{proof}

Thus our study of channels in appropriate graphs $H$ (in particular, by Proposition \ref{outercomp}, all inner semi-Eulerian graphs) reduces to the study of billiard nests in an outer completion $G$. Billiard nests are considerably easier to work with since every face of $G$ is contained in a unique billiard path in the path basis of $G$. In general there is no such basis for the channels of $H$; vertices of $H$ may be contained in no channel and there may be no channel basis for $H$ with pairwise disjoint elements. However, any billiard path can be found by starting with a face of $G$ and adding additional faces as required by the definition of billiard nests. 

\subsection{Arithmetic billiards}\label{arithmetic}
As an application, let us find the billiard paths for the rectangle grid graph $\mathcal R_{m+1\times n+1}$, an outer completion of $\mathcal R_{m-1\times n-1}$. For this graph, we may use our interpretation of billiard paths as the paths through black vertices traced out by billiard balls travelling at 45 degree angles. (Explicitly, a face in the billiard path is interpreted as a billiard ball traveling between the black vertices in that face.)  We begin by straightening out the billiard paths; to do so, we tile the plane with copies of our rectangle.

\begin{figure}[H]
	\centering
	\begin{tikzpicture}
	\draw (0,0) rectangle (3,2);
	
	\draw[thick, red] (0,0) -- (2,2) -- (3,1)
		-- (2,0) -- (0,2);
		
	\draw[fill, blue!50!brown] (0,0) circle (.1);
	\draw (0,0) circle (.1);
	\draw[fill, orange!50!brown] (3,0) circle (.1);
	\draw (3,0) circle (.1);
	\draw[fill, green!50!brown] (0,2) circle (.1);
	\draw (0,2) circle (.1);
	\draw[fill, yellow!50!brown] (3,2) circle (.1);
	\draw (3,2) circle (.1);
	
	\begin{scope}[shift={(6,-2)}]
	\draw (-1.5, 3) node[scale=2]{$\iff$};
	\begin{scope}[xscale=3,yscale=2]
	\draw (0,0) grid (2,3);
	\end{scope}
	\draw[thick, red] (0,0) -- (3*2,2*3);
	\foreach \x in {0,2} {
	\foreach \y in {0,2} {
	\draw[fill, blue!50!brown] (3*\x,2*\y) circle (.1);
	\draw (3*\x,2*\y) circle (.1);
	}}
	\foreach \x in {1} {
		\foreach \y in {0,2} {
			\draw[fill, orange!50!brown] (3*\x,2*\y) circle (.1);
			\draw (3*\x,2*\y) circle (.1);
	}}
	\foreach \x in {0,2} {
		\foreach \y in {1,3} {
			\draw[fill, green!50!brown] (3*\x,2*\y) circle (.1);
			\draw (3*\x,2*\y) circle (.1);
	}}
	\foreach \x in {1} {
		\foreach \y in {1,3} {
			\draw[fill, yellow!50!brown] (3*\x,2*\y) circle (.1);
			\draw (3*\x,2*\y) circle (.1);
	}}
	\end{scope}
	\end{tikzpicture}
	\caption{A rectangle used to tile the plane. The coloring is provided as a visual indicator of the rectangle's orientation. Any billiard path on the left corresponds to a line of slope 1 on the right and vice versa.}
\end{figure}
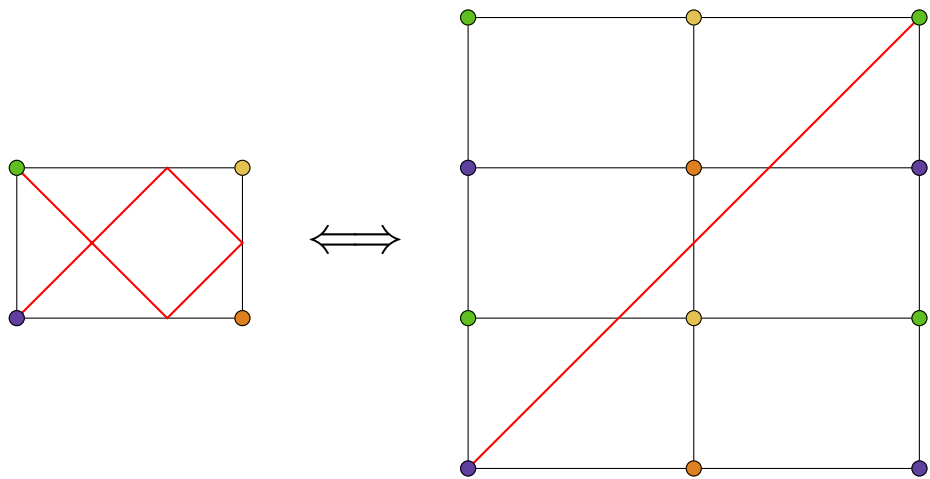
We can then lift the billiard path to a straight line of slope 1 in the tessellation. A billiard path between two corners of $\mathcal R_{m+1\times n+1}$ will be the diagonal of a square in the tessellation. Any such square must have a side length divisible by $m$ and $n$ (the side lengths of $\mathcal R_{m+1\times n+1}$). The minimal square in the tessellation with corners from $\mathcal R_{m+1\times n+1}$ then has side length given by the least common multiple of $m$ and $n$. To determine the number of internal faces of $\mathcal R_{m+1\times n+1}$ through which the path travels, we may count the number of unit squares through which the path travels in the tessellation. Since the straightened billiard path travels along the diagonal of a square with side length $\lcm(m,n)$, this path travels through $\lcm(m,n)$ unit squares. 

For now assume that at least one of $m+1$ and $n+1$ is even. Then exactly two corners of $\mathcal R_{m+1\times n+1}$ are black. Any two distinct billiard paths pass through distinct internal faces of $\mathcal R_{m+1\times n+1}$. We shall count path basis elements by counting the internal faces through which they pass. There is one billiard path through the black corners. From the last paragraph, we know this path uses $\lcm(m,n)$ internal faces. Now, every other path on black vertices uses \emph{twice} as many internal faces. Indeed, since the other paths do not pass through a corner, they must end on their starting point. To reach their starting point in the tessellation, the paths must lift to the diagonal of a square of side length $2\lcm(m,n)$ since one of $m$ and $n$ is odd. Since every internal face is part of a unique path basis element, we may now count the billiard paths for $\mathcal R_{m+1\times n+1}$.

\begin{theorem}
	The rectangle grid graph $\mathcal R_{m+1\times n+1}$ with $(m+1)(n+1)$ even has a path basis of size
	$$\frac{\gcd(m,n)+1}{2} .$$
\end{theorem}
\begin{proof}
	There are $mn$ total internal faces in $\mathcal R_{m+1\times n+1}$. From the above, $mn-\lcm(m,n)$ of these are part of a billiard path not passing through a corner. Each such basis element uses $2\lcm(m,n)$ internal faces. Thus there are
	$$\frac{mn - \lcm(m,n)}{2\lcm(m,n)}
		= \frac{\gcd(m,n)-1}{2}$$
	non-corner billiard paths. Adding back the last billiard path gives the claim.
\end{proof}

\begin{corollary}
	The rectangle grid graph $\mathcal R_{m-1\times n-1}$ with $(m-1)(n-1)$ even has
	$$|\cb(\mathcal R_{m-1\times n-1})| = 2^{\frac{\gcd(m,n)-1}{2}}. $$
\end{corollary}
\begin{proof}
	The inner subgraph of $\mathcal R_{m+1\times n+1}$ is $\mathcal R_{m-1 \times n-1}$, and the reduced dual graph of $\mathcal R_{m\times n}$ is connected. Thus the result follows by Theorem \ref{bounce}.
\end{proof}

\begin{corollary}\label{corarithmetic}
	The number of matchings of $\mathcal R_{m-1\times n-1}$ is divisible by
	$$2^{\frac{\gcd(m,n)-1}{2}}. $$
	Furthermore, the number of matchings is odd if and only if $m$ and $n$ are coprime.
\end{corollary}
\begin{proof}
	If $m-1$ and $n-1$ are both odd, then $\m=0$ and the claim follows. Otherwise, the hypothesis of the previous corollary holds and we may channel 2s (apply Theorem \ref{twos}) to get the divisibility statement. The last claim follows from Proposition \ref{algproof}.
\end{proof}

In the next section, we apply the geometric interpretation of billiard paths for subgraphs of the square lattice to construct an algorithm for finding a path basis for such graphs. 

\subsection{Finding billiard paths in the square lattice}
Let $G$ be a subgraph of the square lattice $\Z\times \Z$ formed by taking a simple cycle and adding all edges and vertices in its interior. Assume we are given the boundary cycle (for instance, as an ordered list of edges).  Let $P$ be the set of exterior black vertices in $G$. From the input data we may determine in $O(|P|)$ time the set of internal faces in $G$ which are incident to the boundary. Our algorithm for constructing a path basis will have complexity $O(|P|\log |P|)$. We shall construct an auxiliary graph $A$ with vertex set $P$ which will have connected components corresponding to billiard paths in $G$. To begin, given the coordinates $(x,y)$ of a vertex $b\in P$ we compute two indices:
\begin{align*}
b^+ &= y-x, \\
b^- &= y+x,
\end{align*}
the positive and negative index, respectively. The positive index indicates which line of slope 1 the vertex lies on, and the negative index indicates which line of slope $-1$ the vertex lies on. Now sort the pairs 
$$\{(b^+, b^-) \mid b\in P\}$$
lexicographically. Let $b_1$ and $b_2$ be two consecutive vertices in this list. If the vertex $b_1$ is incident to an internal face to its upper right, then put an edge between $b_1$ and $b_2$ in $A$. Iterate through the list and do this for each consecutive pair. 

Next, sort the pairs
$$\{(b^-, b^+) \mid b\in P\}$$
lexicographically. Again, consider consecutive vertices $b_1$ and $b_2$ in this list. If $b_1$ is incident to an internal face to its upper left, then put an edge between $b_1$ and $b_2$ in $A$. Do this for all such consecutive pairs. 

After this is complete, the connected components of $A$ will be in correspondence with billiard paths of $\bp(G)$. 
Specifically, a connected component of $A$ corresponds to the minimal nest containing any face of $G$ incident to a vertex in that component.
If the number of connected components in $A$ is $d$, then there are $2^d$ billiard nests in $G$, and if the reduced dual of $G$ is connected, then there are $2^{d-1}$ channels on the black vertices of the inner subgraph of $G$.

Note that computing the indexes, adding edges to $A$, and finding the connected components of $A$ takes $O(|P|)$ time, while the sorts take $O(|P|\log |P|)$ time. Thus the algorithm runs in almost linear time.

\begin{example1}
	Let $G$ be the graph shown in Figure \ref{algex}. Set $A$ to be the graph on the external black vertices of $G$, with no edges.
	
	\begin{figure}[H]
		\centering
	\begin{tikzpicture}
	\draw[ black] (0,0) grid (3,5);
	\draw (0,0) grid (5,3);
	
	\draw[red] (0,0) -- (3,3) -- (1,5) -- (0,4) -- (4,0)--(5,1)--(3,3);
	
	
	\vertices{0}{0}{3}{5}
	\vertices{4}{0}{5}{3}

	\node at (-1,0){$(0,0)$};
	\node at (-1,2){$(2,2)$};
	\node at (-1,4){$(4,4)$};
	\node at (1,5.5){$(4,6)$};
	\node at (2,-.5){$(-2,2)$};
	\node at (3,5.5){$(2,8)$};
	\node at (4,-.5){$(-4,4)$};
	\node at (6,1){$(-4,6)$};
	\node at (6,3){$(-2,8)$};
	\node at (3.5,3.5){$(0,6)$};
	\end{tikzpicture}
	\caption{The graph $G$ with external black vertices labeled by their indexes in the form $(b^+, b^-)$. One of the billiard paths is shown in red.} 
	\label{algex}
	\end{figure}
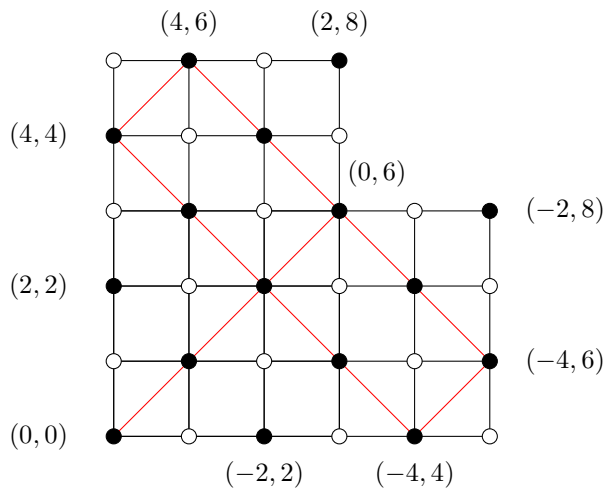
	We list out the index pairs in lexicographical order:
	
	\begin{table}[H]
	\centering
	\addtolength{\tabcolsep}{-3pt}
	\begin{tabular}{|c||c|c|c|c|c|c|c|c|c|c|}
		\hline 
		$(b^+, b^-)$ & $(-4, 4)$& $(-4, 6)$&$(-2,2)$&$(-2,8)$&$(0,0)$&$(0,6)$
		&$(2,2)$&$(2,8)$&$(4,4)$&$(4,6)$  \\ 
		\hline 
		upper right face?& y & n &y &n &y &n &y &n & y& n \\ 
		\hline 
	\end{tabular} 
    \end{table}
We connect each vertex with a $y$ to the next vertex on the list. The following figure shows $A$ after this step is completed.

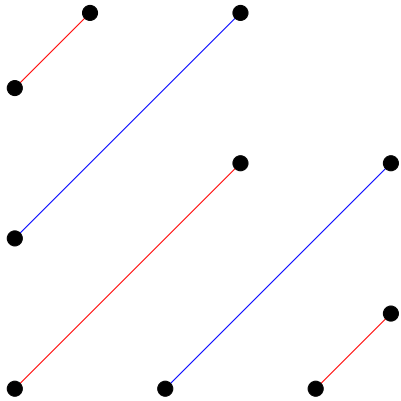
\begin{figure}[H]
	\centering
\begin{tikzpicture}
\draw[red] (0,4)--(1,5);
\draw[blue] (0,2)--(3,5);
\draw[red] (0,0)--(3,3);
\draw[blue] (2,0)--(5,3);
\draw[red] (4,0)--(5,1);

\pgfmathsetmacro{\radius}{0.1}
\draw[fill = black] (0,0) circle (\radius);
\draw[fill = black] (0,2) circle (\radius);
\draw[fill = black] (0,4) circle (\radius);
\draw[fill = black] (1,5) circle (\radius);
\draw[fill = black] (3,5) circle (\radius);
\draw[fill = black] (3,3) circle (\radius);
\draw[fill = black] (5,3) circle (\radius);
\draw[fill = black] (5,1) circle (\radius);
\draw[fill = black] (4,0) circle (\radius);
\draw[fill = black] (2,0) circle (\radius);
\end{tikzpicture}

\end{figure}

The colors here mean nothing at the moment, but once we finish adding edges they will indicate the connected components of $A$. For the second phase, we reorder the index pairs, this time lexicographically based on $(b^-,b^+)$.

\begin{table}[H]
	\centering
	\addtolength{\tabcolsep}{-2.7pt}
	\begin{tabular}{|c||c|c|c|c|c|c|c|c|c|c|}
		\hline 
		$(b^-, b^+)$ & $(0,0)$& $(2,-2)$&$(2,2)$&$(4,-4)$&$(4,4)$&$(6,-4)$
		&$(6,0)$&$(6,4)$&$(8,-2)$&$(8,2)$  \\ 
		\hline 
		upper left face?& n & y &n &y &n &y &y &n & n& n \\ 
		\hline 
	\end{tabular} 
\end{table}

Again we connect each vertex with a y to the following vertex to complete the construction of $A$. The result is shown below, with connected components displayed in different colors.

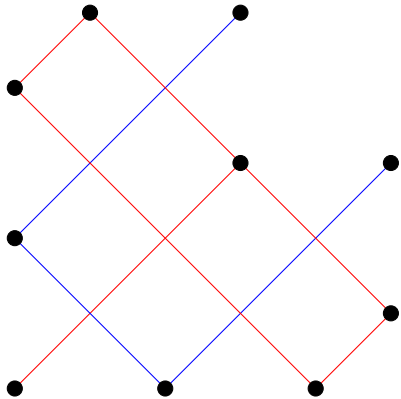
\begin{figure}[H]
	\centering
	\begin{tikzpicture}
	\draw[red] (0,4)--(1,5);
	\draw[blue] (0,2)--(3,5);
	\draw[red] (0,0)--(3,3);
	\draw[blue] (2,0)--(5,3);
	\draw[red] (4,0)--(5,1);
	\draw[red] (4,0)--(0,4);
	\draw[red] (1,5)--(5,1);
	\draw[blue] (2,0)--(0,2);

	\pgfmathsetmacro{\radius}{0.1}
	\draw[fill = black] (0,0) circle (\radius);
	\draw[fill = black] (0,2) circle (\radius);
	\draw[fill = black] (0,4) circle (\radius);
	\draw[fill = black] (1,5) circle (\radius);
	\draw[fill = black] (3,5) circle (\radius);
	\draw[fill = black] (3,3) circle (\radius);
	\draw[fill = black] (5,3) circle (\radius);
	\draw[fill = black] (5,1) circle (\radius);
	\draw[fill = black] (4,0) circle (\radius);
	\draw[fill = black] (2,0) circle (\radius);
	\end{tikzpicture}
	
\end{figure}
There are two connected components of $A$, corresponding to the two billiard paths in $\bp(G)$. Thus there are $2^{2-1}=2$ channels on the black vertices of the inner subgraph $G'$ of $G$. 
\end{example1}
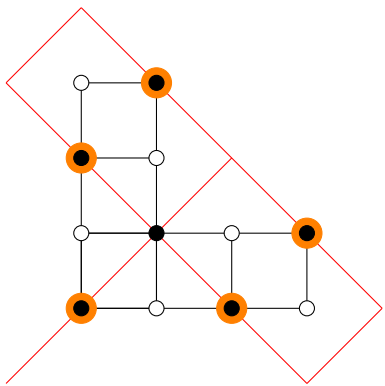
\begin{figure}[H]
	\centering
	\begin{tikzpicture}
	\draw[ black] (0,0) grid (1,3);
	\draw (0,0) grid (3,1);
	
	\draw[red] (-1,3)--(0,4);
	\draw[red] (-1,-1)--(2,2);
	\draw[red] (3,-1)--(4,0);
	\draw[red] (3,-1)--(-1,3);
	\draw[red] (0,4)--(4,0);
	
	\pgfmathsetmacro{\chr}{0.2}
	\draw[fill,orange] (0,0) circle (\chr);
	\draw[fill,orange] (0,2)  circle (\chr);
	\draw[fill,orange] (1,3)  circle (\chr);
	\draw[fill,orange] (2,0)  circle (\chr);
	\draw[fill,orange] (3,1) circle (\chr);

	\vertices{0}{0}{1}{3}
	\vertices{2}{0}{3}{1}
	\end{tikzpicture}
	\caption{The inner subgraph $G'$ formed by deleting the external vertices of $G$, with its nonzero channel highlighted. Either of the billiard nests shown in the previous figure constructs this channel.}
	\end{figure}

\begin{remark}
    Let $G$ be a subgraph of the square lattice formed by taking a simple cycle and adding all edges and vertices in its interior. A variant of the algorithm above may be used to compute the path basis for an outer completion of $G$ (rather than for $G$ itself, as done above). The path basis for an outer completion controls the channels (and therefore parity information) in $G$, rather than in the inner subgraph $G'$. This alteration then lets us understand the channels of graphs that do not have an outer completion fitting in the square lattice. 
    
    To compute the path basis for an outer completion $H$ of $G$, with boundary cycle $\partial H$, make the following changes:
    \begin{itemize}
        \item The auxiliary graph $A$ now has vertex set given by the internal faces of $H$ that are incident to vertices in $P$ (the exterior black vertices of $G$).
        \item Start off by connecting faces in $A$ that are incident to the same black vertex in $\partial H$, as well as every second face surrounding a vertex in $P$.
        \item Instead of connecting the vertices of $P$, we connect the appropriate faces in $A$. 
        \begin{itemize}
        \item In the $(b^+,b^-)$ table, if a vertex $v\in P$ has a $y$ label, then connect the upper right face incident to $v$ to the lower left face incident to the next vertex in the table.
        
        \item In the $(b^-,b^+)$ table, if a vertex $v\in P$ has a $y$ label, then connect the upper left face incident to $v$ to the lower right face incident to the next vertex in the table.
        \end{itemize}
    \end{itemize}
    The connected components of $A$ then are equinumerous with the path basis for $H$. Since an outer completion for $G$ can be constructed in $O(|P|)$ time and can be chosen to add at most $|P|$ internal faces, the altered algorithm remains $O(|P|\log |P|)$.
\end{remark}
\section{Combinatorial arguments}\label{seccomb}

In this section we give a combinatorial proof that existence of a nonempty channel in $G$ is necessary and sufficient for $\m$ to be even (Proposition \ref{algproof}). In the course of this proof, we develop conditions under which deleting two adjacent vertices (a ``vertex pair") results in a graph with the same number of channels. This is called \emph{channel routing} and is described in Lemma \ref{routing}; we apply it to deduce a result about rectangle grid graphs in addition to the main theorem of the section. 

In Section \ref{secprelim} we described Lov\'asz's proof of Proposition \ref{algproof}. That argument was of a surprisingly different nature from the proof of Theorem \ref{twos}, which relied upon a Kasteleyn signing and matrix normal forms and which could fail for graphs without a Kasteleyn signing (cf. Remark \ref{kastsignnec}). Like Lov\'asz's argument, the combinatorial results in the remainder of this section also do not need a Kasteleyn signing. We do, however, require that our graph be bipartite. We prove the forward and backward directions of Proposition \ref{algproof} separately in the two subsections.

\subsection{Existence of a channel implies $m_G$ even}
Our approach for this direction will be to construct a fixed-point free involution on the set of matchings $\mat (G)$, which will imply that $|\mat(G)|=m_G$ is even. To construct such an involution, we employ a technique called \emph{cycle flipping}. This is a commonly used method to build involutions on perfect matchings and show 2-divisibility results. See for example \cite{Ciucu1996} or \cite{Pachter1997} to see this applied to graphs with reflective symmetry, or \cite{Kuo2004} for disjoint unions of two graphs.

Given a perfect matching $\mu\in\mat(G)$, the idea is to find a cycle $Y$ of edges in the graph such that every second edge in the cycle is in $\mu$. We may then construct a new edge set $\mu'$ by replacing the edges of $\mu \cap Y$ with the edges of $Y-\mu$. Since each vertex in $Y$ is contained in exactly one edge in either case, $\mu'$ is also a perfect matching.

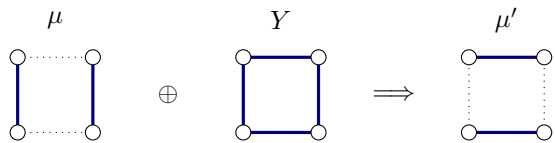
\begin{figure}[htb]
    \centering
    \begin{tikzpicture}
        \pgfmathsetmacro{\radius}{0.1}
        \pgfmathsetmacro{\x}{0}
        \draw (0.5,1.5) node{$\mu$};
        \draw[dotted] (\x,0) -- ++(1,0) (\x,1) -- ++(1,0);
        \draw[blue!50!black,very thick] (\x,0) -- ++(0,1) (\x+1,0) -- ++(0,1);
        \draw[fill=white] (\x,0) circle (\radius);
        \draw[fill=white] (\x,1) circle(\radius);
        \draw[fill=white] (\x+1,1) circle(\radius);
        \draw[fill=white] (\x+1,0) circle(\radius);
        \draw (2,0.5) node{$\oplus$};
        \pgfmathsetmacro{\x}{3}
        \draw (\x+0.5,1.5) node{$Y$};
        \draw[blue!50!black,very thick] (\x,0) -- ++(1,0) (\x,1) -- ++(1,0);
        \draw[blue!50!black,very thick] (\x,0) -- ++(0,1) (\x+1,0) -- ++(0,1);
        \draw[fill=white] (\x,0) circle (\radius);
        \draw[fill=white] (\x,1) circle(\radius);
        \draw[fill=white] (\x+1,1) circle(\radius);
        \draw[fill=white] (\x+1,0) circle(\radius);
        \draw (5,0.5) node{$\implies$};
        \pgfmathsetmacro{\x}{6}
        \draw (\x+0.5,1.5) node{$\mu'$};
        \draw[blue!50!black,very thick] (\x,0) -- ++(1,0) (\x,1) -- ++(1,0);
        \draw[dotted] (\x,0) -- ++(0,1) (\x+1,0) -- ++(0,1);
        \draw[fill=white] (\x,0) circle (\radius);
        \draw[fill=white] (\x,1) circle(\radius);
        \draw[fill=white] (\x+1,1) circle(\radius);
        \draw[fill=white] (\x+1,0) circle(\radius);
    \end{tikzpicture}
    \caption{Flipping the cycle $Y$ from $\mu$ to $\mu'$.}
    \label{fig:sum}
\end{figure}

For the map this produces to be an involution, the same cycle has to be identified for both $\mu$ and $\mu'$. We will produce this cycle using the following tool.

\begin{definition}
A \textit{pairing function} across a vertex set $C\subseteq V$ is a collection of involutions $f_v$ associated to each $v\in V$ which act on the edge set
$$\{\edge{v}{v'}\in E \mid v'\in C\}$$
such that $f_v\circ f_v = \text{id}$ and $f_v$ has no fixed points.
\end{definition}

The existence of a pairing function across $C$ is a combinatorial realization of the statement that, for each vertex $v$, the neighborhood $N(v)$ contains an even number of elements from $C$. Applying this to the definition of a channel gives the following lemma.

\begin{lemma}
Let $G$ be a graph and let $C \subseteq V$ be any vertex set. Then $C$ is a channel if and only if there exists a pairing function across $C$.
\end{lemma}

We will use a pairing function to trace out a walk such that every second edge lies in our matching. Finiteness of the graph will force this path to eventually enter a cycle with the properties we require.

\begin{theorem}\label{cycle}
	Let $G$ be a bipartite graph with a nonempty channel $C\in \mathcal C(G)$. Then to each matching $\mu$ of $G$, one can assign a nonempty set of edges $S(\mu) \subseteq E$ satisfying the following properties.
	
	\begin{enumerate}[label=(\roman*)]
		\item $S(\mu)$ is a simple cycle of even length.
		
		\item Every second edge of $S(\mu)$ lies in $\mu$.
		
		\item $S(\mu)$ depends only on the edges of $\mu$ containing a vertex from $C$. 
		
		\item If $\mu'$ is a matching satisfying $S(\mu) = \mu \oplus \mu'$ then $S(\mu')=S(\mu)$.
	\end{enumerate}
\end{theorem} 
\begin{proof}
	Fix a vertex $v_0\in C$ and a pairing function $f_{\bullet}$ across $C$.  We define a walk \[e_0 = \edge{v_0}{v_1}, \quad e_1 = \edge{v_1}{v_2}, \quad e_2 = \edge{v_2}{v_3}, \quad \dots\] as follows.
	If $n$ is even, take $e_n$ to be the unique edge in $\mu$ incident to $v_n$. If $n$ is odd, take $e_n=f_{v_n}(e_{n-1})$. Since $f_\bullet$ is a pairing function across $C$, $v_{n-1} \in C$ implies $v_{n+1} \in C$. (An example of this walk is shown in Figure \ref{fig:pairfunction}.)
	
	Because $G$ is finite, there exist $n_0$ and $p > 0$ such that $v_{n_0+p}=v_{n_0}$. Take $p$ to be minimal and $n_0$ to be minimal for that choice of $p$. Note that $p$ must be even since $G$ is bipartite. Also $n_0$ must be even: if $n_0$ were odd, then both $e_{n_0-1}$ and $e_{n_0+p-1}$ would be edges in $\mu$ containing $v_{n_0+p} = v_{n_0}$. But this would imply $e_{n_0-1} = e_{n_0+p-1}$ and hence $v_{n_0-1}=v_{n_0+p-1}$, which would contradict minimality of $n_0$.

	We claim that
	\begin{equation*}
	S(\mu) = \left\{e_n \mid n_0\leq n < n_0 + p \right\}
	\end{equation*}
	satisfies the desired properties. Property (i) follows by construction as well as minimality and evenness of $p$; (ii) is clear from construction. For (iii), the only dependence of the walk on $\mu$ involves edges incident to $v_n$ for $n$ even. Since $v_n \in C$ for $n$ even, the cycles produced will be the same for matchings that only differ away from $C$.
	
	\begin{figure}
	\centering
	\begin{tikzpicture}[scale=1.4]
	\begin{scope}[opacity=.4]
	\draw[dotted,gray] (-3,0) grid (5,3);
	\begin{scope}[shift={(-3,0)}]
	\foreach \x in {0} {
		\draw[orange,fill] (7-\x,\x) circle (.2);
		\draw[orange,fill] (8-\x,1+\x) circle (.2);
		\draw[orange,fill] (\x,1+\x) circle (.2);
		\draw[orange,fill] (\x+1,\x) circle (.2);
	}

	\draw[orange,fill] (2,1+2) circle (.2);
	\draw[orange,fill] (2+1,2) circle (.2);
	\draw[orange,fill] (0,3)circle(.2) (3,0)circle(.2);
	\end{scope}
	\vertices{-3}{0}{5}{3};
	\draw[very thick, blue!50!black] (0,1) -- (0,2);
	\draw[very thick, blue!50!black] (1,1) -- (1,2);
	\draw[very thick, blue!50!black] (0,3) -- (1,3);
	\draw[very thick, blue!50!black] (4,1) -- (4,2);
	\draw[very thick, blue!50!black] (-3,0) -- (-3,1);
	\draw[very thick, blue!50!black] (-2,0) -- (-1,0);
	\draw[very thick, blue!50!black] (-2,1) -- (-1,1);
	\draw[very thick, blue!50!black] (-1,2) -- (-1,3);
	\draw[very thick, blue!50!black] (-3,2) -- (-2,2);
	\draw[very thick, blue!50!black] (-3,3) -- (-2,3);
	\draw[very thick, blue!50!black] (3,2) -- (3,1);
	
	\end{scope}
	\begin{scope}[shift={(-3,0)}]
	\foreach \x in {2} {
		\draw[orange,fill] (7-\x,\x) circle (.2);
		\draw[orange,fill] (8-\x,1+\x) circle (.2);
	}

	\draw[orange,fill] (8,1+0) circle (.2);
	\draw[orange,fill] (7,0) circle (.2);
	\draw[orange,fill] (3,0)circle(.2);
	\draw[orange,fill] (8,3)circle(.2) (5,0)circle(.2);
	\end{scope}
	
	\draw (0,0) -- (5,0) -- (5,3) -- (2,3) -- (2,0);
	\draw[very thick, blue!50!black] (0,0) -- (1,0);
	\draw[very thick, blue!50!black] (2,0) -- (3,0);
	\draw[very thick, blue!50!black] (4,0) -- (5,0);
	\draw[very thick, blue!50!black] (5,1) -- (5,2);
	\draw[very thick, blue!50!black] (5,3) -- (4,3);
	\draw[very thick, blue!50!black] (2,3) -- (3,3);
	\draw[very thick, blue!50!black] (2,1) -- (2,2);
	\vertices{0}{0}{5}{0};
	\vertices{5}{0}{5}{3};
	\vertices{2}{3}{4}{3};
	\vertices{2}{1}{2}{2};
	\foreach \x in {0,1,...,5} {
		\node[label=below:{$v_{\x}$}] (\x) at (\x, 0) {};
	} 
	\foreach \x in {0,1,...,4} {
		\node[label=below:{$e_{\x}$}] (e\x) at (.5+\x, 0.1) {};
	} 
	\foreach \x in {6,7} {
		\node[label=right:{$v_{\x}$}] (\x) at (5, \x-5) {};
	} 
	\foreach \x in {8,9,10,11}{
	\draw (5-\x+8,3) node[label=above:{$v_{\x}$}] {};
	}
	\draw (2,2)node[label=left:{$v_{12}$}]{} 
		  (2,1)node[label=left:{$v_{13}$}]{};
	\foreach \x in {5,6,7} {
		\node[label=right:{$e_{\x}$}] (\x) at (4.9, \x-4.5) {};
	} 
	\foreach \x in {8,9,10}{
		\draw (4.5-\x+8,2.9) node[label=above:{$e_{\x}$}] {};
	}
	\foreach \x in {11,12,13} {
		\node[label=left:{$e_{\x}$}] at (2.1, 13.5-\x) {};
	} 
	\draw (2,2)node[label=left:{$v_{12}$}]{} 
	(2,1)node[label=left:{$v_{13}$}]{};

	\begin{scope}[green, thick]
	\foreach \x in {1,3} {
		\draw (\x-.3,0) to[out=90,in=90] (\x+.3,0);
	}
	\draw (5-.3,0) to[out=90,in=180] (5,.3);
	\draw (5,2-.3) to[out=180,in=180] (5,2+.3);
	\draw (4.3,3) to[out=-90,in=-90] (4-.3,3);
	\draw (2.3,3) to[out=-90,in=0] (2,3-.3);
	\draw (2,1.3) to[out=0,in=0] (2,1-.3);
	\end{scope}
	\end{tikzpicture}
	\caption{An illustration of the construction used in the proof of Theorem \ref{cycle}. The matching $\mu$ is shown in purple, and edges which are swapped by a relevant pairing functions are indicated by green lines. In this case, Theorem \ref{cycle} produces the cycle $S(\mu)=\{e_2,...,e_{13}\}$.}
	\label{fig:pairfunction}
\end{figure}
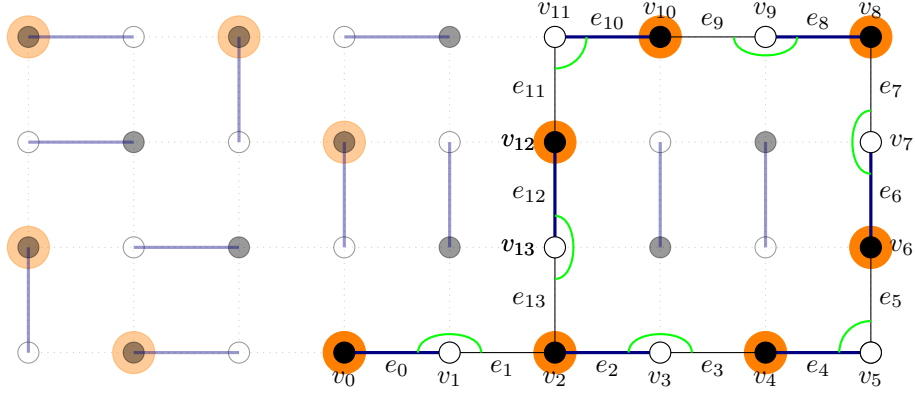

		For (iv), let $S(\mu) = \mu \oplus \mu'$, so that
		\[\mu' = (\mu \setminus C) \cup (C \setminus \mu) = (\mu \setminus C) \cup \{e_{n_0+1}, e_{n_0+3}, \dots, e_{n_0+p-1}\}.\]
		Denote the walk constructed for $\mu'$ by
		\[e'_0 = \edge{v'_0}{v'_1}, \quad e'_1 = \edge{v'_1}{v'_2}, \quad e'_2 = \edge{v'_2}{v'_3}, \quad \dots.\]
		Since $\mu$ and $\mu'$ agree outside of $S(\mu)$, $e_n = e'_n$ for $n < n_0$.
		
		By definition, $e'_{n_0}$ is the unique edge of $\mu'$ containing $v'_{n_0} = v_{n_0}=v_{n_0+p}$, so $e'_{n_0} = e_{n_0+p-1}$ and hence $v'_{n_0+1} = v_{n_0+p-1}$. Then
		\[e'_{n_0+1} = f_{v'_{n_0+1}}(e'_{n_0}) = f_{v_{n_0+p-1}}(e_{n_0+p-1}) = e_{n_0+p-2}\]
		and hence $v'_{n_0+2} = v_{n_0+p-2}$. Repeating this argument, we find that $e'_{n_0+k} = e_{n_0+p-1-k}$ for $0 \leq k < p$. In other words, the walk obtained for $\mu'$ is the same as that for $\mu$ except that the cycle $S(\mu)$ is traversed in the reverse direction. It follows that $S(\mu) = S(\mu')$, as desired.
\end{proof}

\begin{corollary}
	Let $G$ be a bipartite graph. If $\mc(G)$ contains a nonempty channel, then $m_G$ is even. 
\end{corollary}
\begin{proof}
	If $C\in\mathcal C(G)$ is nonempty, then we claim the map on matchings of $G$ given by
	$$ \mu \longmapsto \mu' := \mu\oplus S(\mu) $$
	is an involution with no fixed points. By our discussion at the start of the section, we just need to show that $S(\mu') = S(\mu)$. This follows directly from Theorem \ref{cycle}(iv).
\end{proof}

\begin{remark}\label{action}
	It is interesting (and rather inconvenient) to note that the action of channels on matchings defined above cannot in general be extended to a group action of $\mc(G)$ or $\cb(G)$ since, for instance, the action of two distinct channels need not commute. Such a group action would be a very useful combinatorial tool. We give some thoughts on this at the end of the paper.
\end{remark}

\subsection{Even $m_G$ implies existence of a channel}
Proving the converse statement will take some different machinery, which will turn out to have more general applications. 
The following lemma describes how channels are affected by removal of an edge. We restrict to bipartite $G$ for clarity---one can generalize the argument to all graphs if care is taken about how vertices can appear in channels, but for brevity we will not do so.

Suppose that $G$ has an edge $e=\edge{b}{w}$ such that the edge-deleted graph $G-e$ has a channel containing $b$. The following lemma states that $G$ and the vertex-deleted graph $G-\{b,w\}$ have the same number of channels on their black vertices. (Recall that $\cb(G)$ denotes the set of channels containing only black vertices in a bipartite graph $G$.)

\begin{lemma}[Channel Routing Lemma]\label{routing}
	Let $G=(V,E)$ be a bipartite graph, and fix an edge $e=\edge{b}{w}$ (with vertices the corresponding colors). Define the subgraphs
	$$G^e = G - e$$
	$$\text{and} $$
	$$G' = G - \{b,w\}.$$
	Then the following statements hold:
	\begin{enumerate}[label=(\roman*)]
		\item The channels of $\cb(G)$ not containing $b$ are exactly those of $\cb(G^e)$ not containing $b$.
		
		\item If there is a channel $B\in \cb(G^e)$ such that $b\in B$, then there is a bijection 
		$$\cb(G) \longleftrightarrow \cb(G')$$
		preserving channels which do not contain $b$.
	\end{enumerate}
\end{lemma}
\begin{proof}\ 
	\begin{enumerate}[label=(\roman*)]
		\item 
		Removing/adding the edge $e$ can only change $N(v) \cap C$ for a white vertex $v$ if $v=w$ and $b \in C$. Thus it cannot affect any channel on black vertices that does not contain $b$.
		
		\item 
		We construct a bijection 
		$$\cb(G) \longleftrightarrow \cb(G') $$
		by sending channels in $\cb(G)$ according to
		$$f\colon C \longmapsto \begin{cases}
			C & \text{ if } b\not\in C, \\
			B\oplus C & \text{ if } b\in C.
		\end{cases}$$
		This map is well-defined: for any white vertex $v \neq w$, $N(v)$ is the same in both $G$ and $G^e$. Hence $|N(v) \cap C|$, $|N(v) \cap B|$, and $|N(v) \cap (B \oplus C)|$ are all even by the evenness constraint of channels and properties of symmetric difference, so $|N(v) \cap f(C)|$ is even. Since $f(C)$ never contains $b$, removing vertices $b$ and $w$ does not affect $|N(v) \cap f(C)|$, so $f(C)$ is a channel in $\cb(G')$.
		
		
		We define the inverse map similarly. For a channel $C$ in $\cb(G')$, define
		$$g\colon C \longmapsto \begin{cases}
		C & \text{ if $|N(w)\cap C|$ is even}, \\
		B\oplus C & \text{ if $|N(w)\cap C|$ is odd},
		\end{cases}$$
		where $N(w)$ refers to the neighborhood of $w$ in $G$. This map is again well-defined. To see this, note that at any white vertex $v \neq w$, $|N(v) \cap g(C)|$ is even as before. At $w$, if $g(C)=C$, then $|N(w) \cap g(C)| = |N(w) \cap C|$ is even. Otherwise, $|N(w) \cap C|$ is odd, but $|N(w) \cap B|$ is also odd since $N(w)$ differs from the neighborhood of $w$ in $G^e$ only in $b$, which lies in $B$. Thus again $|N(w) \cap g(C)| = |N(w) \cap (B \oplus C)|$ is even.
		
		
		Finally, we verify that these maps are inverses. Any channel on black vertices that does not contain $b$ satisfies the first condition in both definitions, so $f$ and $g$ act on them as the identity. For channels $C$ in $\cb(G)$ that do contain $b$, the set $B\oplus C$ has odd intersection with $N(w)$ since $B$ does, and therefore $g\circ f$ is the identity. Channels in $\cb(G')$ do not contain $b$, so those that have odd intersection with $N(w)$ map to a channel $B\oplus C$ which contains $b$. Thus $f\circ g$ is the identity as well, and the claim follows. \qedhere
	\end{enumerate}
\end{proof}

Channel routing is a versatile tool. To begin with, let us use it to prove constructively the claim titling this section.

\begin{theorem}
	If a bipartite graph $G$ has an even number of matchings, then it has a nonempty channel.
\end{theorem}
\begin{proof}
    Assume without loss of generality that $G$ has at least as many black vertices as white vertices. We will show that $\cb(G)\neq 0$. First note that $G$ is nonempty since the empty graph has one matching, which is odd. We proceed by induction on $|V|+|E|$. If every white vertex of $G$ has even degree, then we can take the black vertices of $V$ to be our channel. Otherwise, some white vertex $w$ has odd degree. Then for some edge $e=\edge{b}{w}$, the subgraph $G-\{b,w\}$ has an even number of matchings, since otherwise, by Proposition \ref{matchprop},
	$$m_G = \sum_{b\colon \edge{b}{w}\in E}m_{G-\{b,w\}} $$
	would be the sum of an odd number of odd numbers, which would be odd, a contradiction. Define $G^e :=G-e$ and $G':=G-\{b,w\}$. Now, because $m_{G'}$ is even, so is 
	$$m_G- m_{G'} = m_{G^e}.$$
	Note $G'$ cannot be the empty graph since that would imply $G$ is a single edge, which would have an odd number of matchings, and likewise $G^e$ is also nonempty. Thus by the inductive hypothesis both $\cb(G')$ and $\cb(G^e)$ have nonempty channels.
	
	Let $B\in \cb(G^e)$ be a nonempty channel. If $b\not\in B$, then channel routing (Lemma~\ref{routing}(i)) implies $B$ is a channel of $G$. Otherwise, $b\in B$. Let $C\in \cb(G')$ be a nonempty channel. Then applying the map  
	$$g\colon C \longmapsto \begin{cases}
		C & \text{ if $|N(w)\cap C|$ is even} \\
		B\oplus C & \text{ if $|N(w)\cap C|$ is odd}
		\end{cases}$$
	described in channel routing (Lemma~\ref{routing}(ii)) gives us a nonempty channel in $\cb(G)$. In either case, $\cb(G)\neq 0$.
\end{proof}

Channel routing allows us to remove one vertex pair at a time from our graph while keeping track of the available channels. This is particularly useful when we have a channel in $G^e$ containing $b$ so that condition (ii) of channel routing holds. Let us see how this can be used for subgraphs of the square lattice.

\begin{example1}
Refer to Figure \ref{diagdigdiag}. Let $G$ be a subgraph of the square lattice such that each internal face is a unit square. Pick a diagonal of $G$ that starts at a corner vertex $v$ and ends at an opposing side vertex $b$, as in Figure \ref{digdiag}. 

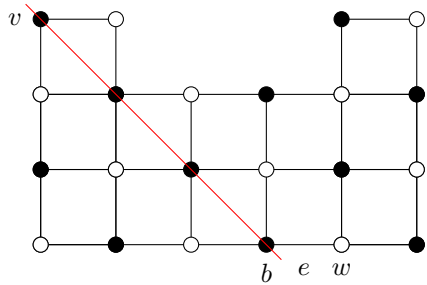
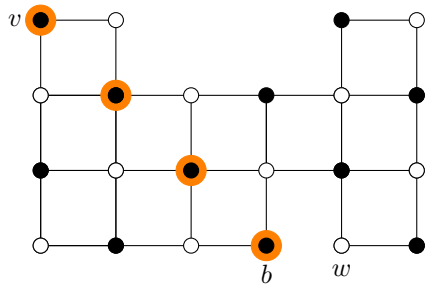
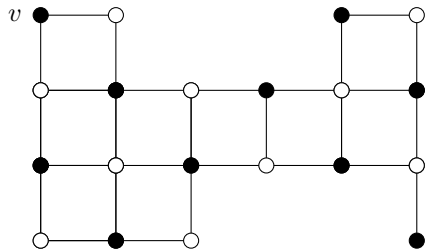
\begin{figure}
    \begin{subfigure}{\textwidth}
    \centering
    \begin{tikzpicture}
        \draw (0,0) grid (1, 3);
        \draw (0,0) grid (5, 2);
        \draw (4,0) grid (5,3);
        \node[label=left:{$v$}] at (0,3) {};
        \node[label=below:{$b$}] at (3,0) {};
        \node[label=below:{$e$}] at (3.5,0) {};
        \node[label=below:{$w$}] at (4,0) {};
        \begin{scope}[shift={(-1,0)}]
        \vertices{1}{0}{2}{3};
        \vertices{1}{0}{6}{2};
        \vertices{5}{0}{6}{3};
        \end{scope}
        \draw[red] (-.2,3.2) -- (3.2,-.2);
    \end{tikzpicture}
    \caption{The graph $G$ with a diagonal selected.}
    \label{digdiag}
    \end{subfigure}

    \vspace{.5cm}
    \begin{subfigure}{\textwidth}
    \centering
    \begin{tikzpicture}
        \draw (0,0) grid (1, 3);
        \draw (0,0) grid (3, 2);
        \draw (3,1) grid (4,2);
        \draw (4,0) grid (5,3);
        \node[label=left:{$v$}] at (0,3) {};
        \node[label=below:{$b$}] at (3,0) {};
        \node[label=below:{$w$}] at (4,0) {};
        \foreach \x in {0, 1,..., 3} {
		    \draw[orange,fill] (\x,3-\x) circle(.2);
	    }
        \begin{scope}[shift={(-1,0)}]
        \vertices{1}{0}{2}{3};
        \vertices{1}{0}{6}{2};
        \vertices{5}{0}{6}{3};
        \end{scope}
        
    \end{tikzpicture}
    \caption{The graph $G^e$. The vertices on the diagonal form a channel for this graph.}
    \label{fig:my_label}
    \end{subfigure}
    
    \vspace{.5cm}
    \begin{subfigure}{\textwidth}
    
    \centering
    \begin{tikzpicture}
        \draw (0,0) grid (1, 3);
        \draw (0,0) grid (2, 2);
        \draw (2,1) grid (4,2);
        \draw (5,0) grid (5,1);
        \draw (4,1) grid (5,3);
        \node[label=left:{$v$}] at (0,3) {};
        \begin{scope}[shift={(-1,0)}]
        \vertices{1}{0}{2}{3};
        \vertices{1}{0}{3}{2};
        \vertices{1}{1}{6}{2};
        \vertices{6}{0}{6}{3};
        \vertices{5}{1}{5}{3};
        \end{scope}
        
    \end{tikzpicture}
    \caption{The graph $G'$. Channel routing implies this graph has the same number of channels as $G$.}
    \label{fig:diaged}
    \end{subfigure}
    \caption{Given a diagonal in a subgraph of the square lattice, we may apply channel routing to the edge $e$ containing the last vertex $b$ and meeting the diagonal at an obtuse angle.}
    \label{diagdigdiag}
\end{figure}

Let $e=\edge{b}{w}$ be the unique edge containing $b$ that forms an obtuse angle with the diagonal. (If there is more than one such edge, then we are not in a situation where the method of this example applies.) Then the graph $G^e = G -e$ has a channel $B$ given by the vertices on the diagonal between $v$ and $b$. 

In particular, $b\in B$. Thus, with $G' := G-\{b,w\}$, channel routing implies that
$$|\cb(G)| = |\cb(G')|. $$
This implies, for instance, that $m_G$ and $m_{G'}$ have the same parity (by Theorem \ref{twos}). \end{example1}

In some cases, repeated application of channel routing can reduce our graph to one with known properties. Rather than deleting each vertex pair individually and examining the channels of each intermediate graph, the following theorem allows us to check an analog of the channel routing condition on a single graph to remove all of the vertex pairs at once.

\begin{theorem}\label{multchannel}
	Let $G$ be a bipartite graph with $n$ vertex disjoint edges $e_1=\edge{b_1}{w_1},...,e_n=\edge{b_n}{w_n}$ selected. 
	Set
	\begin{align*}
	G^e &= G-\{e_1,...,e_n\},\\
	G' &= G-\{b_1,w_1,...,b_n,w_n\}.
	\end{align*}
	Suppose that, for all $i$, there exists a channel $B_i\in\cb(G^e)$ such that 
	\[B_i\cap \{b_1,...,b_n\} \subseteq \{b_i\}.\]
	Then $m_G$ and $m_{G'}$ have the same parity. If equality holds for all $i$, then $|\cb(G)|=|\cb(G')|$.
\end{theorem}
\begin{proof}
	If any nonempty channel in $\cb(G^e)$ uses none of $b_1,...,b_n$, then it is also a channel in $G$ and $G'$, so both $m_G$ and $m_{G'}$ are even. Otherwise there exist channels $B_1,...,B_n$ in $\cb(G^e)$ such that
	$B_i\cap \{b_1,...,b_n\} = \{b_i\}$
	for all $i$.
	
	Since $b_1\in B_1$, channel routing (Lemma~\ref{routing}(ii)) implies $|\cb(G)| = |\cb(G-\{b_1,w_1\})|$. Since $b_1\not\in B_i$ for $i > 1$, the channels $B_2,...,B_n$ are also channels in $\cb(G-\{b_1,w_1\})$. Thus we may proceed by induction to find
	$|\cb(G)| = |\cb(G')|$,
	as desired.
\end{proof}

As one consequence, we can use Theorem~\ref{multchannel} to determine the parity of the number of domino tilings of a rectangle grid graph.

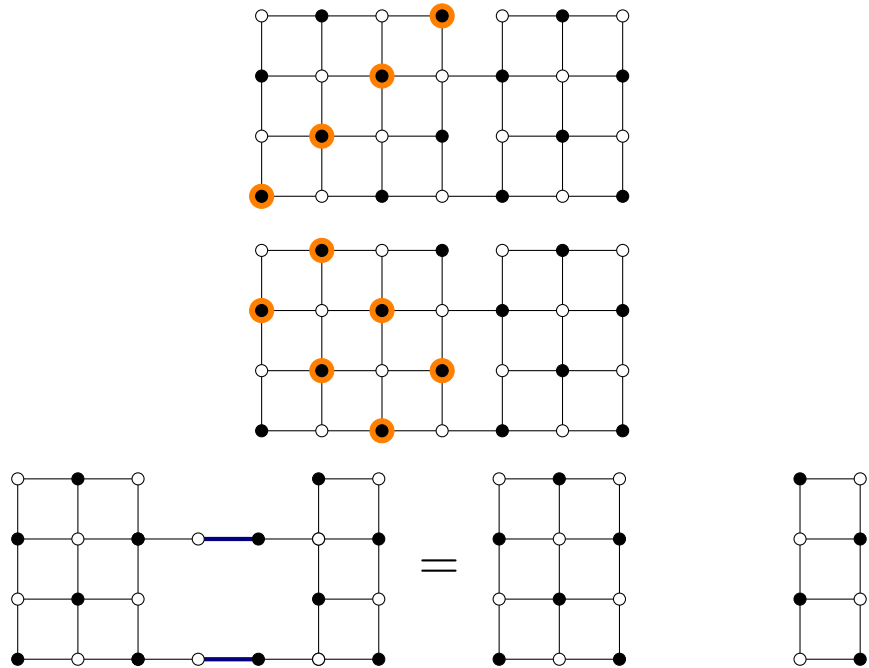
\begin{figure}[H]
	\centering
	\begin{tikzpicture}[scale=0.8,yscale=-1]
	\draw[thin] (0,0) grid (3,3);
	\draw[thin] (4,0) grid (6,3);
	\draw (3,1) -- (4,1);
	\draw (3,3) -- (4,3);
	\foreach \x in {0, 1,..., 3} {
		\draw[orange,fill] (\x,3-\x) circle(.2);
	}
	\begin{scope}[shift={(-1,0)}]
	\vertices{1}{0}{7}{3};
	\end{scope}
	\end{tikzpicture} \\ \vspace{1em}
	\begin{tikzpicture}[scale=0.8, yscale=-1]
	\draw[thin] (0,0) grid (3,3);
	\draw[thin] (4,0) grid (6,3);
	\draw (3,1) -- (4,1);
	\draw (3,3) -- (4,3);
	\foreach \x in {0, 1, 2} {
		\draw[orange,fill] (\x,1+\x) circle(.2);
		\draw[orange,fill] (\x+1,\x) circle(.2);
	}
	\begin{scope}[shift={(-1,0)}]
	\vertices{1}{0}{7}{3};
	\end{scope}
	\end{tikzpicture} \\ \vspace{1em}
	\begin{tikzpicture}[scale=0.8, yscale=-1]
	\draw[thin] (0,0) grid (2,3);
	\draw (2,3) -- (5,3);
	\draw (2,1) -- (5,1);
	\draw[thin] (5,0) grid (6,3);
	\draw[ultra thick,blue!50!black] 
	(3,1)--(4,1) (3,3)--(4,3);
	\begin{scope}[shift={(-1,0)}]
	\vertices{1}{0}{3}{3};
	\vertices{3}{1}{6}{1};
	\vertices{3}{3}{6}{3};
	\vertices{6}{0}{7}{3};
	\end{scope}
	
	\begin{scope}[shift={(8,0)}]
	\draw (-1,1.5) node[scale=2] {$=$};
	
	\draw[thin] (0,0) grid (2,3);
	\draw[thin] (5,0) grid (6,3);
	\begin{scope}[shift={(-1,0)}]
	\vertices{1}{0}{3}{3};
	\vertices{6}{0}{7}{3};
	\end{scope}
	\end{scope}
	\end{tikzpicture}
	\caption{The first two diagrams show $G^e$ and its associated channels in orange. By Theorem \ref{multchannel}, the number of matchings of $G'$ (the third figure) has the same parity as that of the original graph $G=\mathcal R_{4\times 7}$. Any matching of $G'$ will use the purple edges, so this result continues to hold if we remove those vertex pairs.}
	\label{fig:digger}
\end{figure}


\begin{proposition} \label{rectodd}
	The $m \times n$ rectangular grid graph has an odd number of matchings if and only if $gcd(m+1,n+1)=1$.
\end{proposition}
\begin{proof}
    If $m=n = 0$, then $\m= 1$. If $m=n>0$, then $\m$ is even, either by Proposition \ref{step-diag} for $n$ even or since $\m=0$ for $n$ odd. Thus we may assume without loss of generality that $m < n$.

	Refer to Figure \ref{fig:digger}. Declare the lower left vertex of $\mathcal R_{m\times n}$ to be black. Let $e_1,...,e_r$ (where $r=\lceil m/2 \rceil$) be the edges between the black vertices in the $m$th column and the white vertices in the $(m+1)$th column. Set $G^e = G-\{e_1,...,e_r\}$. 
	The channels constructed in Proposition \ref{step-diag} give $r$ independent channels in $\cb(\mathcal R_{m \times m})$. These are also valid channels of $G^e$ which satisfy the hypothesis of Theorem \ref{multchannel}. Thus by that theorem, we may remove black vertices in column $m$ and white vertices in column $m+1$ while preserving the parity of $m_G$. We are left with $G'$, a graph consisting of two rectangle grid graphs connected by bridges as in the third row of Figure \ref{fig:digger}. Since there are an equal number of white and black vertices in each rectangle, any matching of $G'$ must use the middle edge on each bridge. Thus we may remove the rest of the vertices in those columns without changing the number of matchings. The resulting graph will be the disjoint union of a $m\times (m-1)$ rectangle and a $m \times (n-m-1)$ rectangle. The result then follows by induction on the size of the rectangle, since
	\[\gcd(m+1, n-m-1+1) = \gcd(m+1,n+1) \text{ and } \gcd(m+1, m-1+1) = 1. \qedhere\]
\end{proof}
\section{Graph operations (or, how to dig channels)}
Lemma \ref{routing} shows that in some cases deleting a vertex pair preserves the number of channels in a graph. In this section, we describe a set of local graph moves which unconditionally preserve channels. 
In many cases, these operations will allow us to compute the number of channels by reducing to a graph with no edges. 

The nature of these graph operations requires that we allow multiple edges between a pair of vertices. We continue to use $N(v)$ to denote the neighborhood of $v$; however, it may now be a multiset in which a vertex appears with multiplicity equal to the number of edges connecting it to $v$. In this case, for a channel $C$, $N(v)\cap C$ denotes the sub-multiset of $N(v)$ consisting of all vertices which appear in $C$ (with the same multiplicity as they appear in $N(v)$). All other set operations that appear involve standard sets.

We begin by introducing our operations of interest.

\subsection{Channel-preserving moves}

A \textit{2-valent vertex contraction} may be applied to any vertex $v$ of degree two that is adjacent to distinct vertices $v_1,v_2$. The resulting graph is formed by contracting the edges incident to $v$ and deleting self-loops if they occur.
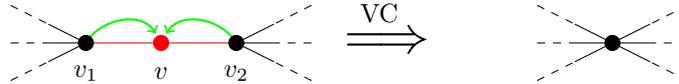
\begin{figure}[H]
	\centering
	\begin{tikzpicture}
	\draw[white] (-1.5,0) (9,0);
	
	\pgfmathsetmacro{\radius}{0.1}
	\draw[red,fill] (1,0) circle (.1);
	\draw[red] (0,0) -- (1,0) -- (2,0);
	\node (a) at (1,0) {};
	\draw[green,thick,->] (0,0) to[out=60,in=120] (a);
	\draw[green,thick,->] (2,0) to[out=120,in=60] (a);
	\draw[fill = black] (0,0) circle (\radius) (2,0) circle (\radius);
	\node at (1,-.4) {$v$};
	\node at (0,-.4) {$v_1$};
	\node at (2,-.4) {$v_2$};
	\draw[dashed] (-1,-.5) -- (-.5,-.25);
	\draw (-.5,-.25) -- (0,0);
	\draw[dashed] (-1,0) -- (-.5,0);
	\draw (-.5,0) -- (0,0);
	\draw[dashed] (-1,.5) -- (-.5,.25);
	\draw (-.5,.25) -- (0,0);
	\draw[dashed] (3,-.5) -- (2.5,-.25);
	\draw (2.5,-.25) -- (2,0);
	\draw[dashed] (3,0) -- (2.5,0);
	\draw (2.5,0) -- (2,0);
	\draw[dashed] (3,.5) -- (2.5,.25);
	\draw (2.5,.25) -- (2,0);
	
	\draw (4,0) node[scale=2]{$\implies$};
	\draw (3.9,.4) node{VC};
	
	\begin{scope}[shift={(7,0)}]
	
	\pgfmathsetmacro{\radius}{0.1}
	\draw[fill = black] (0,0) circle (\radius);
	\draw[dashed] (-1,-.5) -- (-.5,-.25);
	\draw (-.5,-.25) -- (0,0);
	\draw[dashed] (-1,0) -- (-.5,0);
	\draw (-.5,0) -- (0,0);
	\draw[dashed] (-1,.5) -- (-.5,.25);
	\draw (-.5,.25) -- (0,0);
	\draw[dashed] (1,-.5) -- (0.5,-.25);
	\draw (0.5,-.25) -- (0,0);
	\draw[dashed] (1,0) -- (0.5,0);
	\draw (0.5,0) -- (0,0);
	\draw[dashed] (1,.5) -- (0.5,.25);
	\draw (0.5,.25) -- (0,0);
	
	\end{scope}
	
	\end{tikzpicture}
\end{figure}

A \textit{doubled edge deletion} may be applied to any pair of edges $e_1,e_2$ that share the same endpoints. This operation removes $e_1$ and $e_2$ from the graph.
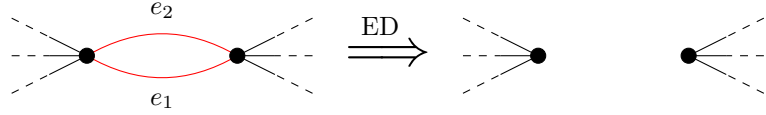
\begin{figure}[H]
	\centering
	\begin{tikzpicture}
	\draw[white] (-1.5,0) (9,0);
	\draw[red] (0,0)  to[out=30, in=150] (2,0);
	\draw[red] (0,0) to[out=-30, in=-150] (2,0);
	\node at (1,-.6) {$e_1$};
	\node at (1,.6) {$e_2$};
	\draw[fill=black] (0,0) circle (.1) (2,0) circle (.1);
	\draw[dashed] (-1,-.5) -- (-.5,-.25);
	\draw (-.5,-.25) -- (0,0);
	\draw[dashed] (-1,0) -- (-.5,0);
	\draw (-.5,0) -- (0,0);
	\draw[dashed] (-1,.5) -- (-.5,.25);
	\draw (-.5,.25) -- (0,0);
	\draw[dashed] (3,-.5) -- (2.5,-.25);
	\draw (2.5,-.25) -- (2,0);
	\draw[dashed] (3,0) -- (2.5,0);
	\draw (2.5,0) -- (2,0);
	\draw[dashed] (3,.5) -- (2.5,.25);
	\draw (2.5,.25) -- (2,0);
	
	\draw (4,0) node[scale=2]{$\implies$};
	\draw (3.9,.4) node{ED};
	
	\begin{scope}[shift={(6,0)}]
	
	\pgfmathsetmacro{\radius}{0.1}
	\draw[fill=black] (0,0) circle (.1) (2,0) circle (.1);
	\draw[dashed] (-1,-.5) -- (-.5,-.25);
	\draw (-.5,-.25) -- (0,0);
	\draw[dashed] (-1,0) -- (-.5,0);
	\draw (-.5,0) -- (0,0);
	\draw[dashed] (-1,.5) -- (-.5,.25);
	\draw (-.5,.25) -- (0,0);
	\draw[dashed] (3,-.5) -- (2.5,-.25);
	\draw (2.5,-.25) -- (2,0);
	\draw[dashed] (3,0) -- (2.5,0);
	\draw (2.5,0) -- (2,0);
	\draw[dashed] (3,.5) -- (2.5,.25);
	\draw (2.5,.25) -- (2,0);
	
	\end{scope}
	\end{tikzpicture}
\end{figure}

A \textit{forced vertex pair removal} may be applied to distinct adjacent vertices $v_1,v_2$ such that $v_1$ has degree one. The resulting graph is formed by removing $v_1,v_2$ and all edges incident to these vertices.
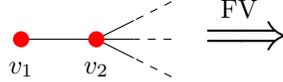
\begin{figure}[H]
	\centering
	\begin{tikzpicture}
	\draw[white] (-1.5,0) (9,0);
	\pgfmathsetmacro{\radius}{0.1}
	\draw (1,0) -- (2,0);
	\node at (1,-.4) {$v_1$};
	\node at (2,-.4) {$v_2$};
	\draw[dashed] (3,-.5) -- (2.5,-.25);
	\draw (2.5,-.25) -- (2,0);
	\draw[dashed] (3,0) -- (2.5,0);
	\draw (2.5,0) -- (2,0);
	\draw[dashed] (3,.5) -- (2.5,.25);
	\draw (2.5,.25) -- (2,0);
	\draw[fill, red] (1,0) circle (\radius) (2,0) circle (\radius);
	
	\draw (4,0) node[scale=2]{$\implies$};
	\draw (3.9,.4) node{FV};

	\end{tikzpicture}
\end{figure}

\begin{definition}
	A \textit{channel-preserving move} is one of:
	\begin{itemize}
		\item[(VC)] 2-valent vertex contraction,
		\item[(ED)] doubled edge deletion, or
		\item[(FV)] forced vertex pair removal.
	\end{itemize}
\end{definition}

As the name suggests, applying a channel-preserving move to a graph preserves the number of channels in that graph. In the following, we write $n\equiv_2 m$ to mean $n-m$ is even.

\begin{lemma}
    Let $G$ be a graph and let $G'$ be the result of applying a channel-preserving move to $G$. Then
    $$|\mc(G)| = |\mc(G')|. $$
    If additionally $G$ is bipartite, then
    \[|\cb(G)| = |\cb(G')| \text{ and } |\cw(G)| = |\cw(G')|. \]
\end{lemma}
\begin{proof}\ 
    First, ED moves clearly preserve the parity of $|N(v)\cap C|$ for all vertices $v$ and vertex sets $C$, implying the claim. We will show the result for VC moves; the argument for FV moves is similar.
    
    Let $v$ be a vertex of degree two with adjacent vertices $v_1$ and $v_2$. Call the resulting contracted vertex $w$ in $G'$. If $C\in \mc(G)$, then in order for the evenness condition to hold at $v$, it follows that $v_1\in C$ if and only if $v_2\in C$. Thus we may define
    \[C' = C \backslash \{v_1, v, v_2\} \cup W, \quad \text{where} \quad W = \begin{cases} \{w\}&\text{if $v_1,v_2 \in C$,}\\ \varnothing&\text{otherwise}. \end{cases}\]
    This preserves the neighborhood size of all unchanged vertices by replacing any occurrence of $v_1$ or $v_2$ in a neighborhood with $w$. Thus evenness holds everywhere except possibly at $w$. To see that $|N(w)\cap C'|$ is even, 
    note that 
    \[|N(w)\cap C'| = 
    |N(v_1)\cap C|+|N(v_2)\cap C|-2|\{v\}\cap C| - k|\{v_1,v_2\}\cap C|  
    ,\]
    where $k$ is the number of edges between $v_1$ and $v_2$. Since
    $| \{v_1,v_2\} \cap C | = |N(v)\cap C|$
    is even, it follows that 
    \[|N(w)\cap C'|\equiv_2 |N(v_1)\cap C| + |N(v_2)\cap C|\]
    is even and hence $C'\in \mc(G')$. 
    
    Conversely, if $C'\in \mc(G')$, then $|N(w)\cap C'|$ is even. Since
    \[ |N(w) \cap C'| = |N(v_1)\cap C'| + |N(v_2)\cap C'|, \]
    it follows that $|N(v_1)\cap C'|\equiv_2 |N(v_2)\cap C'|$. Again let $k$ denote the number of edges between $v_1$ and $v_2$. We may now define $C = C' \backslash \{w\} \cup V_1\cup V_2,$  where 
    \[V_1 = \begin{cases} \{v_1,v_2\}&\text{if $w\in C'$,}\\ \varnothing&\text{otherwise}, \end{cases}
    \quad \text { and }
    V_2 = \begin{cases} \{v\}&\text{if $|N(v_1)\cap C'|+k|\{w\}\cap C'|$ is odd,}\\ \varnothing&\text{otherwise,} \end{cases}
    \]
    which is the inverse to the map $C \mapsto C'$ described above. Again this preserves the neighborhood size of all unchanged vertices. The definition of $C$ ensures that $|N(v_1)\cap C|\equiv_2 |N(v_2)\cap C|$ is even by adding $v$ to $C$ if necessary. Also $|N(v)\cap C| = 2$ if $w\in C'$, and $0$ otherwise. Thus evenness holds at all vertices, and $C\in \mc(G)$.
    
    Further, notice that if $G$ is bipartite, then $w$ has the same color as $v_1$ and $v_2$. So if $C$ uses only vertices of a single color, then $C'$ only uses vertices of the same color. The converse also holds, since for singly colored channels, at most one of $|N(v_1)\cap C'| \neq 0$ and $w\in C'$ can hold. Hence we have a bijection $\mc(G)\leftrightarrow \mc(G')$ that descends to $\cb$ and $\cw$ for bipartite graphs.
\end{proof}

We will call a graph \textit{reducible} if it can be reduced to a set of degree 0 vertices using only channel-preserving moves. 
\begin{theorem}
	Let $G$ be a reducible graph. Then the number of degree 0 vertices remaining after $G$ has been fully reduced is the dimension of $\mc(G)$. In particular, this number is independent of the choice of channel-preserving moves used to reduce the graph. 
\end{theorem}
\begin{proof}
	Since channel-preserving moves preserve channels, we just need to show that a set of $n$ vertices of degree 0 has $2^n$ channels. This is clear, since any subset of these vertices is a valid channel.
\end{proof}

\begin{example1}
	The graph $G$ shown in Figure \ref{fig:channelexample} is a reducible graph. Figure \ref{reduce} shows a possible sequence of channel-preserving moves. Because $G$ reduces to two vertices of degree 0, $G$ must have $2^2$ channels. This is indeed the case; the three nonempty channels are shown in Figure \ref{fig:channelexample}.
\end{example1}
	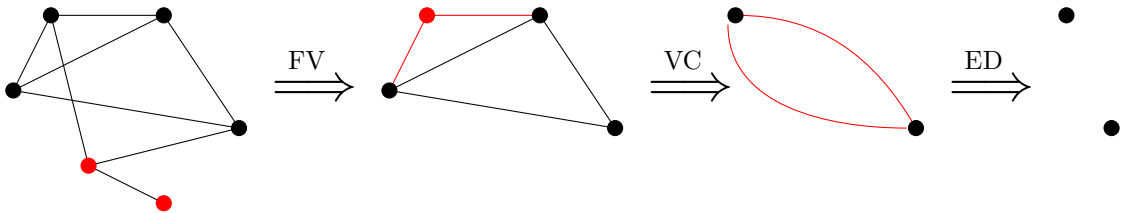
\begin{figure}[H]
		\centering
		\begin{tikzpicture}
		\pgfmathsetmacro{\radius}{0.1}
		\node (a) at (0,0){};
		\node (b) at (0.5,1){};
		\node (c) at (2,1){};
		\node (d) at (3,-0.5){};
		\node (e) at (2,-1.5){};
		\node (f) at (1,-1){};
		\draw (a.center) -- (c.center) -- (d.center) -- (a.center);
		\draw (a.center) -- (b.center);
		\draw (b.center) -- (c.center);
		\draw (b.center) -- (f.center);
		\draw (e.center) -- (f.center);
		\draw (f.center) -- (d.center);
		\draw[fill=black] (a) circle (\radius);
		\draw[fill=black] (b) circle (\radius);
		\draw[fill=black] (c) circle (\radius);
		\draw[fill=black] (d) circle (\radius);
		\draw[red,fill=red] (e) circle (\radius);
		\draw[red,fill=red] (f) circle (\radius);
		
		\draw (4,0) node[scale=2]{$\implies$};
		\draw (3.9,.4) node{FV};
		
		\begin{scope}[shift={(5,0)}]
		
		\node (a) at (0,0){};
		\node (b) at (0.5,1){};
		\node (c) at (2,1){};
		\node (d) at (3,-0.5){};
		\draw (a.center) -- (c.center) -- (d.center) -- (a.center);
		\draw[red] (a.center) -- (b.center);
		\draw[red] (b.center) -- (c.center);
		\draw[fill=black] (a) circle (\radius);
		\draw[red,fill=red] (b) circle (\radius);
		\draw[fill=black] (c) circle (\radius);
		\draw[fill=black] (d) circle (\radius);
		\end{scope}
		
		\draw (9,0) node[scale=2]{$\implies$};
		\draw (8.9,.4) node{VC};
		
		\begin{scope}[shift={(9,0)}]
		
		\node (b) at (0.5,1){};
		\node (d) at (3,-0.5){};
		\draw[red] (d.center) to[out=120, in=0] (b) to[out=-90, in=180] (d);
		\draw[fill=black] (0.6,1) circle (\radius);
		\draw[fill=black] (d) circle (\radius);
		\end{scope}
		
		\draw (13,0) node[scale=2]{$\implies$};
		\draw (12.9,.4) node{ED};
		\draw[fill=black] (14,1) circle (\radius);
		\draw[fill=black] (14.6,-.5) circle (\radius);
		\end{tikzpicture} 
		\caption{A reduction of the graph shown in Figure \ref{fig:channelexample}.}\label{reduce}
	\end{figure}

\begin{example1}
	Not all graphs are reducible. Figure \ref{nonreduce} shows a planar bipartite graph that admits no channel-preserving moves. 
\end{example1}
	
	\begin{figure}[H]
		\centering
		\begin{tikzpicture}
		\pgfmathsetmacro{\radius}{0.1}
		\draw[fill=black] (-2,-2) circle (\radius)--
						  (2,-2) circle (\radius)--
						  (2,2) circle (\radius)--
						  (-2,2) circle (\radius)--(-2,-2);
		\draw[fill=black] (-1,-1) circle (\radius)--
		(1,-1) circle (\radius)--
		(1,1) circle (\radius)--
		(-1,1) circle (\radius)--(-1,-1);
		\draw[fill=black] (1,-1)--(2,-2);
		\draw[fill=black] (1,1)--(2,2);
		\draw[fill=black] (-1,1)--(-2,2);
		\draw[fill=black] (-1,-1)--(-2,-2);
		\draw[fill=white] (-1,1) circle (.1) (1,-1) circle(.1)
		 (2,2) circle (.1) (-2,-2) circle (.1);
		\end{tikzpicture}
		\caption{}
		\label{nonreduce}
	\end{figure}

As this example indicates, it is not clear at first if reducibility occurs often enough to be useful---we would like to have a simple structural property that will imply reducibility. The following lemma will be useful for identifying potential VC or FV moves. Note that the degree of a face in a planar graph is its degree as a vertex in the dual graph.

\begin{lemma}[Corner identification]\label{cornerlemma}
	Let $G$ be a connected planar graph with at least two vertices such that each internal vertex and internal face have degree at least 4. Let $b$ be the number of external vertices in $G$. Then the average degree of the external vertices in $G$ is at most 
	$$3 - \frac{4}{b}. $$
	In particular, there is an external vertex of degree less than 3. 
\end{lemma}
\begin{proof}
	By assumption, the degree of each vertex is at least one. We will use the following notation:
	\begin{align*}
	b &= \text{ the number of external vertices}, \\
	i &= \text{ the number of internal vertices}, \\
	D &= \text{ the total degree of all external vertices}, \\
	e &= \text{ the number of edges}, \\
	f &= \text{ the number of internal faces.}
	\end{align*}
	By planarity of $G$, we have that
	\begin{equation}\label{euler}b + i -e + f =1 .\end{equation}
	Because the total degree of all internal and external vertices is $2e$, and the degree of each internal vertex is at least 4, it follows that
	$$2e \geq 4i + D .$$ 
	The total degree of all internal and external faces is also $2e$. 
	Since the degree of the external face is the total number of vertices counted with multiplicity as one proceeds through its boundary cycle, this number is at least $b$. Thus the total degree of all internal faces is at most $2e-b$, and consequently
	$$2e-b \geq 4f .$$ 
	Adding the two inequalities and utilizing Euler's formula (\ref{euler}) then gives
	$$D \leq 3b - 4 .$$
	Dividing by $b$ reveals the claim.
\end{proof}

A planar graph is called \emph{inner Eulerian} if all internal vertices have even degree. The following results show the utility of this class of graphs.

\begin{lemma}
	Channel-preserving moves preserve the property of being inner Eulerian.
\end{lemma}
\begin{proof}
	Assume $G$ is inner Eulerian. Let $v$ be a vertex of $G$ of degree 2. Let the degrees of the two vertices $v_1$ and $v_2$ adjacent to $v$ be $d_1$ and $d_2$, respectively, and let the number of edges connecting $v_1$ and $v_2$ be $n$. If we perform a VC move on $v$, then we will be left with a vertex of degree $d_1+d_2-2n-2$. If $d_1 + d_2$ is even, then we are done. Otherwise, one of the two vertices adjacent to $v$ has odd degree and must therefore be incident to the external face since $G$ is inner Eulerian. Thus the vertex resulting from contracting $v$ is also external and therefore has unconstrained degree. Thus VC moves preserve being inner Eulerian.
	
	Let $v$ now be a vertex of degree one. Then $v$ is external, so the vertex $v'$ adjacent to $v$ is also external. If we apply an FV move to $v$, then $v$ and $v'$ will be removed. This will change only the degree of vertices adjacent to $v'$. However, once $v'$ is removed all of these vertices will be adjacent to the external face and therefore have unconstrained degree. Thus FV moves preserve being inner Eulerian. 
	
	Finally, ED moves preserve the parity of the degree of every vertex and thus also preserve the property of being inner Eulerian.
\end{proof}

\begin{theorem}
	Let $G$ be an inner Eulerian bipartite graph. Then $G$ is reducible.
\end{theorem}
\begin{proof} 
	By the previous lemma, we are free to perform any channel-preserving move while remaining inner Eulerian. If we can show that it is always possible to perform a channel-preserving move on an inner Eulerian bipartite graph $G$ with at least one edge, then the result will follow by induction on the number of edges.

	Without loss of generality, we may assume that $G$ is connected, so every vertex has degree at least one. If any internal vertex has degree $2$, then we can perform a VC move. Otherwise every internal vertex has degree at least $4$ (it must be even as $G$ is inner Eulerian). If any internal face has degree $2$, then we can perform an ED move. Otherwise every internal face has degree at least $4$ (it must be even as $G$ is bipartite).
	
	Thus if neither a VC or ED move is possible in the graph interior, then by Lemma~\ref{cornerlemma}, there exists an external vertex of degree $1$ or $2$ to which an FV or VC move can be applied.
\end{proof}

\subsection{Contracting diagonals}
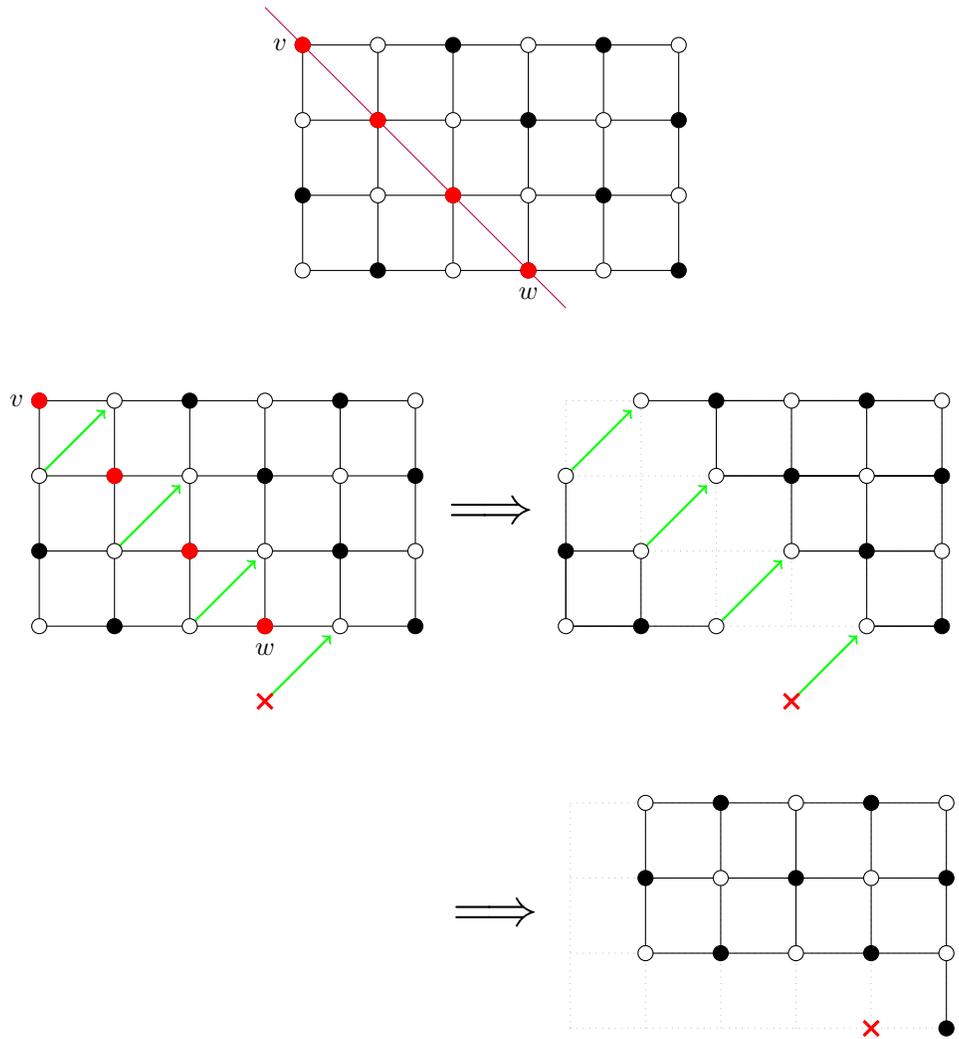
\begin{figure}
	\centering
	\begin{tikzpicture}
	\draw (0,0) grid (5,3);
	\begin{scope}[shift={(-1,0)}]
	\vertices{1}{0}{6}{3};
	\end{scope}
	
	\draw[purple] (-.5,3.5) -- (3.5,-.5);
	
	\draw[red,fill] (0,3) circle (.1);
	\draw[red,fill] (1,2) circle (.1);
	\draw[red,fill] (2,1) circle (.1);
	\draw[red,fill] (3,0) circle (.1);
	
	\node at (-.3,3) {$v$};
	\node at (3,-.3) {$w$};
	\end{tikzpicture}
	\vspace{1cm}
	
	\begin{tikzpicture}
	\node (a) at (1,3){};
	\draw[thick, green,->] (0,2) -- (a);
	\node (b) at (2,2){};
	\draw[thick, green,->] (1,1) -- (b);
	\node (c) at (3,1){};
	\draw[thick, green,->] (2,0) -- (c);
	\node (d) at (4,0){};
	\draw[thick, green,->] (3,-1) -- (d);
	\draw (0,0) grid (5,3);
	\begin{scope}[shift={(-1,0)}]
	\vertices{1}{0}{6}{3};
	\end{scope}
	\draw[red,fill] (0,3) circle (.1);
	\draw[red,fill] (1,2) circle (.1);
	\draw[red,fill] (2,1) circle (.1);
	\draw[red,fill] (3,0) circle (.1);
	\node at (-.3,3) {$v$};
	\node at (3,-.3) {$w$};
	
	\draw[red, very thick] (3-.1, -1-.1) -- (3+.1,-1+.1);
	\draw[red, very thick] (3+.1, -1-.1) -- (3-.1,-1+.1);
	
	\draw (6,1.5) node[scale=2]{$\implies$};
	
	\begin{scope}[shift={(7,0)}]
	\draw[opacity=.5,gray, dotted] (0,0) grid (5,3);
	\node (a) at (1,3){};
	\draw[thick, green,->] (0,2) -- (a);
	\node (b) at (2,2){};
	\draw[thick, green,->] (1,1) -- (b);
	\node (c) at (3,1){};
	\draw[thick, green,->] (2,0) -- (c);
	\node (d) at (4,0){};
	\draw[thick, green,->] (3,-1) -- (d);
	\draw (0,0) grid (1,1);
	\draw (1,3) grid (5,3);
	\draw (2,2) grid (5,2);
	\draw (3,1) grid (5,1);
	\draw (4,0) grid (5,0);
	\draw (0,0) grid (0,2);
	\draw (0,0) grid (2,0);
	
	\draw (2,3) grid (5,2);
	\draw (3,2) grid (5,1);
	\draw (4,1) grid (5,0);
	\begin{scope}[shift={(-1,0)}]
	\vertices{1}{0}{2}{1};
	\vertices{1}{2}{1}{2};
	\vertices{3}{0}{3}{0};
	\vertices{2}{3}{6}{3};
	\vertices{3}{2}{6}{2};
	\vertices{4}{1}{6}{1};
	\vertices{5}{0}{6}{0};
	\end{scope}
	
	\draw[red, very thick] (3-.1, -1-.1) -- (3+.1,-1+.1);
	\draw[red, very thick] (3+.1, -1-.1) -- (3-.1,-1+.1);
	\end{scope}
	\end{tikzpicture}
	\vspace{1cm}
	
		\centering
	\begin{tikzpicture}
	\draw[opacity=0] (0,0) grid (5,3);
	\node[opacity=0] at (-.3,3) {$v$};
	\draw (6,1.5) node[scale=2]{$\implies$};
	\draw[opacity=0] (7,1.5) node[scale=2]{$\implies$};
	\begin{scope}[shift={(7,0)}]
	\draw[opacity=.5,gray,dotted] (0,0) grid (5,3);
	\draw (1,1) grid (5,3);
	\draw (5,0) grid (5,1);
	\begin{scope}[shift={(-1,0)}]
	\vertices{2}{1}{6}{3};
	\vertices{6}{0}{6}{0};
	\end{scope}
	\draw[red, very thick] (4-.1, 0-.1) -- (4+.1,0+.1);
	\draw[red, very thick] (4+.1, 0-.1) -- (4-.1,0+.1);
	\end{scope}
	\end{tikzpicture}
	\caption{Contracting the highlighted diagonal by deleting the diagonal vertices and merging the vertices immediately opposite.}\label{diagid}
\end{figure}

When our graph is a suitably nice subgraph of the square lattice, there is often a useful sequence of channel-preserving moves available called a \textit{diagonal contraction}. Pick a degree 2 vertex $v$ that is a corner of the graph. Then $v$ defines a unique diagonal passing through it, as in the top of Figure \ref{diagid}. 

We say that the diagonal is \textit{contractible} if each internal vertex and each internal face it intersects have degree 4. Such a diagonal can be contracted as follows. Consider all vertices on the diagonal from $v$ to $w$, the last vertex on the diagonal before it reaches the external face. For each such vertex $v_i$, delete $v_i$ and combine each neighbor of $v_i$ with its mirror image across the diagonal, as shown in Figure \ref{diagid}. If a vertex is to combine with a missing vertex (denoted by a red $\times$ in the figure), then that vertex is instead deleted.

\begin{figure}
	\centering
	\begin{tikzpicture}
	\draw[opacity=.5,gray,dotted] (0,0) grid (2,2);
	\node (a) at (0,1){};
	\draw[purple,thick,dashed,opacity=.5] (-.2,-.2) -- (1.2,1.2);
	\draw[thick,green,->] (0,1) -- (a);
	\draw (0,0) grid (1,2);
	\draw (1,0) grid (2,1);
	\draw[red] (0,1) -- (0,0) -- (1,0);
	
	\node (a) at (0,1){};
	\draw[thick,green,->] (1,0) -- (a);
	
	\vertices{0}{0}{1}{2};
	\vertices{2}{0}{2}{1};
	\draw[fill,red] (0,0) circle (.1);
	\node at (0,-.4) {$v$};
	\node at (1.3,1.3) {$w$};

	\begin{scope}[shift={(4,0)}]
	\draw (-1, 1) node[scale=2]{$\implies$};
	\draw (-1.1,1.4) node{VC};
	
	\draw[opacity=.5,gray,dotted] (0,0) grid (2,2);
	\draw (0,1) grid (1,2);
	\draw (2,0) -- (2,1) -- (1,1);
	
	\draw[red] (0,1) -- (1,1);
	\draw[red] (0,1) to[out=-60,in=-120] (1,1);
	\draw (0,1) to[out=-90, in=180] (2,0);
	\node at (1.3,1.3) {$w$};
	
	\vertices{0}{1}{1}{2};
	\vertices{2}{0}{2}{0};
	\vertices{2}{1}{2}{1};
	
	\end{scope}
	
	\begin{scope}[shift={(8,0)}]
	\node (a) at (1,2){};
	\draw[thick,green,->] (2,1) -- (a);
	\draw (-1, 1) node[scale=2]{$\implies$};
	\draw (-1.1,1.4) node{ED};
	
	\draw[opacity=.5,gray,dotted] (0,0) grid (2,2);
	\draw (2,0) -- (2,1);
	\draw (0,1) -- (0,2) -- (1,2);
	\draw[red] (2,1) -- (1,1) -- (1,2);
	
	\draw (0,1) to[out=-90, in=180] (2,0);
	
	\vertices{0}{1}{1}{2};
	\vertices{2}{0}{2}{0};
	\vertices{2}{1}{2}{1};

	\node at (1.3,1.3) {$w$};
	\draw[red,fill] (1,1) circle (.1);
	\end{scope}
	
	\begin{scope}[shift={(12,0)}]
	\draw[opacity=.5,gray,dotted] (0,0) grid (2,2);
	\draw (-1, 1) node[scale=2]{$\implies$};
	\draw (-1.1,1.4) node{VC};
	
	\draw (0,1) grid (1,2);
	\vertices{0}{1}{1}{2};
	
	\end{scope}
	\end{tikzpicture}
	\caption{Diagonal contraction (from $v$ to $w$) that ends on a degree 4 vertex. Channels are preserved.}\label{op3}
\end{figure}
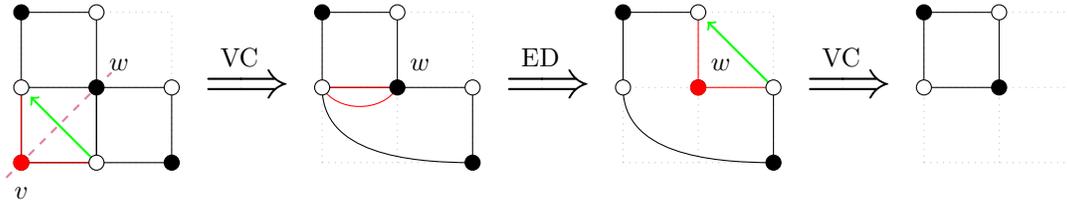

\begin{theorem}\label{diagthm}
	Let $G$ be a subgraph of the square lattice. Let $G'$ be the result of applying a diagonal contraction to a contractible diagonal from a black corner vertex $v$ to a vertex $w$. If $w$ has degree 2, then 
	$$|\cb(G)| = 2|\cb(G')|  \text{ and } |\cw(G)|=|\cw(G')|.$$
	Otherwise,
	$$|\cb(G)| = |\cb(G')|  \text{ and } |\cw(G)|=|\cw(G')|.$$
\end{theorem}
\begin{proof}
	We will show that a diagonal contraction move consists of a sequence of channel-preserving moves. 
	If $v=w$, then diagonal contraction is just a VC move. Otherwise, the diagonal passes through an internal face which $v$ shares with exactly one other black vertex $v'$ and two white vertices $v_1$ and $v_2$. Applying a VC move to $v$ combines $v$, $v_1$, and $v_2$. Since $v_1$ and $v_2$ were both adjacent to $v'$, there are now two edges between $v$ and $v'$. Thus we can apply an ED move to this edge pair. Now, if $v'$ has degree $2$ at this point, then we may repeat this argument starting with $v'$.
	
	Otherwise, $v'$ now has degree $1$ or $0$, so it originally had degree 3 or 2, so it must have been $w$, the last vertex on the diagonal. We can either apply an FV move to the degree 1 vertex or remove the degree 0 vertex from the graph. (The possible ending scenarios are shown in Figures \ref{op1}, \ref{op2}, and \ref{op3}.) Our resulting graph is the diagonal contraction $G'$ essentially by definition. All of our moves were channel preserving, except for removing the degree 0 vertex $w$; this occurs if and only if $w$ had degree 2 in $G$. Removing a degree 0 black vertex from the graph halves the number of channels on black vertices and preserves the number of channels on white vertices. Thus the claim is shown.
\end{proof}

\begin{figure}
	\centering
	\begin{tikzpicture}
	\node at (-.4,2){$v$};
	\node at (2,-.4){$w$};
	\node (a) at (0,1){};
	\draw[purple,thick,dashed,opacity=.5] (-.3,2.3) -- (2.3,-.3);
	\draw[thick,green,->] (1,2) -- (a);
	\begin{scope}[xscale=-1, shift={(-2,0)}]
	\draw[opacity=.5,gray,dotted] (0,0) grid (2,2);
	\draw (0,0) grid (2,1);
	\draw (1,0) grid (2,2);
	\draw[red] (1,2) -- (2,2) -- (2,1);
	\vertices{0}{0}{2}{1};
	\vertices{1}{1}{2}{2};
	\draw[fill,red] (2,2) circle (.1);
	\end{scope}
	
	\draw (3, 1) node[scale=2]{$\implies$};
	\draw (2.9,1.4) node{VC};

	\begin{scope}[shift={(4,0)}]
	\node at (2,-.4){$w$};
	\begin{scope}[xscale=-1, shift={(-2,0)}]
	\draw[opacity=.5,gray,dotted] (0,0) grid (2,2);
	\draw (0,0) grid (2,1);
	\draw[red] (1,1) to[out=60, in=120] (2,1);
	\draw[red] (1,1) -- (2,1);
	\vertices{0}{0}{2}{1};
	\end{scope}
	\end{scope}
	
	\draw (7, 1) node[scale=2]{$\implies$};
	\draw (6.9,1.4) node{ED};
	
	\begin{scope}[shift={(8,0)}]
	\node (a) at (1,0){};
	\draw[thick,green,->] (2,1) -- (a);
	\node at (2,-.4){$w$};
	\begin{scope}[xscale=-1, shift={(-2,0)}]
	\draw[opacity=.5,gray,dotted] (0,0) grid (2,2);
	\draw (0,0) grid (1,1);
	\draw (1,0) -- (2,0) -- (2,1);
	\draw[red] (0,1) -- (1,1) -- (1,0);
	\vertices{0}{0}{2}{1};
	\draw[fill,red] (1,1) circle (.1);
	\end{scope}
	\end{scope}
	
	\draw (-1, -2) node[scale=2]{$\implies$};
	\draw (-1.1,-1.6) node{VC};
	
	\begin{scope}[shift={(0,-3)}]
	\node at (2,-.4){$w$};
	\begin{scope}[xscale=-1, shift={(-2,0)}]
	\draw[opacity=.5,gray,dotted] (0,0) grid (2,2);
	\draw[red] (0,0)--(1,0);
	\draw[red] (0,0) to[out=60,in=120] (1,0);
	\draw (1,0)-- (2,0) -- (2,1);
	\vertices{0}{0}{2}{0};
	\vertices{2}{1}{2}{1};
	\end{scope}
	\end{scope}
	
	\draw (3, -2) node[scale=2]{$\implies$};
	\draw (2.9,-1.6) node{ED};
	
	\begin{scope}[shift={(4,-3)}]
	\node at (2,-.4){$w$};
	\begin{scope}[xscale=-1, shift={(-2,0)}]
	\draw[opacity=.5,gray,dotted] (0,0) grid (2,2);
	\draw (1,0)-- (2,0) -- (2,1);
	\vertices{0}{0}{2}{0};
	\vertices{2}{1}{2}{1};
	\draw[fill,red] (0,0) circle (.1);
	\end{scope}
	\end{scope}
	
	\draw (7, -2) node[scale=2]{$\implies$};
	
	\begin{scope}[shift={(8,-3)}]
	\begin{scope}[xscale=-1, shift={(-2,0)}]
	\draw[opacity=.5,gray,dotted] (0,0) grid (2,2);
	\draw (1,0)-- (2,0) -- (2,1);
	\vertices{1}{0}{2}{0};
	\vertices{2}{1}{2}{1};
	\end{scope}
	\end{scope}
	\end{tikzpicture}
	\caption{Diagonal contraction (from $v$ to $w$) that ends on a degree 2 vertex. Removing the degree 0 vertex $w$ removes a basis channel for $G$.}\label{op1}
\end{figure}
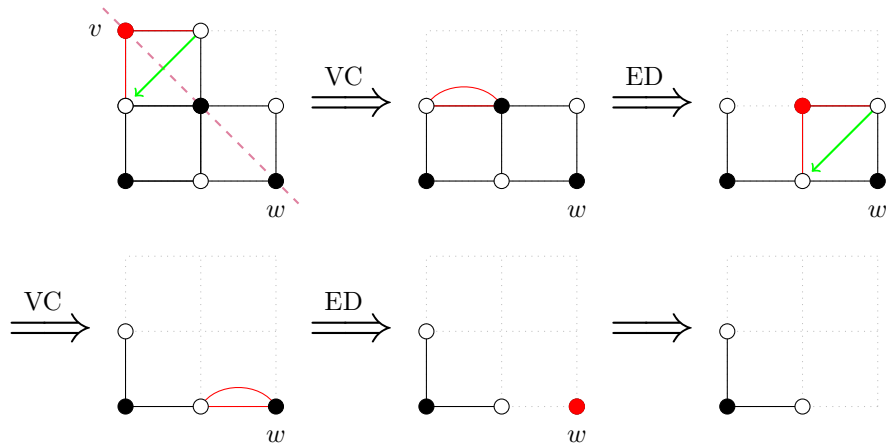
 
\begin{figure}
	\centering
	\begin{tikzpicture}
	\node at (2.3,1){$v$};
	\node at (1,-.4){$w$};
	\node (a) at (1,1){};
	\draw[purple,thick,dashed,opacity=.5] (.7,-.3) -- (2.3,1.3);
	\draw[thick,green,->] (2,0) -- (a);
	\begin{scope}[xscale=-1, shift={(-2,0)}]
	\draw[opacity=.5,gray,dotted] (0,0) grid (2,2);
	\draw (0,0) grid (2,1);
	\draw (1,0) grid (2,2);
	\draw[red] (0,0) -- (0,1) -- (1,1);
	\vertices{0}{0}{2}{1};
	\vertices{1}{1}{2}{2};
	\draw[fill,red] (0,1) circle (.1);
	\end{scope}
	
	\draw (3, 1) node[scale=2]{$\implies$};
	\draw (2.9,1.4) node{VC};
	
	\begin{scope}[shift={(4,0)}]
	\node at (1,-.4){$w$};
	\begin{scope}[xscale=-1, shift={(-2,0)}]
	\draw[opacity=.5,gray,dotted] (0,0) grid (2,2);
	\draw (1,0) grid (2,2);
	\draw[red] (1,0) to[out=150, in=-150] (1,1);
	\draw[red] (1,0) -- (1,1);
	\vertices{1}{0}{2}{2};
	\end{scope}
	\end{scope}
	
	\draw (7, 1) node[scale=2]{$\implies$};
	\draw (6.9,1.4) node{ED};
	
	\begin{scope}[shift={(8,0)}]
	\node at (1,-.4){$w$};
	\begin{scope}[xscale=-1, shift={(-2,0)}]
	\draw[opacity=.5,gray,dotted] (0,0) grid (2,2);
	\draw (1,1) grid (2,2);
	\draw (1,0) -- (2,0) -- (2,1);
	\vertices{1}{0}{2}{2};
	\draw[fill,red] (1,0) circle (.1) (2,0) circle (.1);
	\end{scope}
	\end{scope}
	
	\draw (11, 1) node[scale=2]{$\implies$};
	\draw (10.9,1.4) node{FV};
	
	\begin{scope}[shift={(12,0)}]
	\begin{scope}[xscale=-1, shift={(-2,0)}]
	\draw[opacity=.5,gray,dotted] (0,0) grid (2,2);
	\draw (1,1) grid (2,2);
	\vertices{1}{1}{2}{2};
	\end{scope}
	\end{scope}
	\end{tikzpicture}
	\caption{Diagonal contraction (from $v$ to $w$) that ends on a degree 3 vertex. In this case we end with an FV move that removes the vertex adjacent to $w$. Channels are preserved.}\label{op2}
\end{figure}
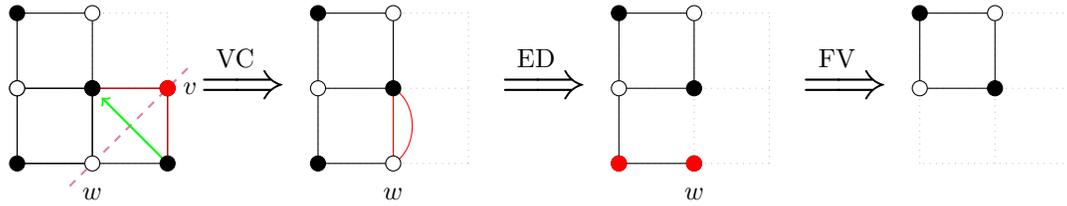

\begin{example1}\label{aztec}
	Let us apply diagonal contraction to a well-known class of graphs. The \textit{Aztec diamond} $G_n$ of rank $n$ is a diamond of side length $n$ in the square lattice. The Aztec diamonds $G_n$ for $n=1$, $2$, and $2$ are shown in Figure \ref{aztecex}.
	
	\begin{figure}[H]
		\centering
		\begin{tikzpicture}
		\draw (0,0) grid (1,1);
		\vertices{0}{0}{1}{1};
		
		\begin{scope}[shift={(4,0)}]
		\draw (0,-1) grid (1,2);
		\draw (-1,0) grid (2,1);
		\vertices{0}{-1}{1}{2};
		\vertices{-1}{0}{2}{1};
		\end{scope}
		
		\begin{scope}[shift={(10,0)}]
		\draw (0,-2) grid (1,3);
		\draw (-2,0) grid (3,1);
		\draw (-1,-1) grid (2,2);
		\vertices{0}{-2}{1}{3};
		\vertices{-2}{0}{3}{1};
		\vertices{-1}{-1}{2}{2};
		\end{scope}
		\end{tikzpicture}
		\caption{The Aztec diamonds of rank 1, 2, and 3.}
		\label{aztecex}
	\end{figure}
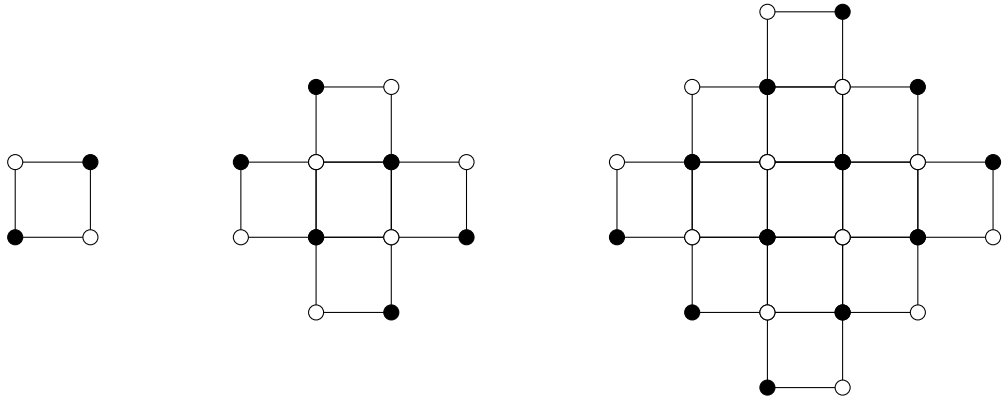

	We will show that $|\mathcal C (G_n)| = 2^{2n}$
	by induction on $n$. The rank 1 Aztec diamond has 4 channels. For the rank $n$ Aztec diamond with $n>1$, we perform the following diagonal contractions:
	
	
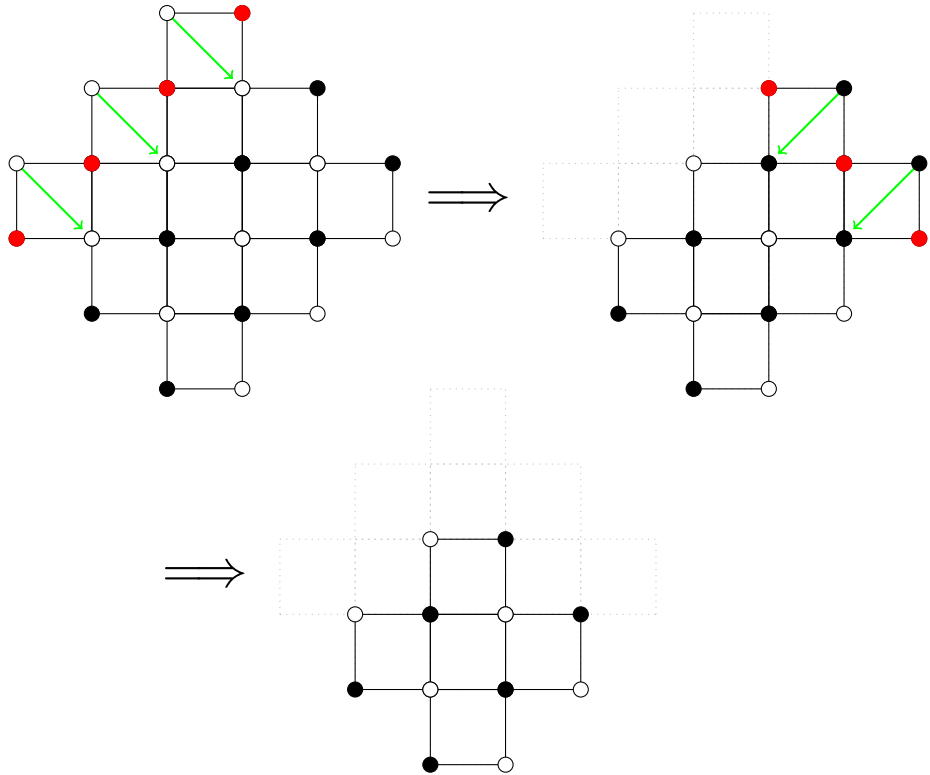
\begin{figure}[H]
		\centering
	\begin{tikzpicture}
		
		\node (a) at (-1,0){};
		\draw[green,thick, ->] (-2,1)-- (a);
		\node (a) at (0, 1){};
		\draw[green,thick,->] (-1,2) -- (a);
		\node (a) at (1, 2){};
		\draw[green,thick,->] (0,3) -- (a);
		
		\draw (0,-2) grid (1,3);
		\draw (-2,0) grid (3,1);
		\draw (-1,-1) grid (2,2);
		\vertices{0}{-2}{1}{3};
		\vertices{-2}{0}{3}{1};
		\vertices{-1}{-1}{2}{2};
		
		\draw[fill,red] (-2,0) circle (.1) (-1,1) circle (.1) (0, 2) circle (.1)
		(1, 3) circle (.1);
		
		\begin{scope}[shift={(7,0)}]
		\draw (-3, .5) node[scale=2]{$\implies$};
		
		\node (a) at (1,1) {};
		\draw[green, thick, ->] (2,2) -- (a);
		\node (a) at (2,0) {};
		\draw[green, thick, ->] (3,1) -- (a);
		
		\draw[gray,opacity=.5,dotted] (0,-2) grid (1,3);
		\draw[gray,opacity=.5,dotted] (-2,0) grid (3,1);
		\draw[gray,opacity=.5,dotted] (-1,-1) grid (2,2);
		
		\draw (0,-2) grid (1,1);
		\draw (1,0) grid (2,2);
		\draw (-1,-1) grid (2,0);
		\draw (2,0) grid (3,1);
		\vertices{0}{-2}{1}{1};
		\vertices{1}{0}{2}{2};
		\vertices{2}{0}{3}{1};
		\vertices{-1}{-1}{2}{0};
		
		\draw[fill,red] (1,2) circle (.1) (2,1) circle (.1) (3,0) circle (.1);
		\end{scope}
		
		\begin{scope}[shift={(3.5,-5)}]
		\draw (-3, .5) node[scale=2]{$\implies$};
		
		\draw[gray,opacity=.5,dotted] (0,-2) grid (1,3);
		\draw[gray,opacity=.5,dotted] (-2,0) grid (3,1);
		\draw[gray,opacity=.5,dotted] (-1,-1) grid (2,2);
		
		\draw (0,-2) grid (1,1);
		\draw (-1,-1) grid (2,0);
		\vertices{0}{-2}{1}{1};
		\vertices{-1}{-1}{2}{0};
		\end{scope}

	\end{tikzpicture}
\end{figure}
	
This produces the Aztec diamond of rank $n-1$. Since both diagonal contractions ended on a vertex of degree 2, 
$$|\mathcal C(G_n)| = 2^2|\mathcal C(G_{n-1})| .$$
The result follows by induction. Since there are $2^{2n}$ channels in $G_n$, by Theorem \ref{twos} it follows that $2^n$ divides the number of matchings of $G_n$. Indeed, it is well-known that $G_n$ has $2^{\binom{n+1}{2}}$ matchings.
\end{example1}

The fact that diagonal contraction preserves channels has some implications that have been noticed previously in the literature. For instance, Tenner \cite{Tenner2005} proves a Tiling Parity Theorem and uses it to great effect. When stated in our language, the Tiling Parity Theorem is a statement about diagonal contraction for certain diagonals that do not end on a corner. 
Here we describe a generalization of Tenner's theorem. (In her language, we cover the case $k \geq 3$, though the cases $k < 3$ follow by a similar argument.) Recall that $n\equiv_2 m$ means that $n-m$ is even.

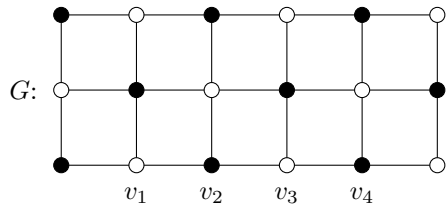
\begin{figure}[H]
	\centering
	\begin{tikzpicture}
	\node at (-.5,1) {$G$:};
	\draw (0,0) grid (5, 2);
	\vertices{0}{0}{5}{2};
	\draw (1,-.4) node{$v_1$} (2,-.4) node{$v_2$}
	(3,-.4) node{$v_3$} (4,-.4) node{$v_4$};
	\end{tikzpicture}
\end{figure}
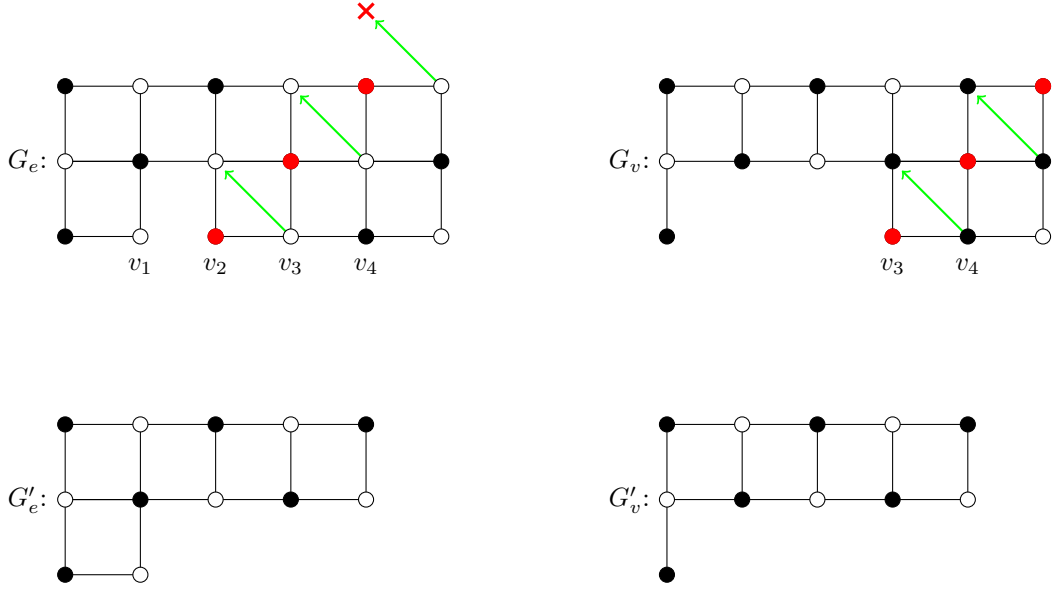
\begin{figure}[H]
	\centering
	\begin{tikzpicture}
	\node at (-.5,1) {$G_e$:};
	
	\node (a) at (2,1){};
	\draw[green,thick,->] (3,0)--(a);
	\node (a) at (3,2){};
	\draw[green,thick,->] (4,1)--(a);
	\node (a) at (4,3){};
	\draw[green,thick,->] (5,2)--(a);
	
	\draw (0,1) grid (5,2);
	\draw (0,0)grid(1,1);
	\draw (2,0) grid(5,1);
	\vertices{0}{0}{5}{2};
	
	\draw[red,fill] (2,0) circle (.1) (3,1) circle (.1) (4,2)circle(.1);
	\draw (1,-.4) node{$v_1$} (2,-.4) node{$v_2$}
	(3,-.4) node{$v_3$} (4,-.4) node{$v_4$};
	\draw[red, very thick] (4-.1, 3-.1) -- (4+.1,3+.1);
	\draw[red, very thick] (4+.1, 3-.1) -- (4-.1,3+.1);
	
	\begin{scope}[shift = {(8,0)}]
	\node at (-.5,1) {$G_v$:};
	\draw (0,1) grid (5,2);
	\draw (0,0)grid(0,1);
	\draw (3,0) grid(5,1);
	
	\node (a) at (3,1){};
	\draw[green,thick,->] (4,0)--(a);
	\node (a) at (4,2){};
	\draw[green,thick,->] (5,1)--(a);
	
	\vertices{3}{0}{5}{2};
	\vertices{0}{1}{2}{2};
	\vertices{0}{0}{0}{0};
	\draw[red,fill] (3,0) circle (.1) (4,1) circle (.1) (5,2)circle(.1);
	\draw (3,-.4) node{$v_3$} (4,-.4) node{$v_4$};
	\end{scope}
	
	
	\begin{scope}[shift={(0,-4.5)}]
	\node at (-.5,1) {$G_e'$:};
	
	\draw (0,0)grid(1,1);
	\draw (0,1) grid(4,2);
	\vertices{0}{0}{1}{2};
	\vertices{2}{1}{4}{2};
	
	\begin{scope}[shift = {(8,0)}]
	\node at (-.5,1) {$G_v'$:};
	\draw (0,1) grid (4,2);
	\draw (0,0)grid(0,1);
	
	\vertices{0}{1}{4}{2};
	\vertices{0}{0}{0}{0};
	\end{scope}
	\end{scope}
	
	\end{tikzpicture}
	\caption[]{A graph to which the Parity Theorem is applicable. By Theorem \ref{parity},
	
	\begin{minipage}{\linewidth}
	\[41 = m_G \equiv_2 m_{G_e'} + 2m_{G_v'} \equiv_2 m_{G_e'} = 11.\]
	\end{minipage}
	}
\end{figure}

\begin{theorem}[Parity Theorem]\label{parity}
	Let $G$ be a subgraph of the square lattice and let $v_1,v_2,v_3,v_4$ be collinear consecutive external vertices such that $v_2$ and $v_3$ have degree $3$. Let $G_e = G-\edge{v_1}{v_2}$ and $G_v=G-\{v_1,v_2\}$ be the edge-deleted and vertex-deleted subgraphs of $G$, so that $v_2$ defines a unique diagonal in $G_e$ and $v_3$ defines a unique diagonal in $G_v$. If both of these diagonals are contractible, then
	$$m_G \equiv_2 2^{\delta_e} m_{G_e'} + 2^{\delta_v}m_{G'_v},$$
	where $G_e'$ and $G_v'$ are the graphs resulting from contracting the diagonals in $G_e$ and $G_v$ respectively, and $\delta_i$ is $1$ or $0$ depending on if the diagonal in $G_i$ ends on a vertex of degree $2$ or of degree at least $3$, respectively.
\end{theorem}
\begin{proof}
	By Proposition \ref{matchprop},
	$$m_G = m_{G_e} + m_{G_v} .$$
	By Theorem \ref{diagthm}, channels are preserved by diagonal contractions if and only if the diagonal ends on a vertex of degree at least $3$. Specifically,
	$$|\mathcal C(G_e)| = 2^{\delta_e}|\mathcal C(G_e')|  \text{ and } |\mathcal C(G_v)|=2^{\delta_v}|\mathcal C(G_v')|.$$
	Now by Proposition \ref{algproof}, $m_{G_i}$ is even if and only if $|\mathcal C(G_i)|>1$. Thus 
	$$m_{G_e} \equiv_2 2^{\delta_e}m_{G_e'} \text{ and } m_{G_v} \equiv_2 2^{\delta_v}m_{G_v'}.$$
	The claim follows.
\end{proof}
\begin{remark}
	In fact, if $G$ has at least one matching, then 
	$$m_G \equiv_2 m_{G_e'} + m_{G_v'} .$$
	This holds because $G$ has a matching only if $G$ has an equal number of black and white vertices, whereas a diagonal contraction that ends on a degree 2 vertex removes unequal numbers of white and black vertices. Thus $\delta_i=1$ only if $m_{G_i'}=0$. 
\end{remark}

We conclude this section by showing that diagonal contraction is always possible for certain subgraphs of the square lattice.
\begin{theorem}
	Let $G$ be a subgraph of the square lattice such that every internal face of $G$ is a unit square and such that every edge bounds an internal face. Then $G$ has a contractible diagonal starting at a degree 2 vertex.
\end{theorem}
\begin{proof}
	Each internal vertex is incident only to internal faces and therefore to exactly four unit squares. Hence all internal faces and vertices have degree 4, so any diagonal in $G$ will be contractible. Any degree 2 vertex in $G$ will be a corner, since by assumption the two incident edges bound a unit square. Thus we just need to show that there is a degree 2 vertex in $G$. But this follows from Lemma \ref{cornerlemma} since every internal face and vertex have degree four.
\end{proof}
\section{Conclusion}
As we have seen, channels provide an effective lower bound on the power of two dividing a matching count. In addition, when there are no nonempty channels, they tell us that the number of matchings is odd. It would be nice to find exact powers of two more generally. This prompts a natural question.

\begin{problem}\label{prob1}
	Determine when the number of channels is the exact power of two dividing $m_G^2$. Even better, determine how many additional powers of two are carried by each channel more generally. Alternatively, provide a method for determining an upper bound on powers of two dividing $m_G$.
\end{problem}
For more on how additional powers of two are distributed in the Smith normal form of the Kasteleyn matrix (and therefore among channels), see \cite{Kuperberg2002}. Additional powers of two may be associated to a result such as Ciucu's Factorization Theorem (Proposition \ref{ciucu}). Indeed, graphs where this theorem applies tend to have additional powers of two beyond what channels would predict (for instance the Aztec Diamond in Example \ref{aztec}, cf.\ \cite{Ciucu2003}). 
One possible route for approaching Problem \ref{prob1} is to study equivalence classes of perfect matchings under the action of channels described in Theorem \ref{cycle}. Another route may arise by solving the following problem.

\begin{problem}\label{prob2}
	Find a combinatorial proof of Theorem \ref{twos}.
\end{problem}
Since Theorem \ref{twos} requires the Kasteleyn signing of $G$, such a proof would likely invoke planarity. As mentioned in Remark \ref{action}, one possible approach to this is constructing a free action of $\cb(G)$ on the set of matchings of $G$. (For the general, non-bipartite case, we would want an action of $\mc(G)$ on pairs of matchings.) Because the definition of channels involves neighborhoods of even size, searching for an action that uses properties of Eulerian circuits may yield productive results. A related problem is to construct an action of billiard nests on matchings, for graphs which admit them. This could be more tractable due to the canonical path basis for the space of nests.

The reducible graphs described in Section 6 may provide a tractable entry point to these problems. For such graphs, the problem of constructing an action of channels on matchings reduces to understanding how such an action plays with the channel-preserving moves. Since vertex contraction and forced vertex deletion both preserve matchings, this further reduces to studying the impact of doubled edge deletion on matchings.

\section*{Acknowledgments}
The first author was supported by a grant from the North Carolina State University Office of Undergraduate Research. The second author was partially supported by grants from the National Science Foundation, DMS-1700302 and CCF-1900460.




\begin{bibdiv}
\begin{biblist}

\bib{Ciucu1996}{article}{
      author={Ciucu, Mihai},
       title={Enumeration of {Perfect} {Matchings} in {Graphs} with
  {Reflective} {Symmetry}},
        date={1997-01},
        ISSN={00973165},
     journal={Journal of Combinatorial Theory, Series A},
      volume={77},
      number={1},
       pages={67\ndash 97},
         url={https://linkinghub.elsevier.com/retrieve/pii/S0097316596927259},
}

\bib{Ciucu2003}{article}{
      author={Ciucu, Mihai},
       title={Perfect matchings and perfect powers},
        date={2003-01},
        ISSN={09259899},
     journal={Journal of Algebraic Combinatorics},
      volume={17},
      number={3},
       pages={335\ndash 375},
}

\bib{Cohn1999}{article}{
      author={Cohn, Henry},
       title={2-adic behavior of numbers of domino tilings},
        date={1999-08},
        ISSN={10778926},
     journal={Electronic Journal of Combinatorics},
      volume={6},
      number={1},
       pages={\#R14},
}

\bib{Eppstein2018}{inproceedings}{
      author={Eppstein, David},
      author={Vazirani, Vijay~V.},
       title={{NC} {Algorithms} for {Computing} a {Perfect} {Matching}, the
  {Number} of {Perfect} {Matchings}, and a {Maximum} {Flow} in
  {One}-{Crossing}-{Minor}-{Free} {Graphs}},
        date={2019},
   booktitle={The 31st {ACM} on {Symposium} on {Parallelism} in {Algorithms}
  and {Architectures} - {SPAA} '19},
   publisher={ACM Press},
     address={Phoenix, AZ, USA},
       pages={23\ndash 30},
         url={http://dl.acm.org/citation.cfm?doid=3323165.3323206},
}

\bib{Fraser2017}{article}{
      author={Fraser, Chris},
      author={Lam, Thomas},
      author={Le, Ian},
       title={From dimers to webs},
        date={2019-01},
        ISSN={0002-9947, 1088-6850},
     journal={Transactions of the American Mathematical Society},
      volume={371},
      number={9},
       pages={6087\ndash 6124},
         url={http://www.ams.org/tran/2019-371-09/S0002-9947-2019-07641-9/},
}

\bib{Kasteleyn1961}{article}{
      author={Kasteleyn, P.~W.},
       title={The statistics of dimers on a lattice. i. the number of dimer
  arrangements on a quadratic lattice},
        date={1961-12},
        ISSN={00318914},
     journal={Physica},
      volume={27},
      number={12},
       pages={1209\ndash 1225},
}

\bib{Kasteleyn1963}{article}{
      author={Kasteleyn, P.~W.},
       title={Dimer statistics and phase transitions},
        date={1963-02},
        ISSN={00222488},
     journal={Journal of Mathematical Physics},
      volume={4},
      number={2},
       pages={287\ndash 293},
}

\bib{Kenyon1999}{article}{
      author={Kenyon, Richard~W.},
      author={Propp, James~G.},
      author={Wilson, David~B.},
       title={Trees and matchings.},
        date={2000},
     journal={The Electronic Journal of Combinatorics},
      volume={7},
      number={1},
       pages={\#R25},
         url={http://eudml.org/doc/120355},
}

\bib{Kuo2004}{article}{
      author={Kuo, Eric~H.},
       title={Applications of graphical condensation for enumerating matchings
  and tilings},
        date={2004},
        ISSN={03043975},
     journal={Theoretical Computer Science},
      volume={319},
      number={1-3},
       pages={29\ndash 57},
}

\bib{Kuperberg2002}{article}{
      author={Kuperberg, Greg},
       title={Kasteleyn cokernels},
        date={2002-08},
        ISSN={10778926},
     journal={Electronic Journal of Combinatorics},
      volume={9},
      number={1},
       pages={\#R29},
}

\bib{Lam2015}{article}{
      author={Lam, Thomas},
       title={Dimers, webs, and positroids},
        date={2015},
        ISSN={14697750},
     journal={Journal of the London Mathematical Society},
      volume={92},
      number={3},
       pages={633\ndash 656},
}

\bib{Lovasz1993}{book}{
      author={Lov{\'{a}}sz, L.},
       title={Combinatorial problems and exercises, second edition},
   publisher={Elsevier Science},
        date={1993},
        ISBN={9780444815040},
}

\bib{Pachter1997}{article}{
      author={Pachter, Lior},
       title={Combinatorial approaches and conjectures for 2-divisibility
  problems concerning domino tilings of polyominoes},
        date={1997},
        ISSN={10778926},
     journal={Electronic Journal of Combinatorics},
      volume={4},
       pages={{\#}R29},
}

\bib{Percus1969}{article}{
      author={Percus, Jerome~K.},
       title={One more technique for the dimer problem},
        date={1969-10},
        ISSN={00222488},
     journal={Journal of Mathematical Physics},
      volume={10},
      number={10},
       pages={1881\ndash 1884},
}

\bib{Propp1999}{article}{
      author={Propp, James},
       title={Enumeration of matchings: Problems and progress},
        date={1999-04},
     journal={New Perspectives in Geometric Combinatorics},
      volume={38},
       pages={255\ndash 291},
}

\bib{Stanley2016}{article}{
      author={Stanley, Richard~P.},
       title={Smith normal form in combinatorics},
        date={2016-11},
        ISSN={10960899},
     journal={Journal of Combinatorial Theory. Series A},
      volume={144},
       pages={476\ndash 495},
}

\bib{Tenner2005}{article}{
      author={Tenner, Bridget~Eileen},
       title={Tiling parity results and the holey square solution},
        date={2005},
     journal={Electronic Journal of Combinatorics},
      volume={11},
      number={2},
       pages={{\#}R17},
}

\bib{Tenner2009}{article}{
      author={Tenner, Bridget~Eileen},
       title={Domino tiling congruence modulo 4},
        date={2009-11},
        ISSN={0911-0119},
     journal={Graphs and Combinatorics},
      volume={25},
      number={4},
       pages={625\ndash 638},
}

\bib{Tomei2002}{article}{
      author={Tomei, Carlos},
      author={Vieira, Tania},
       title={The kernel of the adjacency matrix of a rectangular mesh},
        date={2002-08},
        ISSN={01795376},
     journal={Discrete and Computational Geometry},
      volume={28},
      number={3},
       pages={411\ndash 425},
}

\end{biblist}
\end{bibdiv}
\end{document}